\documentclass{amsart}%
\usepackage{amssymb}
\usepackage{amsfonts}
\usepackage{amsmath}
\usepackage{graphicx}%
\newtheorem{theorem}{Theorem}
\theoremstyle{plain}
\newtheorem{theorema}{Theorem A}
\newtheorem{theoremb}{Theorem B}
\newtheorem{theoremc}{Theorem C}

\newtheorem{lemma}{Lemma}

\newtheorem{proposition}{Proposition}
\newtheorem{remark}{Remark}

\numberwithin{equation}{section}
\begin{document}
\title{Onto Interpolating Sequences for the Dirichlet Space}
\author{Nicola Arcozzi}
\address{Dipartimento do Matematica, Universita di Bologna, 40127 Bologna, ITALY}
\author{Richard Rochberg}
\address{Campus Box 1146, Washington University, 1 Brookings Drive, St. Louis, MO 63130}
\author{Eric Sawyer}
\address{Department of Mathematics \& Statistics, McMaster University Hamilton,
Ontairo, L8S4K1, CANADA}
\thanks{(Arcozzi) Partially supported by the COFIN project Analisi Armonica, funded by
the Italian Minister for Research.}
\thanks{(Rochberg) This material is based upon work supported by the National Science
Foundation under Grant No. 0400962.}
\thanks{(Sawyer) This material based upon work supported by the National Science and
Engineering Council of Canada}

\begin{abstract}
We describe two new classes of onto interpolating sequences for the Dirichlet
space, in particular resolving a question of Bishop. We also give a complete
description of the analogous sequences for a discrete model of the Dirichlet space.

\end{abstract}
\maketitle
\tableofcontents

\section{Overview}

We begin with an informal overview; precise definitions and statements will be
given later.

\subsection{Interpolating Sequences for the Hardy Space}

The study of interpolating sequences for the Dirichlet space evolved from the
study of interpolating sequences for the Hardy space and we begin by recalling
that work.

The Hardy space, $H^{2},$ is a Hilbert space of holomorphic functions on the
unit disk. It is normed by, for $f=\sum_{n\leq0}a_{n}z^{n},$ $\left\Vert
f\right\Vert _{H^{2}}^{2}=\sum\left\vert a_{n}\right\vert ^{2}.$ It has
reproducing kernels $k_{z}^{H^{2}}=(1-\bar{z}w)^{-1}$. The kernel has norm
$\left\Vert k_{z}^{H^{2}}\right\Vert _{H^{2}}=(1-\left\vert z\right\vert
^{2})^{-1/2}$ and we denote the normalized kernels by $\hat{k}_{z}^{H^{2}};$
$\hat{k}_{z}^{H^{2}}=(1-\left\vert z\right\vert ^{2})^{1/2}k_{z}^{H^{2}}.$

Given a sequence $Z=\left\{  z_{i}\right\}  $ contained in the disk we define
a weighted restriction map $R_{H^{2}}=R_{H^{2},Z}$ mapping functions on the
disk to functions on $Z$ by $\left(  R_{H^{2}}f\right)  (z_{i})=\left\Vert
k_{z_{i}}^{H^{2}}\right\Vert _{H^{2}}^{-1}f(z_{i}).$ Straightforward Hilbert
space considerations insure that $R_{H^{2}}$ maps $H^{2}$ into $\ell^{\infty
}(Z).$ If, in fact, $R_{H^{2}}$ maps $H^{2}$ into and onto $\ell^{2}(Z)$ then
we say that $Z$ is an \textit{interpolating sequence} for the Hardy space.

If $R_{H^{2}}$ maps onto $\ell^{2}(Z)$ then norm control is possible; there is
a $C>0$ so that for each $\gamma\in\ell^{2}(Z)$ there is a $\Gamma\in H^{2}$
with $R_{H^{2}}\Gamma=\gamma$ and $\left\Vert \Gamma\right\Vert _{H^{2}}\leq
C\left\Vert \gamma\right\Vert .$ To see this note that, letting $V_{Z}$ be the
closed subspace of functions in $H^{2}$ that vanish on $Z,$ there is a well
defined linear map $\Lambda$ of $\ell^{2}(Z)$ into $V_{Z}^{\perp}$ such that
$R_{H^{2}}\Lambda$ is the identity on $\ell^{2}(Z).$ The closed graph theorem
insures that $\Lambda$ is continuous and thus, given $\gamma\in\ell^{2}(Z),$
the choice $\Gamma=\Lambda\gamma$ meets the requirements. In particular, if
$R_{H^{2}}$ is onto then there is a $C>0$ so that for any $i,j,i\neq j$ we can
find $f_{ij}\in H^{2}$ with
\begin{equation}
R_{H^{2}}f_{ij}(z_{i})=1,\text{ }R_{H^{2}}f_{ij}(z_{j})=0 \label{size}%
\end{equation}
and $\left\Vert f_{ij}\right\Vert _{H^{2}}\leq C.$ On the other hand, given
$z_{i},z_{j}\in\mathbb{D},$ any $f_{ij}\in H^{2}$ which satisfies (\ref{size})
then it must satisfy%
\[
\left\Vert f_{ij}\right\Vert \geq\left(  1-\left\vert \left\langle \hat
{k}_{z_{i}}^{H^{2}},\hat{k}_{z_{j}}^{H^{2}}\right\rangle \right\vert
^{2}\right)  ^{-1/2}.
\]
If $Z$ is an interpolating sequence then we can combine this with the previous
estimate for $\left\Vert f_{ij}\right\Vert $ and find that the points of $Z$
satisfy a separation condition: $\exists\varepsilon>0,\forall i,j,i\neq j$%
\begin{equation}
\left\vert \left\langle \hat{k}_{z_{i}}^{H^{2}},\hat{k}_{z_{j}}^{H^{2}%
}\right\rangle \right\vert \leq1-\varepsilon. \tag{Sep(Hardy)}\label{seph}%
\end{equation}
An equivalent condition, written using $\beta,$ the hyperbolic metric on the
disk, is $\exists\varepsilon^{\prime}>0,\forall i,j,i\neq j$%
\begin{equation}
\beta(z_{i},z_{j})\geq\varepsilon^{\prime}. \tag{Sep'(Hardy)}\label{seph'}%
\end{equation}

The interpolating sequence $Z$ must satisfy (\ref{seph}) because $R_{H^{2}}$
is onto. The requirement that $R_{H^{2}}$ be into, that is, that $R_{H^{2}}$
be bounded, gives a different requirement. Associate to the sequence $Z$ the
measure $\lambda_{Z}$ given by
\begin{equation}
\lambda_{Z}=\sum_{j=1}^{\infty}\left\Vert \tilde{k}_{z}\right\Vert _{H^{2}%
}^{-2}\delta_{z_{j}}. \label{measure}%
\end{equation}
The condition that $R_{H^{2}}$ be bounded, that there is a $C>0$ so that for
all $f\in H^{2}$ $\left\Vert R_{H^{2}}f\right\Vert \leq\sqrt{C}\left\Vert
f\right\Vert _{H^{2}},$ is equivalent to the condition that $\lambda_{Z}$ be a
\emph{Carleson measure }(for the Hardy space). That is, there is a $C>0$ so
that for $\lambda=\lambda_{Z}$ and all $f\in H^{2}$%
\begin{equation}
\int\left\vert f\right\vert ^{2}d\lambda\leq C\left\Vert f\right\Vert _{H^{2}%
}^{2}. \tag{Car(Hardy)}\label{carh}%
\end{equation}
If this condition holds for some $\lambda$ then \textit{a fortiori} it holds
when $f$ is a reproducing kernel. This implies the \textit{simple condition:
}$\exists C>0,$ $\forall z\in\mathbb{D}$
\begin{equation}
\lambda(T(z))\leq C\left\vert I_{z}\right\vert . \tag{CarSimp(Hardy)}%
\label{carsimph}%
\end{equation}
Here $I_{z}$ is the boundary interval with center $z/\left\vert z\right\vert $
and length $2\pi(1-\left\vert z\right\vert ),$ $\left\vert I_{z}\right\vert $
is its length, and $T(z)$ is the \textit{tent over }$I_{z},$ the convex hull
of $z$ and $I_{z}.$

In fact, these conditions characterize interpolating sequences. Combining the
results of Carleson in 1958 \cite{C} and Shapiro and Shields in 1961 \cite{SS}
we have

\begin{theorem}
\label{hardy}Given $Z,$ the following conditions are equivalent:

\begin{enumerate}
\item $R_{H^{2},Z}$ maps $H^{2}$ into and onto $\ell^{2}\left(  Z\right)  $,
that is, $Z$ is an interpolating sequence for the Hardy space.

\item $Z$ satisfies (\ref{seph}) and $\lambda_{Z}$ satisfies (\ref{carh}),

\item $Z$ satisfies (\ref{seph}) and $\lambda_{Z}$ satisfies (\ref{carsimph}),

\item $R_{H^{2},Z}$ maps $H^{2}$ onto $\ell^{2}\left(  \lambda_{Z}\right)  .$
\end{enumerate}
\end{theorem}

The equivalence of the first two statements is the traditional description of
interpolating sequences for the Hardy space. The second and third are
equivalent because for any positive measure $\lambda$ the conditions
(\ref{carh}) and (\ref{carsimph}) are equivalent. The first condition
certainly implies the last; the converse of that implication is a subtle
consequence of the details of the analysis.

\subsection{Interpolating Sequences for the Dirichlet Space}

The Dirichlet space, $B_{2},$ is a Hilbert space of holomorphic functions on
the unit disk. It is normed by, for $f=\sum_{n\leq0}a_{n}z^{n},$ $\left\Vert
f\right\Vert ^{2}=\left\vert a_{0}\right\vert ^{2}+\sum_{n>0}n\left\vert
a_{n}\right\vert ^{2}.$ It has reproducing kernels $k_{z}=-\bar{z}w\log
(1-\bar{z}w).$ Note that $\left\Vert k_{z}\right\Vert ^{2}=-\left\vert
z\right\vert ^{2}\log(1-\left\vert z\right\vert ^{2})$ and that for $z$ near
the boundary, the only case of interest for us,
\[
\left\Vert k_{z}\right\Vert ^{2}\sim-\log(1-\left\vert z\right\vert ^{2}).
\]
We denote the normalized kernels by $\hat{k}_{z};$ $\hat{k}_{z}=\left\Vert
k_{z}\right\Vert ^{-1}k_{z}.$

Given a sequence $Z=\left\{  z_{i}\right\}  \subset\mathbb{D}$ we now define a
weighted restriction map $R=R_{z}$ by $\left(  Rf\right)  (z_{i})=\left\Vert
k_{z_{i}}\right\Vert ^{-1}f(z_{i}).$ As before it is automatic that $R$ maps
$B_{2}$ into $\ell^{\infty}(Z).$ If, in fact, $R$ maps $B_{2}$ into and onto
$\ell^{2}(Z)$ then we say that $Z$ is an \textit{interpolating sequence} for
the Dirichlet space.

As in the Hardy space case, there are two natural necessary conditions for $R
$ to be an interpolating sequence. The fact that $R$ maps onto $\ell^{2}(Z)$
insures that the points of $Z$ satisfy a separation condition: $\exists
\varepsilon>0,\forall i,j,i\neq j$%
\begin{equation}
\left\vert \left\langle \hat{k}_{z_{i}},\hat{k}_{z_{j}}\right\rangle
\right\vert \leq1-\varepsilon. \tag{Sep}\label{sep}%
\end{equation}
This condition can also be given an equivalent reformulation using $\beta:$
$\exists C>0,\forall i,j,i\neq j$%
\begin{equation}
\beta(z_{i},z_{j})\geq C(1+\beta(0,z_{j})). \tag{Sep'}\label{sep'}%
\end{equation}

The map $R$ is bounded if and only if the measure $\mu_{Z}$ defined by
\begin{equation}
\mu_{Z}=\sum_{j=1}^{\infty}\left\Vert k_{z}\right\Vert ^{-2}\delta_{z_{j}}.
\tag{Associated measure}\label{mu}%
\end{equation}
is a \emph{Carleson measure, }but now a Carleson measure for the Dirichlet
space. That is, there is a $C>0$ so that for $\mu=\mu_{Z}$ and all $f\in
B_{2}$%
\begin{equation}
\int\left\vert f\right\vert ^{2}d\mu\leq C\left\Vert f\right\Vert ^{2}.
\tag{Car}\label{car}%
\end{equation}
With the choice $f=k_{z}$ this estimate implies the Dirichlet space version of
the \textit{simple condition: }$\exists C>0,$ $\forall z\in\mathbb{D}$,
$\left\vert z\right\vert \sim1$%
\begin{equation}
\mu(T(z))\leq C(-\log(1-\left\vert z\right\vert ^{2}))^{-1}. \tag{CarSimp}%
\label{carsimp}%
\end{equation}

In unpublished work in the early 1990's Bishop \cite{Bi} and, independently,
Marshall-Sundberg \cite{MS} characterized the interpolating sequences for the
Dirichlet space. The first published proof was given by B\"{o}e.\cite{Bo} in
2002 using different techniques.

\begin{theorem}
\label{dirichlet}Given $Z,$ the following conditions are equivalent:

\begin{enumerate}
\item $R_{Z}$ maps $B_{2}$ into and onto $\ell^{2}\left(  Z\right)  $, that
is, $Z$ is an interpolating sequence for the Dirichlet space.

\item $Z$ satisfies (\ref{sep}) and $\mu_{Z}$ satisfies (\ref{car}),
\end{enumerate}
\end{theorem}

This is a very satisfying analogy with the equivalence of the first two
statements in Theorem \ref{hardy}. The first statements of the theorems are
similar by design, (\ref{seph'}) and (\ref{sep'}) differ in detail but are
similar in spirit, and (\ref{carh}) and (\ref{car}) are formally the same.
However there are also fundamental differences between the two theorems. The
subtle condition (\ref{carh}) is equivalent to the more straightforward
geometric condition (\ref{carsimph}). On the other hand the simple condition
(\ref{carsimp}), while necessary for (\ref{car}) is not equivalent to it.

The fact that condition (5) in Theorem \ref{hardy} is equivalent to the others
is one of the deeper parts of that theorem and the analogous statement fails
for the Dirichlet space. Bishop had noted that there are sequences for which
the restriction map is onto, i.e. $\ell^{2}\left(  \mu_{Z}\right)  \subset
R\left(  B_{2}\right)  $, but the restriction map is not bounded; and hence
those sequences are not interpolating sequences. We call sequences $Z$ for
which the restriction map is onto but not necessarily bounded \emph{onto
interpolating sequences} (for the Dirichlet space). The study of such
sequences is the theme of this paper.

\subsection{Onto Interpolating Sequences for the Dirichlet Space}

If $Z$ is an onto interpolating sequence for the Dirichlet space then the
arguments that led to (\ref{sep}) and (\ref{sep'}) still apply and $Z$ must be
satisfy those conditions. On the other hand, condition (\ref{car}) was a
reformulation of the requirement that $R_{Z}$ be bounded. That requirement is
not imposed in this case and (\ref{car}) is not necessary; in fact it plays no
role in what follows.

However there are three conditions on $\mu_{Z}$ that we will be working with.
The first condition is that $\mu_{Z}$ be finite, $\left\Vert \mu
_{Z}\right\Vert <\infty.$ The second condition is that $\mu_{Z}$ satisfy the
simple condition (\ref{carsimp}).

The third condition is a weaker variation on the simple condition. To describe
it we introduce Bergman trees which will be a basic tool in our analysis. We
describe them now informally, the detailed description is in \cite{ArRoSa},
\cite{ArRoSa2}, or \cite{Sa}. A Bergman tree is a set $\mathcal{T}=\left\{
\alpha_{i}\right\}  \subset\mathbb{D}$ for which there is a positive lower
bound on the hyperbolic distances between distinct points and so that for some
constant $C$ the union of hyperbolic balls, $\bigcup_{i}B(\alpha_{i},C),$
cover $\mathbb{D}$. We regard $\left\{  \alpha_{i}\right\}  $ as the vertices
of a rooted tree with $o$, the vertex nearest the origin, as root. Each point
$\alpha_{i}$ of $\mathcal{T}$, except $o$, is connected by an edge to its
predecessor $\alpha_{i}^{-}$, a nearby point closer to the origin. Being a
rooted tree, $\mathcal{T}$ has a partial ordering; $\alpha\leq\beta$ if
$\alpha$ is on the geodesic connecting $\beta$ to $o$. For $\alpha
\in\mathcal{T}$ let $d(\alpha)$ be the number of vertices on the tree geodesic
connecting $\alpha$ to $o;$ in particular, for any $\alpha,$ $d(\alpha)\geq1.$
We will assume that each $\alpha\in\mathcal{T}$ is the predecessor of exactly
two other points of $\mathcal{T}$, the successors, $\alpha_{\pm}.$ This
assumption is a notational convenience; our trees automatically have an upper
bound on their branching number and all our discussions extend to that case by
just adding notation.

If $Z\ $is separated then there is a positive lower bound on the hyperbolic
distance between distinct points of $Z$ and, given this, it is easy to see
that we can construct a Bergman tree $\mathcal{T}$ for the disk which contains
$Z$ among its nodes. So, without loss of generality, we may assume that
$Z\subset\mathcal{T}$. When the points of $Z$ are regarded as elements of
$\mathcal{T}$ we will often denote them with lower case Greek letters.

For any $\alpha\in\mathcal{T}$ we denote the \textit{shadow of }$\alpha$ by
$S(\alpha):$%
\begin{equation}
S(\alpha)=\left\{  \beta\in\mathcal{T}:\beta\geq\alpha\right\}  . \tag{Shadow}%
\end{equation}
It is the tree analog of the tent $T(\alpha).$ If $\mu_{Z}$ satisfies
(\ref{carsimp}) then, regarded as a measure on $\mathcal{T},$ the measure
satisfies $\exists C>0$ $\forall\alpha\in\mathcal{T}$
\begin{equation}
\mu(S(\alpha))\leq Cd\left(  \alpha\right)  ^{-1}. \tag{TreeSimp}%
\label{treesimple}%
\end{equation}
In this case will say that the measure or the sequence satisfy the
\textit{simple condition}. We say they satisfy the \emph{weak simple condition
}if:\emph{\ }$\exists C>0,$ $\forall\alpha\in\mathcal{T}$
\begin{equation}
\sum_{\substack{\beta\in Z,\text{ }\beta\geq\alpha\\\mu\left(  \gamma\right)
=0\text{ for }\alpha<\gamma<\beta}}\mu\left(  \beta\right)  \leq Cd\left(
\alpha\right)  ^{-1}. \tag{WeakSimp}\label{weaksimple}%
\end{equation}
That is, in this case, for each $\alpha$ we now only consider points of $Z\cap
S(\alpha)$ which have an unobstructed view of $\alpha.$

The following was shown by Bishop \cite{Bi} and, as was noted in \cite{Bo}, it
is also a consequence of the proof in \cite{Bo}.

\begin{theorem}
If the sequence $Z\subset\mathbb{D}$ is separated, the measure $\mu_{Z}$ is
finite, and the measure satisfies the simple condition; then $Z$ is an onto
interpolating sequence.
\end{theorem}

We prove two theorems which extend this result.

\begin{theorema}
If the sequence $Z\subset\mathbb{D}$ is separated, the measure $\mu_{Z}$ is
finite, and the measure satisfies the weak simple condition; then $Z$ is an
onto interpolating sequence.
\end{theorema}

Our next theorem removes the hypothesis that the measure is finite thus, in
particular, answering a question of Bishop who had asked if an onto
interpolating sequence had to be associated with a finite measure.

Our proof requires additional geometric structure for $Z$, that it be
\textit{tree-like. }We say that a sequence $Z$ which is separated and
satisfies the weak simple condition is tree-like if whenever $\alpha,\beta\in
Z$ and $\alpha$ is in a certain expanded version of $S(\beta)$ then, in fact
$\alpha\in S(\beta).$ Specifically, with $C$ the constant from (\ref{sep'})
there is a $\beta\in\left(  1-C/2,1\right)  \ $so that
\begin{align}
&  \text{if }z_{j},z_{k}\in Z,\text{ }\left\vert z_{j}\right\vert
\geq\left\vert z_{k}\right\vert ,\text{ and }\left\vert z_{j}-\left\vert
z_{k}\right\vert ^{-1}z_{k}\right\vert \leq\left(  1-\left\vert z_{k}%
\right\vert ^{2}\right)  ^{\beta}\qquad\text{ } \tag{Tree-like}\label{special}%
\\
&  \qquad\qquad\qquad\qquad\text{then \ }z_{j}\in T\left(  z_{k}\right)
.\nonumber
\end{align}

\begin{theoremb}
If the sequence $Z\subset\mathbb{D}$ is separated, the measure $\mu_{Z}$
satisfies the weak simple condition; and $Z$ is tree-like, then $Z$ is an onto
interpolating sequence.
\end{theoremb}

The proofs of these theorems uses an elaboration of the constructive
techniques of \cite{Bo} and \cite{ArRoSa2}.

\subsection{The B\"{o}e space}

The building blocks for our constructions are a type of function introduced by
B\"{o}e (Lemma \ref{analyticcon} below), and the interpolating functions we
construct all lie in a closed subspace of $B_{2}$ spanned by those functions.
We call that span the B\"{o}e space. We will, in fact, prove refinements of
Theorems A and B which also include the result that if onto interpolation is
possible for a separated sequence $Z$ \textit{using only functions from the
B\"{o}e space}$\mathit{\ }$then $Z$ must be separated and satisfy the weak
simple condition.

\subsection{Tree Interpolation and Tree Capacities}

Recall that the Dirichlet space $B_{2}\left(  \mathcal{T}\right)  $ of the
tree $\mathcal{T}$ consists of all complex-valued functions $f$ on
$\mathcal{T}$ for which the norm
\begin{equation}
\left\Vert f\right\Vert _{B_{2}\left(  \mathcal{T}\right)  }=\left\{
\left\vert f\left(  o\right)  \right\vert ^{2}+\sum\nolimits_{\beta
\in\mathcal{T}}\left\vert \bigtriangleup f\left(  \beta\right)  \right\vert
^{2}\right\}  ^{1/2}<\infty, \tag{$B_2\left(     \mathcal{T}\right)     $}%
\end{equation}
Here $\bigtriangleup f\left(  \beta\right)  =f(\beta)-f(\beta^{-})$. The
analog of Theorem \ref{dirichlet} for $B_{2}\left(  \mathcal{T}\right)  $ is
Theorem 26 in \cite{ArRoSa}. Here we complement that and characterize the onto
interpolating sequences for $B_{2}\left(  \mathcal{T}\right)  ,$ that is,
sequences for which the restriction map taking functions on $\mathcal{T}$ to
functions on $Z$ takes $B_{2}\left(  \mathcal{T}\right)  $ onto $\ell
^{2}\left(  \mu\right)  $.

These sequences are analogs of onto interpolation sequences for $B_{2}.$
Furthermore, it was shown in \cite[Sec. 7.2]{ArRoSa2} that the restriction map
$f\rightarrow\left\{  f\left(  \alpha\right)  \right\}  _{\alpha\in
\mathcal{T}}$ is bounded from $B_{2}$ to $B_{2}\left(  \mathcal{T}\right)  .$
As a consequence%
\begin{align}
&  \text{Every onto interpolating sequence for }B_{2} \tag{Restriction}%
\label{every}\\
&  \text{is an onto interpolating sequence for }B_{2}\left(  \mathcal{T}%
\right)  .\text{\ }\nonumber
\end{align}

We have a characterization of the onto interpolating sequence for
$B_{2}\left(  \mathcal{T}\right)  .$\ We say a subset $Z$ of $\mathcal{T}$
satisfies the \emph{tree capacity condition }if\emph{\ }$\exists
C>0,\forall\alpha\in Z$%
\begin{equation}
\inf\left\{  \sum\limits_{\beta\in\mathcal{T}}\left\vert \bigtriangleup
f\left(  \beta\right)  \right\vert ^{2}:f\left(  \alpha\right)  =1,\text{
}f\left(  \gamma\right)  =0\text{ }\forall\gamma\in Z\setminus\left\{
\alpha\right\}  \right\}  \leq\frac{C}{d(\alpha)}. \tag{TreeCap}%
\label{treecap}%
\end{equation}

\begin{theoremc}
The sequence $Z\subset\mathcal{T}$ is an onto interpolating sequence for
$B_{2}\left(  \mathcal{T}\right)  $ if and only if it satisfies the tree
capacity condition.
\end{theoremc}

This result is analogous to the capacitary characterization of onto
interpolating sequences for the Dirichlet space $B_{2}$ given by Bishop in
\cite{Bi}, namely that $Z$ is onto interpolating for the Dirichlet space
$B_{2}$ if and only if for any for $z,w\in Z$. there is an $F_{w}\in
B_{2}\mathbb{\ }$with $F_{w}\left(  z\right)  =\delta_{z,w}$ and $\left\Vert
F_{w}\right\Vert _{B_{2}}^{2}\leq C(-\log(1-\left\vert w\right\vert
^{2}))^{-1}.$This condition is called \emph{weak interpolation} by Schuster
and Seip \cite{ScS}. Thus, for both $B_{2}\left(  \mathcal{T}\right)  $ and
$B_{2}$ a weak interpolating sequence is an onto interpolating sequence.

\subsection{Relations Between the Conditions}

We construct a sequence covered by Theorem A but not by Theorem 3 and one
covered by Theorem B but not Theorem A. That second example answers Bishop's
question by giving an example of an onto interpolating $Z$ with $\left\Vert
\mu_{Z}\right\Vert =\infty.$

In analysis on the tree we will show that the tree separation condition
(\ref{treesep}) and the weak simple condition (\ref{weaksimple}) imply the
tree capacity condition (\ref{treecap}). On the other hand, using a recursive
scheme for computing tree capacities which is developed while proving Theorem
C, we construct an example of an onto interpolating sequence for $B_{2}\left(
\mathcal{T}\right)  $ that fails not only the weak simple condition, but also
fails to be contained in any separated sequence satisfying the weak simple condition.

Taking into account (\ref{every}) we see the tree capacity condition is
\emph{necessary} for $Z$ to be onto interpolating for the classical Dirichlet
space $B_{2}$. This gives a geometric condition stronger the separation
condition (\ref{sep}) which any onto interpolating sequence must satisfy.
Bishop \cite{Bi} gave a similar necessary condition involving logarithmic
capacity on the circle.

We also construct an example of a separated sequence with $\left\Vert \mu
_{Z}\right\Vert <\infty$ for which $Z$ fails the tree capacity condition
(\ref{treecap}) and hence, by Theorem C and (\ref{every}), is not an onto
interpolating for $B_{2}\left(  \mathcal{T}\right)  $ or for $B_{2}.$ Thus
separation and \emph{finite} measure alone are not enough for onto interpolation.

\subsection{Contents}

Theorems A, and B are proved in Sections 2 and 3 but the proof of the
necessity of the conditions under additional hypotheses is postponed to
Section 7. Theorem C is proved in Sections 5. Sections 4 and 6 contain
examples and further discussion of the relationship between the various conditions.

\section{Theorem A\label{interpolation}}

In fact, the hypotheses (\ref{sep}), (\ref{weaksimple}) and $\left\Vert
\mu\right\Vert <\infty$ yield more than Theorem A. The more general theorem
includes a converse, that under certain extra conditions one can conclude that
$Z$ satisfies both (\ref{sep}) and (\ref{weaksimple}). In this section we
state the general theorem and prove the half that implies Theorem A. The other
half is proved in Section \ref{neces}.

Suppose $Z$ is given and fixed. For $w\in Z$ we denote by $\varphi_{w}$ the
function introduced by B\"{o}e in \cite{Bo} in his work on interpolation. By
construction $\varphi_{w}$ is a function in $B_{2}$ which is essentially 1 on
the tent $T(w)$ and small away from that region. The details of the
construction and properties are recalled in Lemma \ref{analyticcon} below.
Actually there are various choices in Lemmas \ref{Mars} and \ref{analyticcon}
below. We assume that allowable choices have been made once and for all. Also,
we further require that the chosen parameters satisfy
\begin{align}
\beta &  <\alpha<\frac{2\beta\eta}{(\eta+1)}\label{parameters}\\
s  &  >\frac{(\alpha-\rho)}{(\rho-\beta)}\nonumber
\end{align}
as we need that in our proof of the necessity of (\ref{weaksimple}) for a
certain type of interpolation, Proposition \ref{Riesz} below.

We define the \emph{B\"{o}e space,} $B_{2,Z}$, to be the closed linear span in
$B_{2}$ of the functions $\left\{  \varphi_{w}\right\}  _{w\in Z}.$ It follows
from (\ref{bounds}) in Lemma \ref{inductivecon} and Proposition \ref{Riesz}
below that for an appropriate cofinite subset $\left\{  \zeta_{j}\right\}
_{j=1}^{\infty}$ of $Z$ we have%
\begin{align}
B_{2,Z}  &  =\left\{  \varphi=\sum_{j=1}^{\infty}a_{j}\varphi_{\zeta_{j}}%
:\sum_{j=1}^{\infty}\left\vert a_{j}\right\vert ^{2}\mu\left(  \zeta
_{j}\right)  <\infty\right\}  ,\label{atomicBoe}\\
\left\Vert \varphi\right\Vert _{B_{2,Z}}  &  \approx\left\Vert \left\{
a_{j}\right\}  _{j=1}^{\infty}\right\Vert _{\ell^{2}\left(  \mu\right)
}.\nonumber
\end{align}
We will say $Z$ is an onto interpolating sequence for $B_{2,Z}$ if it is an
onto interpolating sequence for the Dirichlet space $B_{2}$ and if, further,
the interpolating functions can all be selected from $B_{2,Z}$. In the case
that $\mu_{Z}$ is finite, the next theorem completely characterizes onto
interpolation for $B_{2,Z}$ in terms of separation and the weak simple
condition. Theorem A in an immediate corollary.

\begin{theorem}
\label{Boeonto}Let $Z\subset\mathbb{D}$ and suppose $\left\Vert \mu
_{Z}\right\Vert <\infty$. Then $Z$ is an onto interpolating sequence for the
B\"{o}e space $B_{2,Z}$ if and only if both the separation condition
(\ref{sep}) and the weak simple condition (\ref{weaksimple}) hold.
\end{theorem}

We will need the following lemma from \cite{MS}, (see also \cite{ArRoSa2}).
Let $C$ be the constant in (\ref{sep}). For $w\in\mathbb{D}$ and
$1-C/2<\beta<1$, define%
\[
V_{w}=V_{w}^{\beta}=\left\{  z\in\mathbb{D}:\left\vert z-\left\vert
w\right\vert ^{-1}w\right\vert \leq(1-\left\vert w\right\vert ^{2})^{\beta
}\right\}  .
\]

\begin{lemma}
\label{Mars}Suppose the separation condition in (\ref{sep}) holds. Then for
every $\beta$ satisfying $1-C/2<\beta<1$ there is $\eta>\beta\eta>1$ such that
if $V_{z_{i}}^{\beta}\cap V_{z_{j}}^{\beta}\neq\phi$ and $\left\vert
z_{j}\right\vert \geq\left\vert z_{i}\right\vert $, then $z_{i}\notin
V_{z_{j}}^{\beta}$ and
\begin{equation}
\left(  1-\left\vert z_{j}\right\vert \right)  \leq\left(  1-\left\vert
z_{i}\right\vert \right)  ^{\eta}. \label{ji}%
\end{equation}

\end{lemma}

We have the following useful consequence of Lemma \ref{Mars}. If $\sigma>0$
and $\mu$ satisfies (\ref{weaksimple}), then
\begin{equation}
\sum_{z_{j}\geq z_{k}}\left(  1-\left\vert z_{j}\right\vert \right)  ^{\sigma
}\leq C_{\sigma}\left(  1-\left\vert z_{k}\right\vert \right)  ^{\sigma}.
\label{sigmasum}%
\end{equation}
Indeed, if $\mathcal{G}_{1}\left(  z_{k}\right)  =\left\{  \alpha_{m}%
^{1}\right\}  $ consists of the minimal elements in $\left[  S\left(
z_{k}\right)  \setminus\left\{  z_{k}\right\}  \right]  \cap Z$,
$\mathcal{G}_{2}\left(  z_{k}\right)  =\cup_{m}\mathcal{G}_{1}\left(
\alpha_{m}^{1}\right)  $, etc., we have using (\ref{ji}) and (\ref{weaksimple}%
),
\begin{align*}
\sum_{z_{j}\geq z_{k}}\left(  1-\left\vert z_{j}\right\vert \right)
^{\sigma}  &  =\left(  1-\left\vert z_{k}\right\vert \right)  ^{\sigma}%
+\sum_{\ell=1}^{\infty}\sum_{\beta\in\mathcal{G}_{\ell-1}\left(  z_{k}\right)
}\sum_{\alpha\in\mathcal{G}_{1}\left(  \beta\right)  }\left(  1-\left\vert
\alpha\right\vert \right)  ^{\sigma}\\
&  \leq\left(  1-\left\vert z_{k}\right\vert \right)  ^{\sigma}+C_{\delta}%
\sum_{\ell=1}^{\infty}\sum_{\beta\in\mathcal{G}_{\ell-1}\left(  z_{k}\right)
}\sum_{\alpha\in\mathcal{G}_{1}\left(  \beta\right)  }\left(  1-\left\vert
\alpha\right\vert \right)  ^{\sigma-\delta}\left(  \log\frac{1}{1-\left\vert
\alpha\right\vert }\right)  ^{-1}\\
&  \leq\left(  1-\left\vert z_{k}\right\vert \right)  ^{\sigma}+C_{\delta}%
\sum_{\ell=1}^{\infty}\sum_{\beta\in\mathcal{G}_{\ell-1}\left(  z_{k}\right)
}\left(  1-\left\vert \beta\right\vert \right)  ^{\left(  \sigma
-\delta\right)  \eta}C\left(  \log\frac{1}{1-\left\vert \beta\right\vert
}\right)  ^{-1}\\
&  \leq\left(  1-\left\vert z_{k}\right\vert \right)  ^{\sigma}+C_{\delta}%
\sum_{z_{j}>z_{k}}\left(  1-\left\vert z_{j}\right\vert \right)  ^{\left(
\sigma-\delta\right)  \eta}.
\end{align*}
Now we can choose $\delta>0$ so small that $\left(  \sigma-\delta\right)
\eta-\sigma=\theta>0$, and $R$ such that $C_{\delta}\left(  1-R\right)
^{\theta}=\frac{1}{2}$, so that for $\left\vert z_{k}\right\vert \geq R$ we
have
\[
C_{\delta}\sum_{z_{j}>z_{k}}\left(  1-\left\vert z_{j}\right\vert \right)
^{\left(  \sigma-\delta\right)  \eta}\leq\left\{  C_{\delta}\sup_{j\geq
1}\left(  1-\left\vert z_{j}\right\vert \right)  ^{\theta}\right\}
\sum_{z_{j}>z_{k}}\left(  1-\left\vert z_{j}\right\vert \right)  ^{\sigma}%
\leq\frac{1}{2}\sum_{z_{j}>z_{k}}\left(  1-\left\vert z_{j}\right\vert
\right)  ^{\sigma}.
\]
Thus $\sum_{z_{j}\geq z_{k}}\left(  1-\left\vert z_{j}\right\vert \right)
^{\sigma}\leq2\left(  1-\left\vert z_{k}\right\vert \right)  ^{\sigma}$,
proving (\ref{sigmasum}) for $\left\vert z_{k}\right\vert \geq R$. Now the
number of points $z_{k}$ in the ball $B\left(  0,R\right)  $ depends only on
$R $ and the separation constant $C$ in (\ref{sep}), and it is now easy to
obtain (\ref{sigmasum}) in general.

We will also use a lemma from \cite{Bo} which constructs a holomorphic
function $\varphi_{w}=\Gamma_{s}g_{w}$, where $\Gamma_{s}$ is the projection
operator below, that is close to $1$ on the Carleson region associated to a
point $w\in\mathbb{D}$, and decays appropriately away from the Carleson
region. Again let $1-C/2<\beta<1$ where $C$ is as in (\ref{sep}). Given
$\beta<\rho<\alpha<1$, we will use the cutoff function $c_{\rho,\alpha}$
defined by
\begin{equation}
c_{\rho,\alpha}\left(  \gamma\right)  =\left\{
\begin{array}
[c]{lll}%
0 & \text{for} & \gamma<\rho\\
\frac{\gamma-\rho}{\alpha-\rho} & \text{for} & \rho\leq\gamma\leq\alpha\\
1 & \text{for} & \alpha<\gamma
\end{array}
\right.  . \label{defc}%
\end{equation}

\begin{lemma}
\label{analyticcon}(Lemma 4.1 in \cite{Bo}) Suppose $s>-1$, $C$ is as in
(\ref{sep}), and $1-C/2<\beta<1$. There are $\beta_{1}$, $\rho$ and $\alpha$
satisfying $\beta<\beta_{1}<\rho<\alpha<1$ such that for every $w\in
\mathbb{D}$, we can find a function $g_{w}$ so that
\[
\varphi_{w}\left(  z\right)  =\Gamma_{s}g_{w}\left(  z\right)  =\int
_{\mathbb{D}}\frac{g_{w}\left(  \zeta\right)  (1-\left\vert \zeta\right\vert
^{2})^{s}}{\left(  1-\overline{\zeta}z\right)  ^{1+s}}d\zeta
\]
satisfies, with $c_{\rho,\alpha}$ is as in (\ref{defc}), and $\gamma
_{w}\left(  z\right)  $ is defined by $\left\vert z-\left\vert w\right\vert
^{-1}w\right\vert =\left(  1-\left\vert w\right\vert ^{2}\right)  ^{\gamma
_{w}\left(  z\right)  },$
\begin{equation}
\left\{
\begin{array}
[c]{llll}%
\varphi_{w}\left(  w\right)  & = & 1 & \\
\varphi_{w}\left(  z\right)  & = & c_{\rho,\alpha}\left(  \gamma_{w}\left(
z\right)  \right)  +O\left(  \left(  \log\frac{1}{1-\left\vert w\right\vert
^{2}}\right)  ^{-1}\right)  , & z\in V_{w}^{\beta}\\
\left\vert \varphi_{w}\left(  z\right)  \right\vert  & \leq & C\left(
\log\frac{1}{1-\left\vert w\right\vert ^{2}}\right)  ^{-1}\left(  1-\left\vert
w\right\vert ^{2}\right)  ^{\left(  \rho-\beta_{1}\right)  \left(  1+s\right)
}, & z\notin V_{w}^{\beta_{1}}%
\end{array}
\right.  , \label{satisfies}%
\end{equation}
Furthermore we have the estimate
\begin{equation}
\int_{\mathbb{D}}\left\vert g_{w}\left(  \zeta\right)  \right\vert ^{2}%
d\zeta\leq C\left(  \log\frac{1}{1-\left\vert w\right\vert ^{2}}\right)
^{-1}. \label{maincon'}%
\end{equation}

\end{lemma}

\subsubsection{The Sufficiency Proof\label{sp}}

We now prove that the hypotheses of Theorem \ref{Boeonto} are sufficient for interpolation.

Order the points $\left\{  z_{j}\right\}  _{j=1}^{\infty}$ so that
$1-\left\vert z_{j+1}\right\vert \leq1-\left\vert z_{j}\right\vert $ for
$j\geq1$. We now define a \textquotedblleft forest structure\textquotedblright%
\ on the index set $\mathbb{N}$ by declaring that $j$ is a child of $i$ (or
that $i$ is a parent of $j$) provided that
\begin{align}
i  &  <j,\label{daughter}\\
V_{z_{j}}  &  \subset V_{z_{i}},\nonumber\\
V_{z_{j}}  &  \varsubsetneq V_{z_{k}}\text{ for }i<k<j.\nonumber
\end{align}
Note if we have competing indices $i$ and $i^{\prime}$ with $V_{z_{j}}\subset
V_{z_{i}}\cap V_{z_{i^{\prime}}}$then the child $j$ chooses the
\textquotedblleft nearest\textquotedblright\ parent $i$. We define a partial
order associated with this parent-child relationship by declaring that $j$ is
a successor of $i$ (or that $i$ is a predecessor of $j$) if there is a
\textquotedblleft chain\textquotedblright\ of indices $\left\{  i=k_{1}%
,k_{2},...,k_{m}=j\right\}  \subset\mathbb{N}$ such that $k_{\ell+1}$ is a
child of $k_{\ell}$ for $1\leq\ell<m$. Under this partial ordering,
$\mathbb{N}$ decomposes into a disjoint union of trees. Thus associated to
each index $\ell\in\mathbb{N}$, there is a unique tree containing $\ell$ and,
unless $\ell$ is the root of the tree, a unique parent $P\left(  \ell\right)
$ of $\ell$ in that tree. Denote by $\mathcal{G}_{\ell}$ the unique geodesic
joining the root of the tree to $\ell$. We will usually identify $\ell$ with
$z_{\ell}$ and thereby transfer the forest structure $\mathcal{F}$ to $Z$ as well.

If $f(z)\in B_{2}$ and $f(z_{0})=0$ for some $z_{0}\in\mathbb{D}$ then
$f(z)/(z-z_{0})\in B_{2}.$ Using this it is easy to show that $Z$ is an onto
interpolating sequence if and only if some cofinite subsequence if $Z$ is.
With this observation, and recalling the hypothesis that $\left\Vert
\mu\right\Vert <\infty,$ we see that it suffices to do the proof under the
additional assumption that
\begin{equation}
\left\Vert \mu\right\Vert =\sum_{j=1}^{\infty}\mu\left(  z_{j}\right)
=\sum_{j=1}^{\infty}\left(  \log\frac{1}{1-\left\vert z_{j}\right\vert ^{2}%
}\right)  ^{-1}<\varepsilon. \label{logeps''}%
\end{equation}
Where $\varepsilon$ is a small quantity to be specified later. With this done
we now further suppose that the sequence $\left\{  z_{j}\right\}  _{j=1}^{J}$
is finite, and obtain an appropriate estimate independent of $J\geq1$. Fix
$\alpha,s>-1$ and a sequence of complex numbers $\left\{  \xi_{j}\right\}
_{j=1}^{J}$ in $\ell^{2}\left(  \mu\right)  $ where
\[
\left\Vert \left\{  \xi_{j}\right\}  _{j=1}^{J}\right\Vert _{\ell^{2}\left(
\mu\right)  }=\left\Vert \left\{  \frac{\xi_{j}}{\left\Vert k_{z_{j}}%
^{\alpha,2}\right\Vert _{B_{2}}}\right\}  _{j=1}^{J}\right\Vert _{\ell^{2}}.
\]
We will define a function $\varphi=\mathcal{S}\xi$ on the disk $\mathbb{D}$
by
\begin{equation}
\varphi\left(  z\right)  =\mathcal{S}\xi\left(  z\right)  =\sum_{j=1}^{J}%
a_{j}\varphi_{z_{j}}\left(  z\right)  ,\;\;\;\;\;z\in\mathbb{D},
\label{defofM}%
\end{equation}
that will be our candidate for the interpolating function of $\xi$. We follow
the inductive scheme of B\"{o}e that addresses the main difficulty in
interpolating holomorphic functions, namely that on the sequence $Z$ the
building blocks $\varphi_{z_{j}}$ take a large set of values (rather than just
$0$ and $1$ as in the tree analogue).

Recall that $Pz_{j}$ denotes the parent of $z_{j}$ in the forest structure
$\mathcal{F}$ and that $\mathcal{G}_{\ell}$ is the geodesic from the root to
$z_{\ell}$ in the tree containing $z_{\ell}$. In order to define the
coefficients $a_{j}$ we will use the doubly indexed sequence $\left\{
\beta_{i,j}\right\}  $ of numbers given by%
\begin{equation}
\beta_{i,j}=\varphi_{Pz_{j}}\left(  z_{i}\right)  . \label{defbetakmod}%
\end{equation}
We consider separately the indices in each tree of the forest $\left\{
1,2,...,J\right\}  $, and define the coefficients inductively according to the
\emph{natural} ordering of the integers. So let $\mathcal{Y}$ be a tree in the
forest $\left\{  1,2,...,J\right\}  $ with root $k_{0}$. Define $a_{k_{0}}%
=\xi_{k_{0}}$. Suppose that $k\in\mathcal{Y}\setminus\left\{  k_{0}\right\}  $
and that the coefficients $a_{j}$ have been defined for $j\in\mathcal{Y}$ and
$j<k$. Let%
\[
\mathcal{G}_{k}=\left[  k_{0},k\right]  =\left\{  k_{0},k_{1},...,k_{m-1}%
,k_{m}=k\right\}
\]
be the geodesic $\mathcal{G}_{k}$ in $\mathcal{Y}$ joining $k_{0}$ to $k$, and
note that $\mathcal{G}_{k}=\mathcal{G}_{k_{m-1}}\cup\left\{  k\right\}  $.
Define
\[
f_{k}\left(  z\right)  =f_{k_{m}}\left(  z\right)  =\sum_{i=1}^{m}a_{k_{i}%
}\varphi_{z_{k_{i}}}\left(  z\right)  =f_{k_{m-1}}\left(  z\right)
+a_{k}\varphi_{z_{k}}\left(  z\right)
\]
and
\[
\omega_{k}=f_{k_{m-1}}\left(  z_{k}\right)  =\sum_{i=1}^{m}a_{k_{i-1}}%
\varphi_{z_{k_{i-1}}}\left(  z_{k}\right)  =\sum_{i=1}^{m}\beta_{k,i}%
a_{k_{i-1}},\;\;\;\;\;k\geq1.
\]
Then define the coefficient $a_{k}$ by
\begin{equation}
a_{k}=\xi_{k}-\omega_{k},\;\;\;\;\;k\geq1. \label{amdef}%
\end{equation}
This completes the inductive definition of the sequence $\left\{
a_{k}\right\}  _{k\in\mathcal{Y}}$, and hence defines the entire sequence
$\left\{  a_{i}\right\}  _{i=1}^{J}$.

We first prove the following $\ell^{2}\left(  d\mu\right)  $ estimate for the
sequence $\left\{  a_{j}^{m}\right\}  _{j=1}^{J}$ given in terms of the data
$\left\{  \xi_{j}^{m}\right\}  _{j=1}^{J}$ by the scheme just introduced. This
is the difficult step in the proof of sufficiency.

\begin{lemma}
\label{weightcon}The sequence $\left\{  a_{i}\right\}  _{i=1}^{J}$ constructed
in (\ref{amdef}) above satisfies
\begin{equation}
\left\Vert \left\{  a_{j}\right\}  _{j=1}^{J}\right\Vert _{\ell^{2}\left(
d\mu\right)  }\leq C\left\Vert \left\{  \xi_{j}\right\}  _{j=1}^{J}\right\Vert
_{\ell^{2}\left(  d\mu\right)  }. \label{control}%
\end{equation}

\end{lemma}

\textbf{Proof}: Without loss of generality, we may assume for the purposes of
this proof that the forest of indices $\left\{  j\right\}  _{j=1}^{J}$ is
actually a single tree $\mathcal{Y}$. Now fix $\ell$. At this point it will be
convenient notation to momentarily relabel the points $\left\{  z_{j}\right\}
_{j\in\mathcal{G}_{\ell}}=\left\{  z_{k_{0}},z_{k_{1}},...,z_{k_{m}}\right\}
$ as $\left\{  z_{0},z_{1},...,z_{m}\right\}  $, with similar relabeling of
the and similarly relabel the $\left\{  \alpha_{j}\right\}  ,$ $\left\{
\zeta_{j}\right\}  ,$ and $\left\{  \beta_{j}\right\}  $ so that
\[
a_{k}=\xi_{k}-\sum_{i=1}^{k}\beta_{k,i}a_{i-1},\;\;\;\;\;0\leq k\leq\ell.
\]
In other words, we are restricting attention to the geodesic $\mathcal{G}%
_{\ell}$ and relabeling sequences so as to conform to the ordering in the
geodesic. We also rewrite $f_{k}\left(  z\right)  $ and $\omega_{k}$ as%
\[
f_{k}\left(  z\right)  =\sum_{i=1}^{k}a_{i}\varphi_{z_{i}}\left(  z\right)
=f_{k-1}\left(  z\right)  +a_{k}\varphi_{z_{k}}\left(  z\right)
\]
and
\begin{equation}
\omega_{k}=f_{k-1}\left(  z_{k}\right)  =\sum_{i=1}^{k}a_{i-1}\varphi
_{z_{i-1}}\left(  z_{k}\right)  =\sum_{i=1}^{k}\beta_{k,i}a_{i-1}%
,\;\;\;\;\;k\geq1. \label{defomega}%
\end{equation}
so that the coefficients $a_{k}$ are given by
\begin{align}
a_{0}  &  =\xi_{0},\label{defa}\\
a_{k}  &  =\xi_{k}-\omega_{k},\;\;\;\;\;k\geq1.\nonumber
\end{align}

We now claim that
\begin{equation}
\left(
\begin{array}
[c]{c}%
\omega_{1}\\
\omega_{2}\\
\vdots\\
\omega_{k}%
\end{array}
\right)  =\left[
\begin{array}
[c]{cccc}%
b_{1,1} & 0 & \cdots & 0\\
b_{2,1} & b_{2,2} & \cdots & 0\\
\vdots & \vdots & \ddots & \vdots\\
b_{k,1} & b_{k,2} & \cdots & b_{k,k}%
\end{array}
\right]  \left(
\begin{array}
[c]{c}%
\xi_{0}\\
\xi_{1}\\
\vdots\\
\xi_{k-1}%
\end{array}
\right)  ,\;\;\;\;\;1\leq k\leq\ell, \label{indclaim}%
\end{equation}
where
\begin{align}
b_{i,j}  &  =0,\;\;\;\;\;i<j,\label{bformula}\\
b_{i,i}  &  =\beta_{i,i},\nonumber\\
b_{i,j}  &  =b_{i-1,j}^{\ast}-b_{i-1,j}\beta_{i,i},\;\;\;\;\;i>j,\nonumber
\end{align}
and the $b_{i,j}^{\ast}$ are defined in the following calculations. We also
claim that the $b_{i,j}$ are bounded:%
\begin{equation}
\left\vert b_{i,j}\right\vert \leq C. \label{bijbounded}%
\end{equation}
For this we will use the estimate (see Lemma \ref{fax} below)
\begin{equation}
\left\vert \varphi_{w}^{\prime}\left(  z\right)  \right\vert \leq\left(
1-\left\vert w\right\vert ^{2}\right)  ^{-\alpha},\;\;\;\;\;z\in\mathbb{D}.
\label{derest}%
\end{equation}
Note first that
\[
b_{1,1}=\beta_{1,1}=\varphi_{z_{0}}\left(  z_{1}\right)
\]
since then (\ref{defomega}) and (\ref{defa}) yield
\[
b_{1,1}\xi_{0}=\omega_{1},
\]
which is (\ref{indclaim}) for $k=1$. We also have (\ref{bijbounded}) for
$1\leq j\leq i=1$ since (\ref{satisfies}) yields%
\[
\left\vert b_{1,1}\right\vert \leq1+\lambda\left(  z_{0}\right)  ,
\]
where we have introduced the convenient notation%
\[
\lambda\left(  z_{j}\right)  =\left(  \log\frac{1}{1-\left\vert z_{j}%
\right\vert ^{2}}\right)  ^{-1}.
\]
We now define a function $b_{1,1}\left(  z\right)  $ by
\[
b_{1,1}\left(  z\right)  =\varphi_{z_{0}}\left(  z\right)  ,
\]
i.e. we replace $z_{1}$ by $z$ throughout the formula for $b_{1,1}$. If we
then define
\[
b_{1,1}^{\ast}=b_{1,1}\left(  z_{2}\right)  =\varphi_{z_{0}}\left(
z_{2}\right)  ,
\]
we readily obtain
\begin{align*}
b_{2,1}  &  =b_{1,1}^{\ast}-b_{1,1}\beta_{2,2}=\varphi_{z_{0}}\left(
z_{2}\right)  -\varphi_{z_{0}}\left(  z_{1}\right)  \varphi_{z_{1}}\left(
z_{2}\right)  ,\\
b_{2,2}  &  =\beta_{2,2}=\varphi_{z_{1}}\left(  z_{2}\right)  .
\end{align*}
Indeed, from (\ref{defomega}), (\ref{defa}) and the equality $\xi_{1}%
=a_{1}+\omega_{1}=a_{1}+\varphi_{z_{0}}\left(  z_{1}\right)  a_{0}$, we have
\begin{align*}
b_{2,1}\xi_{0}+b_{2,2}\xi_{1}  &  =\left[  \varphi_{z_{0}}\left(
z_{2}\right)  -\varphi_{z_{0}}\left(  z_{1}\right)  \varphi_{z_{1}}\left(
z_{2}\right)  \right]  a_{0}+\varphi_{z_{1}}\left(  z_{2}\right)  \left[
a_{1}+\varphi_{z_{0}}\left(  z_{1}\right)  a_{0}\right] \\
&  =\varphi_{z_{0}}\left(  z_{2}\right)  a_{0}+\varphi_{z_{1}}\left(
z_{2}\right)  a_{1}\\
&  =\omega_{2},
\end{align*}
which proves (\ref{indclaim}) for $k=2$. We also have (\ref{bijbounded}) for
$1\leq j\leq i=2$ since the bound%
\[
\left\vert b_{2,2}\right\vert \leq1+\lambda\left(  z_{1}\right)
\]
is obvious from (\ref{satisfies}), and the bound for $b_{2,1}$ follows from
(\ref{satisfies}), (\ref{derest}) and Lemma \ref{Mars}:%
\begin{align*}
\left\vert b_{2,1}\right\vert  &  \leq\left\vert \varphi_{z_{0}}\left(
z_{2}\right)  -\varphi_{z_{0}}\left(  z_{1}\right)  \right\vert +\left\vert
\varphi_{z_{0}}\left(  z_{1}\right)  \right\vert \left\vert 1-\varphi_{z_{1}%
}\left(  z_{2}\right)  \right\vert \\
&  \leq\left\vert \varphi_{z_{0}}^{\prime}\left(  \zeta_{0}\right)
\right\vert \left\vert z_{2}-z_{1}\right\vert +\left(  1+\lambda\left(
z_{0}\right)  \right)  \left(  1+\lambda\left(  z_{1}\right)  \right)  ,
\end{align*}
and since%
\begin{align*}
\left\vert \varphi_{z_{0}}^{\prime}\left(  \zeta_{0}\right)  \right\vert
\left\vert z_{2}-z_{1}\right\vert  &  \leq\left(  1-\left\vert z_{0}%
\right\vert ^{2}\right)  ^{-\alpha}\left\vert z_{2}-z_{1}\right\vert \\
&  \leq\left(  1-\left\vert z_{0}\right\vert ^{2}\right)  ^{-\alpha}\left(
1-\left\vert z_{1}\right\vert ^{2}\right)  ^{\beta}\\
&  \leq\left(  1-\left\vert z_{0}\right\vert ^{2}\right)  ^{\beta\eta-\alpha},
\end{align*}
we obtain%
\begin{equation}
\left\vert b_{2,1}\right\vert \leq\left(  1-\left\vert z_{0}\right\vert
^{2}\right)  ^{\beta\eta-\alpha}+e^{\lambda\left(  z_{0}\right)
+\lambda\left(  z_{1}\right)  }. \label{b21}%
\end{equation}

We now define functions $b_{2,1}\left(  z\right)  $ and $b_{2,2}\left(
z\right)  $ by
\begin{align*}
b_{2,1}\left(  z\right)   &  =\varphi_{z_{0}}\left(  z\right)  -b_{1,1}%
\varphi_{z_{1}}\left(  z\right)  ,\\
b_{2,2}\left(  z\right)   &  =\varphi_{z_{1}}\left(  z\right)  ,
\end{align*}
i.e. we replace $z_{2}$ by $z$ throughout the formulas for $b_{2,1}$ and
$b_{2,2}$. If we then set
\begin{align*}
b_{2,1}^{\ast}  &  =b_{2,1}\left(  z_{3}\right)  ,\\
b_{2,2}^{\ast}  &  =b_{2,2}\left(  z_{3}\right)  ,
\end{align*}
we obtain as above that
\begin{align*}
b_{3,1}  &  =b_{2,1}^{\ast}-b_{2,1}\beta_{3,3}=\left[  \varphi_{z_{0}}\left(
z_{3}\right)  -b_{1,1}\varphi_{z_{1}}\left(  z_{3}\right)  \right]  -\left[
\varphi_{z_{0}}\left(  z_{2}\right)  -b_{1,1}\varphi_{z_{1}}\left(
z_{2}\right)  \right]  \varphi_{z_{2}}\left(  z_{3}\right)  ,\\
b_{3,2}  &  =b_{2,2}^{\ast}-b_{2,2}\beta_{3,3}=\varphi_{z_{1}}\left(
z_{3}\right)  -\varphi_{z_{1}}\left(  z_{2}\right)  \varphi_{z_{2}}\left(
z_{3}\right)  ,\\
b_{3,3}  &  =\beta_{3,3}=\varphi_{z_{2}}\left(  z_{3}\right)  ,
\end{align*}
which proves (\ref{indclaim}) for $k=3$. Moreover, we again have
(\ref{bijbounded}) for $1\leq j\leq i=3$. Indeed,%
\[
\left\vert b_{3,3}\right\vert \leq1+\lambda\left(  z_{2}\right)  ,
\]
and the arguments used above to obtain (\ref{b21}) show that%
\[
\left\vert b_{3,2}\right\vert \leq\left(  1-\left\vert z_{1}\right\vert
^{2}\right)  ^{\beta\eta-\alpha}+e^{\lambda\left(  z_{1}\right)
+\lambda\left(  z_{2}\right)  }.
\]
Finally,%
\begin{align*}
\left\vert b_{3,1}\right\vert  &  \leq\left\vert b_{2,1}^{\ast}-b_{2,1}%
\right\vert +\left\vert b_{2,1}\right\vert \left\vert 1-\beta_{3,3}\right\vert
\\
&  \leq\left\{  \left\vert \varphi_{z_{0}}^{\prime}\left(  \zeta_{0}\right)
\right\vert +\left\vert b_{1,1}\right\vert \left\vert \varphi_{z_{1}}^{\prime
}\left(  \zeta_{1}\right)  \right\vert \right\}  \left\vert z_{2}%
-z_{3}\right\vert \\
&  +\left(  1-\left\vert z_{0}\right\vert ^{2}\right)  ^{\beta\eta-\alpha
}+e^{\lambda\left(  z_{0}\right)  +\lambda\left(  z_{1}\right)  +\lambda
\left(  z_{2}\right)  },
\end{align*}
and since%
\begin{align*}
&  \left\{  \left\vert \varphi_{z_{0}}^{\prime}\left(  \zeta_{0}\right)
\right\vert +\left\vert b_{1,1}\right\vert \left\vert \varphi_{z_{1}}^{\prime
}\left(  \zeta_{1}\right)  \right\vert \right\}  \left\vert z_{2}%
-z_{3}\right\vert \\
&  \ \ \ \ \ \leq\left\{  \left(  1-\left\vert z_{0}\right\vert ^{2}\right)
^{-\alpha}+\left(  1+\lambda\left(  z_{0}\right)  \right)  \left(
1-\left\vert z_{1}\right\vert ^{2}\right)  ^{-\alpha}\right\}  \left\vert
z_{2}-z_{3}\right\vert \\
&  \ \ \ \ \ \leq\left\{  \left(  1-\left\vert z_{1}\right\vert ^{2}\right)
^{-\frac{\alpha}{\eta}}+\left(  1+\lambda\left(  z_{0}\right)  \right)
\left(  1-\left\vert z_{1}\right\vert ^{2}\right)  ^{-\alpha}\right\}  \left(
1-\left\vert z_{2}\right\vert ^{2}\right)  ^{\beta}\\
&  \ \ \ \ \ \leq\left\{  \left(  1+A\lambda\left(  z_{0}\right)  \right)
\left(  1-\left\vert z_{1}\right\vert ^{2}\right)  ^{-\alpha}\right\}  \left(
1-\left\vert z_{1}\right\vert ^{2}\right)  ^{\beta\eta},\\
&  \ \ \ \ \ \leq\left(  1+A\lambda\left(  z_{0}\right)  \right)  \left(
1-\left\vert z_{1}\right\vert ^{2}\right)  ^{\beta\eta-\alpha},
\end{align*}
for a large constant $A$, we have%
\[
\left\vert b_{3,1}\right\vert \leq e^{A\lambda\left(  z_{0}\right)  }\left(
1-\left\vert z_{1}\right\vert ^{2}\right)  ^{\beta\eta-\alpha}+e^{\lambda
\left(  z_{0}\right)  +\lambda\left(  z_{1}\right)  +\lambda\left(
z_{2}\right)  }.
\]
Continuing in this way with
\begin{align}
b_{i,j}\left(  z\right)   &  =b_{i-1,j}\left(  z\right)  -b_{i-1,j}%
\varphi_{z_{i-1}}\left(  z\right)  ,\label{continuing}\\
b_{i,j}  &  =b_{i,j}\left(  z_{i}\right)  ,\nonumber\\
b_{i,j}^{\ast}  &  =b_{i,j}\left(  z_{i+1}\right)  ,\nonumber
\end{align}
we can prove (\ref{indclaim}) and (\ref{bijbounded}) by induction on $k$ and
$i$ (see below). The bound $C$ in (\ref{bijbounded}) will use the fact that%
\begin{equation}
\lambda\left(  z_{0}\right)  +\lambda\left(  z_{1}\right)  +\lambda\left(
z_{2}\right)  +...\leq C\lambda\left(  z_{0}\right)  . \label{lambdas}%
\end{equation}
To see (\ref{lambdas}) we use Lemma \ref{Mars}.

Now if $\mathcal{G}_{\ell}=\left[  k_{0},k_{1},...,k_{m-1},k_{m}\right]  $,
then by applying (\ref{ji}) repeatedly, we obtain
\[
\left(  1-\left\vert z_{k_{i}}\right\vert ^{2}\right)  \leq\left(
1-\left\vert z_{k_{0}}\right\vert ^{2}\right)  ^{\eta^{i}},
\]
and so combining these estimates we have
\begin{align}
\lambda\left(  z_{0}\right)  +\lambda\left(  z_{1}\right)  +\lambda\left(
z_{2}\right)  +...  &  \leq C\sum_{i\in\mathcal{G}_{\ell}\setminus\left\{
k_{0}\right\}  }\left(  \log\frac{1}{1-\left\vert z_{P\left(  i\right)
}\right\vert ^{2}}\right)  ^{-1}\label{sumlogs}\\
&  \leq C\left(  \sum_{j=0}^{m-1}\eta^{-j}\right)  \left(  \log\frac
{1}{1-\left\vert z_{k_{0}}\right\vert ^{2}}\right)  ^{-1}\nonumber\\
&  \leq C_{\eta}\left(  \log\frac{1}{1-\left\vert z_{k_{0}}\right\vert ^{2}%
}\right)  ^{-1}=C_{\eta}\lambda\left(  z_{0}\right) \nonumber
\end{align}
since $\eta>1$, which yields (\ref{lambdas}).

\bigskip

We now give the induction details for proving (\ref{indclaim}) and
(\ref{bijbounded}). The proof of (\ref{indclaim}) is straightforward by
induction on $k$, so we concentrate on proving (\ref{bijbounded}) by induction
on $i$. If we denote the $i^{th}$ row
\[
\left[
\begin{array}
[c]{ccccccc}%
b_{i,1} & b_{i,2} & \cdots & b_{i,i} & 0 & \cdots & 0
\end{array}
\right]
\]
of the matrix in (\ref{indclaim}) by $\mathbf{B}_{i}$, the corresponding row
of starred components
\[
\left[
\begin{array}
[c]{ccccccc}%
b_{i,1}^{\ast} & b_{i,2}^{\ast} & \cdots & b_{i,i}^{\ast} & 0 & \cdots & 0
\end{array}
\right]
\]
by $\mathbf{B}_{i}^{\ast}$, and the row having all zeroes except a one in the
$i^{th}$ place by $\mathbf{E}_{i}$, then we have the recursion formula
\begin{align}
\mathbf{B}_{i}  &  =\mathbf{B}_{i-1}^{\ast}-\mathbf{B}_{i-1}+\left(
1-\beta_{i,i}\right)  \mathbf{B}_{i-1}+\beta_{i,i}\mathbf{E}_{i}%
\label{recfor}\\
&  =\left\{  \left(  1-\beta_{i,i}\right)  \mathbf{B}_{i-1}+\beta
_{i,i}\mathbf{E}_{i}\right\}  -\left(  \mathbf{B}_{i-1}-\mathbf{B}_{i-1}%
^{\ast}\right) \nonumber
\end{align}
which expresses $\mathbf{B}_{i}$ as a \textquotedblleft convex
combination\textquotedblright\ of the previous row and the unit row
$\mathbf{E}_{i}$, minus the difference of the previous row and its starred
counterpart. In terms of the components of the rows, we have
\begin{equation}
b_{i,j}=\left[  b_{i-1,j}^{\ast}-b_{i-1,j}\right]  +\left(  1-\beta
_{i,i}\right)  b_{i-1,j}+\beta_{i,i}\delta_{i,j}. \label{recforcomp}%
\end{equation}

For a large constant $A$ that will be chosen later so that the induction step
works, we prove the following estimate by induction on $i$:%
\begin{equation}
\left\vert b_{i,j}\right\vert \leq e^{A\left\{  \lambda\left(  z_{j-1}\right)
+...+\lambda\left(  z_{i-1}\right)  \right\}  },\ \ \ \ \ i\geq j.
\label{iest}%
\end{equation}
The initial case $i=j$ follows from%
\[
\left\vert b_{j,j}\right\vert =\left\vert \varphi_{z_{j-1}}\left(
z_{j}\right)  \right\vert \leq1+\lambda\left(  z_{j-1}\right)  \leq
e^{\lambda\left(  z_{j-1}\right)  }.
\]
Now (\ref{continuing}) yields%
\[
b_{i,j}^{\prime}\left(  z\right)  =b_{i-1,j}^{\prime}\left(  z\right)
-b_{i-1,j}\varphi_{z_{i-1}}^{\prime}\left(  z\right)  ,
\]
and so by the induction assumption for indices smaller than $i$, we have from
(\ref{derest}) that%
\begin{align}
\left\Vert b_{i,j}^{\prime}\right\Vert _{L^{\infty}}  &  \leq\left\vert
b_{i-1,j}\right\vert \left(  1-\left\vert z_{i-1}\right\vert ^{2}\right)
^{-\alpha}+\left\Vert b_{i-1,j}^{\prime}\right\Vert _{L^{\infty}%
}\label{bijprime}\\
&  \leq\left\vert b_{i-1,j}\right\vert \left(  1-\left\vert z_{i-1}\right\vert
^{2}\right)  ^{-\alpha}+\left\vert b_{i-2,j}\right\vert \left(  1-\left\vert
z_{i-2}\right\vert ^{2}\right)  ^{-\alpha}+\left\Vert b_{i-2,j}^{\prime
}\right\Vert _{L^{\infty}}\nonumber\\
&  \vdots\nonumber\\
&  \leq\left\{  \sup_{j\leq k\leq i-1}\left\vert b_{k,j}\right\vert \right\}
\left[  \left(  1-\left\vert z_{i-1}\right\vert ^{2}\right)  ^{-\alpha
}+...+\left(  1-\left\vert z_{j}\right\vert ^{2}\right)  ^{-\alpha}\right]
+\left(  1-\left\vert z_{j-1}\right\vert ^{2}\right)  ^{-\alpha}\nonumber\\
&  \leq e^{A\left\{  \lambda\left(  z_{j-1}\right)  +...+\lambda\left(
z_{i-2}\right)  \right\}  }\left[  \left(  1-\left\vert z_{i-1}\right\vert
^{2}\right)  ^{-\alpha}+...+\left(  1-\left\vert z_{j-1}\right\vert
^{2}\right)  ^{-\alpha}\right] \nonumber\\
&  \leq e^{A\left\{  \lambda\left(  z_{j-1}\right)  +...+\lambda\left(
z_{i-2}\right)  \right\}  }\left[  \left(  1-\left\vert z_{i-1}\right\vert
^{2}\right)  ^{-\alpha}+...+\left(  1-\left\vert z_{i-1}\right\vert
^{2}\right)  ^{-\frac{\alpha}{\eta^{i-j}}}\right]  ,\nonumber
\end{align}
where the last line uses (\ref{ji}). Thus we have from (\ref{bijprime}),
(\ref{recforcomp}) and (\ref{continuing}),%
\begin{align}
\left\vert b_{i,j}\right\vert  &  \leq\left\vert b_{i-1,j}^{\ast}%
-b_{i-1,j}\right\vert +\left\vert \left(  1-\beta_{i,i}\right)  b_{i-1,j}%
\right\vert \label{bijpre}\\
&  \leq\left\vert b_{i-1,j}\left(  z_{i+1}\right)  -b_{i-1,j}\left(
z_{i}\right)  \right\vert +\left\vert 1-\beta_{i,i}\right\vert \left\vert
b_{i-1,j}\right\vert \nonumber\\
&  \leq\left\Vert b_{i,j}^{\prime}\right\Vert _{L^{\infty}}\left\vert
z_{i+1}-z_{i}\right\vert +\left(  1+\lambda\left(  z_{i-1}\right)  \right)
\left\vert b_{i-1,j}\right\vert \nonumber\\
&  \leq e^{A\left\{  \lambda\left(  z_{j-1}\right)  +...+\lambda\left(
z_{i-2}\right)  \right\}  }\times\nonumber\\
&  \left\{  1+\lambda\left(  z_{i-1}\right)  +\left[  \left(  1-\left\vert
z_{i-1}\right\vert ^{2}\right)  ^{-\alpha}+...+\left(  1-\left\vert
z_{i-1}\right\vert ^{2}\right)  ^{-\frac{\alpha}{\eta^{i-j}}}\right]  \left(
1-\left\vert z_{i-1}\right\vert ^{2}\right)  ^{\beta\eta}\right\}  ,\nonumber
\end{align}
upon using the inequality $\left\vert z_{i+1}-z_{i}\right\vert \leq\left(
1-\left\vert z_{i}\right\vert ^{2}\right)  ^{\beta}\leq\left(  1-\left\vert
z_{i-1}\right\vert ^{2}\right)  ^{\beta\eta}$, which follows from Lemma
\ref{Mars}.

Finally we use the inequality (see below for a proof)%
\begin{equation}
\left(  1-\left\vert z_{i-1}\right\vert ^{2}\right)  ^{-\alpha}+...+\left(
1-\left\vert z_{i-1}\right\vert ^{2}\right)  ^{-\frac{\alpha}{\eta^{i-j}}}%
\leq\left\{  C_{\eta}\left(  1-\left\vert z_{i-1}\right\vert ^{2}\right)
^{-\alpha}+2\left(  i-j+1\right)  \right\}  \label{thein}%
\end{equation}
to obtain that%
\begin{align}
&  \left[  \left(  1-\left\vert z_{i-1}\right\vert ^{2}\right)  ^{-\alpha
}+...+\left(  1-\left\vert z_{i-1}\right\vert ^{2}\right)  ^{-\frac{\alpha
}{\eta^{i-j}}}\right]  \left(  1-\left\vert z_{i-1}\right\vert ^{2}\right)
^{\beta\eta}\label{dompower}\\
&  \ \ \ \ \ \leq2\left\{  C_{\eta}\left(  1-\left\vert z_{i-1}\right\vert
^{2}\right)  ^{\beta\eta-\alpha}+2i\left(  1-\left\vert z_{i-1}\right\vert
^{2}\right)  ^{\beta\eta}\right\} \nonumber\\
&  \ \ \ \ \ \leq\left(  A-1\right)  \left(  \log\frac{1}{1-\left\vert
z_{i-1}\right\vert ^{2}}\right)  ^{-1}\nonumber\\
&  \ \ \ \ \ =\left(  A-1\right)  \lambda\left(  z_{i-1}\right)  ,\nonumber
\end{align}
for all $i$ if $A$ is chosen large enough. With such a choice of $A$,
(\ref{bijpre}) yields%
\begin{align*}
\left\vert b_{i,j}\right\vert  &  \leq e^{A\left\{  \lambda\left(
z_{j-1}\right)  +...+\lambda\left(  z_{i-2}\right)  \right\}  }\left\{
1+\lambda\left(  z_{i-1}\right)  +\left(  A-1\right)  \lambda\left(
z_{i-1}\right)  \right\} \\
&  \leq e^{A\left\{  \lambda\left(  z_{j-1}\right)  +...+\lambda\left(
z_{i-2}\right)  +\lambda\left(  z_{i-1}\right)  \right\}  },
\end{align*}
which proves (\ref{iest}), and hence (\ref{bijbounded}) by (\ref{lambdas}). To
see (\ref{thein}), we rewrite it as%
\[
\sum_{\ell=0}^{N}R^{\eta^{-\ell}}\leq C_{\eta}R+2N+2,
\]
and to prove this, note that for $R^{\eta^{-\ell}}>2$ the ratio of the
consecutive terms $R^{\eta^{-\left(  \ell+1\right)  }}$ and $R^{\eta^{-\ell}}$
is $\left(  R^{\eta^{-\ell}}\right)  ^{1-\eta}<2^{1-\eta}$, i.e. this portion
of the series is supergeometric. Thus we have%
\begin{align*}
\sum_{\ell=0}^{N}R^{\eta^{-\ell}}  &  \leq\sum_{\ell\geq0:R^{\eta^{-\ell}}%
>2}R^{\eta^{-\ell}}+\sum_{\ell\leq N:R^{\eta^{-\ell}}\leq2}R^{\eta^{-\ell}}\\
&  \leq R\sum_{j=0}^{\infty}\left(  2^{1-\eta}\right)  ^{j}+2\left(
N+1\right)  .
\end{align*}

We now claim the following crucial property. Recall that $P\left(  m\right)
=m-1$. If $\sigma>0$ and $\gamma_{m-1}\left(  z_{m}\right)  >\alpha+\sigma$,
i.e. $z_{m}\in V_{z_{m-1}}^{\alpha+\sigma}$, then
\begin{equation}
\left\vert b_{i,j}\right\vert \leq C\left(  1-\left\vert z_{m-1}\right\vert
^{2}\right)  ^{\sigma}\text{ for all }j<m\leq i. \label{crucialprop}%
\end{equation}
We first note that from (\ref{derest}), we have for $z_{m}\in V_{z_{m-1}%
}^{\alpha+\sigma}$,
\begin{align}
\beta_{m,m}  &  =\varphi_{z_{m-1}}\left(  z_{m}\right)  =\varphi_{z_{m-1}%
}\left(  z_{m-1}\right)  +\left[  \varphi_{z_{m-1}}\left(  z_{m}\right)
-\varphi_{z_{m-1}}\left(  z_{m-1}\right)  \right] \label{betamm}\\
&  =1+O\left(  \left(  1-\left\vert z_{m-1}\right\vert ^{2}\right)  ^{-\alpha
}\left\vert z_{m}-z_{m-1}\right\vert \right) \nonumber\\
&  =1+O\left(  \left(  1-\left\vert z_{m-1}\right\vert ^{2}\right)  ^{-\alpha
}\left(  1-\left\vert z_{m-1}\right\vert ^{2}\right)  ^{\alpha+\sigma}\right)
\nonumber\\
&  =1+O\left(  \left(  1-\left\vert z_{m-1}\right\vert ^{2}\right)  ^{\sigma
}\right)  .\nonumber
\end{align}
From (\ref{recfor}) we then obtain
\[
\left\Vert \mathbf{B}_{m}-\beta_{m,m}\mathbf{E}_{m}\right\Vert _{\infty}%
\leq\left\Vert \mathbf{B}_{m-1}^{\ast}-\mathbf{B}_{m-1}\right\Vert _{\infty
}+O\left(  \left(  1-\left\vert z_{m-1}\right\vert ^{2}\right)  ^{\sigma
}\right)  \left\Vert \mathbf{B}_{m-1}\right\Vert _{\infty}.
\]
Next, the estimate
\[
\left\Vert \mathbf{B}_{m-1}^{\ast}-\mathbf{B}_{m-1}\right\Vert _{\infty}\leq
C\left(  1-\left\vert z_{m-1}\right\vert ^{2}\right)  ^{\sigma},
\]
follows from (\ref{continuing}), (\ref{bijprime}) and (\ref{dompower}) with
$\beta\eta$ replaced with $\alpha+\sigma$:
\begin{align*}
\left\vert b_{m-1,j}^{\ast}-b_{m-1,j}\right\vert  &  =\left\vert
b_{m-1,j}\left(  z_{m}\right)  -b_{m-1,j}\left(  z_{m-1}\right)  \right\vert
\\
&  \leq\left\Vert b_{m-1,j}^{\prime}\right\Vert _{L^{\infty}}\left\vert
z_{m}-z_{m-1}\right\vert \\
&  \leq\sum_{k=1}^{m}\left\{  C_{\eta}\left(  1-\left\vert z_{k-1}\right\vert
^{2}\right)  ^{-\alpha}+2k\right\}  \left(  1-\left\vert z_{m-1}\right\vert
^{2}\right)  ^{\alpha+\sigma}\\
&  \leq C_{\sigma}\left(  1-\left\vert z_{m-1}\right\vert ^{2}\right)
^{\sigma}.
\end{align*}
Thus altogether we have proved that the top row of the rectangle
$\mathbf{R}_{m}=\left[  b_{i,j}\right]  _{j<m\leq i}$ satisfies
(\ref{crucialprop}), i.e. $b_{m,j}\leq C\left(  1-\left\vert z_{m-1}%
\right\vert ^{2}\right)  ^{\sigma}$ for $j<m$. The proof for the remaining
rows is similar using (\ref{recforcomp}).

For convenience in notation we now define
\[
\Gamma=\left\{  m:z_{m}\in V_{z_{m-1}}^{\alpha+\sigma}\right\}  .
\]
If we take $0<\sigma\leq\left(  \eta-1\right)  \alpha$ and iterate the proof
of (\ref{crucialprop}) and use (\ref{bijbounded}), we obtain the improved
estimate
\begin{equation}
\left\vert b_{i,j}\right\vert \leq C\prod_{m\in\Gamma:j<m\leq i}\left(
1-\left\vert z_{m-1}\right\vert ^{2}\right)  ^{\sigma},\;\;\;\;\;i>j.
\label{impest}%
\end{equation}
To see this we first look at the simplest case when $2,3\in\Gamma$ and
establish the corresponding inequality%
\begin{equation}
\left\vert b_{3,1}\right\vert \leq C\left(  1-\left\vert z_{1}\right\vert
^{2}\right)  ^{\sigma}\left(  1-\left\vert z_{2}\right\vert ^{2}\right)
^{\sigma}. \label{23simple}%
\end{equation}
We have from (\ref{recforcomp}) that%
\[
b_{3,1}=\left[  b_{2,1}^{\ast}-b_{2,1}\right]  +\left(  1-\beta_{3,3}\right)
b_{2,1}.
\]
From (\ref{betamm}) and (\ref{crucialprop}) we have%
\[
\left\vert b_{2,1}\right\vert \left\vert 1-\beta_{3,3}\right\vert \leq
C\left(  1-\left\vert z_{1}\right\vert ^{2}\right)  ^{\sigma}\left(
1-\left\vert z_{2}\right\vert ^{2}\right)  ^{\sigma}.
\]
From (\ref{continuing}) we have%
\begin{align*}
\left\vert b_{2,1}-b_{2,1}^{\ast}\right\vert  &  =\left\vert b_{2,1}\left(
z_{2}\right)  -b_{2,1}\left(  z_{3}\right)  \right\vert \\
&  =\left\vert \left[  b_{1,1}\left(  z_{2}\right)  -b_{1,1}\varphi_{z_{1}%
}\left(  z_{2}\right)  \right]  -\left[  b_{1,1}\left(  z_{3}\right)
-b_{1,1}\varphi_{z_{1}}\left(  z_{3}\right)  \right]  \right\vert \\
&  \leq\left\vert b_{1,1}\left(  z_{2}\right)  -b_{1,1}\left(  z_{3}\right)
\right\vert +\left\vert b_{1,1}\right\vert \left\vert \varphi_{z_{1}}\left(
z_{2}\right)  -\varphi_{z_{1}}\left(  z_{3}\right)  \right\vert \\
&  \leq\left\Vert b_{1,1}^{\prime}\right\Vert _{\infty}\left\vert z_{2}%
-z_{3}\right\vert +\left\vert b_{1,1}\right\vert \left\Vert \varphi_{z_{1}%
}^{\prime}\right\Vert _{\infty}\left\vert z_{2}-z_{3}\right\vert
\end{align*}
where%
\begin{align*}
\left\Vert b_{1,1}^{\prime}\right\Vert _{\infty}\left\vert z_{2}%
-z_{3}\right\vert  &  \leq C\left(  1-\left\vert z_{0}\right\vert ^{2}\right)
^{-\alpha}\left(  1-\left\vert z_{2}\right\vert ^{2}\right)  ^{\alpha+\sigma
}\\
&  \leq C\left(  1-\left\vert z_{1}\right\vert ^{2}\right)  ^{-\frac{\alpha
}{\eta}}\left(  1-\left\vert z_{1}\right\vert ^{2}\right)  ^{\eta\alpha
}\left(  1-\left\vert z_{2}\right\vert ^{2}\right)  ^{\sigma}%
\end{align*}
and%
\begin{align*}
\left\vert b_{1,1}\right\vert \left\Vert \varphi_{z_{1}}^{\prime}\right\Vert
_{\infty}\left\vert z_{2}-z_{3}\right\vert  &  \leq C\left(  1-\left\vert
z_{1}\right\vert ^{2}\right)  ^{-\alpha}\left(  1-\left\vert z_{2}\right\vert
^{2}\right)  ^{\alpha+\sigma}\\
&  \leq C\left(  1-\left\vert z_{1}\right\vert ^{2}\right)  ^{\left(
\eta-1\right)  \alpha}\left(  1-\left\vert z_{2}\right\vert ^{2}\right)
^{\sigma}%
\end{align*}
are both dominated by $C\left(  1-\left\vert z_{1}\right\vert ^{2}\right)
^{\sigma}\left(  1-\left\vert z_{2}\right\vert ^{2}\right)  ^{\sigma}$ if
$0<\sigma\leq\left(  \eta-1\right)  \alpha$. Altogether we have proved
(\ref{23simple}).

Now we suppose that $4\in\Gamma$ as well and prove the estimate%
\begin{equation}
\left\vert b_{4,1}\right\vert \leq C\left(  1-\left\vert z_{1}\right\vert
^{2}\right)  ^{\sigma}\left(  1-\left\vert z_{2}\right\vert ^{2}\right)
^{\sigma}\left(  1-\left\vert z_{3}\right\vert ^{2}\right)  ^{\sigma}.
\label{234simple}%
\end{equation}
Again we have from (\ref{recforcomp}) that%
\[
b_{4,1}=\left[  b_{3,1}^{\ast}-b_{3,1}\right]  +\left(  1-\beta_{4,4}\right)
b_{3,1},
\]
and from (\ref{betamm}) and (\ref{23simple}) we have%
\[
\left\vert b_{3,1}\right\vert \left\vert 1-\beta_{4,4}\right\vert \leq
C\left(  1-\left\vert z_{1}\right\vert ^{2}\right)  ^{\sigma}\left(
1-\left\vert z_{2}\right\vert ^{2}\right)  ^{\sigma}\left(  1-\left\vert
z_{3}\right\vert ^{2}\right)  ^{\sigma}.
\]
From (\ref{continuing}) we have%
\begin{align*}
\left\vert b_{3,1}-b_{3,1}^{\ast}\right\vert  &  =\left\vert b_{3,1}\left(
z_{3}\right)  -b_{3,1}\left(  z_{4}\right)  \right\vert \\
&  =\left\vert \left[  b_{2,1}\left(  z_{3}\right)  -b_{2,1}\varphi_{z_{2}%
}\left(  z_{3}\right)  \right]  -\left[  b_{2,1}\left(  z_{4}\right)
-b_{2,1}\varphi_{z_{2}}\left(  z_{4}\right)  \right]  \right\vert \\
&  \leq\left\vert b_{2,1}\left(  z_{3}\right)  -b_{2,1}\left(  z_{4}\right)
\right\vert +\left\vert b_{2,1}\right\vert \left\vert \varphi_{z_{2}}\left(
z_{3}\right)  -\varphi_{z_{2}}\left(  z_{4}\right)  \right\vert \\
&  \leq\left\Vert b_{2,1}^{\prime}\right\Vert _{\infty}\left\vert z_{3}%
-z_{4}\right\vert +\left\vert b_{2,1}\right\vert \left\Vert \varphi_{z_{2}%
}^{\prime}\right\Vert _{\infty}\left\vert z_{3}-z_{4}\right\vert ,
\end{align*}
where%
\begin{align*}
\left\Vert b_{2,1}^{\prime}\right\Vert _{\infty}\left\vert z_{3}%
-z_{4}\right\vert  &  \leq C\left(  1-\left\vert z_{1}\right\vert ^{2}\right)
^{-\alpha}\left(  1-\left\vert z_{3}\right\vert ^{2}\right)  ^{\alpha+\sigma
}\\
&  \leq C\left(  1-\left\vert z_{1}\right\vert ^{2}\right)  ^{-\alpha}\left(
1-\left\vert z_{2}\right\vert ^{2}\right)  ^{\eta\alpha}\left(  1-\left\vert
z_{3}\right\vert ^{2}\right)  ^{\sigma}\\
&  \leq C\left(  1-\left\vert z_{1}\right\vert ^{2}\right)  ^{\left(
\eta-1\right)  \alpha}\left(  1-\left\vert z_{2}\right\vert ^{2}\right)
^{\left(  \eta-1\right)  \alpha}\left(  1-\left\vert z_{3}\right\vert
^{2}\right)  ^{\sigma}%
\end{align*}
and%
\begin{align*}
\left\vert b_{2,1}\right\vert \left\Vert \varphi_{z_{2}}^{\prime}\right\Vert
_{\infty}\left\vert z_{3}-z_{4}\right\vert  &  \leq C\left(  1-\left\vert
z_{1}\right\vert ^{2}\right)  ^{\sigma}\left(  1-\left\vert z_{2}\right\vert
^{2}\right)  ^{-\alpha}\left(  1-\left\vert z_{3}\right\vert ^{2}\right)
^{\alpha+\sigma}\\
&  \leq C\left(  1-\left\vert z_{1}\right\vert ^{2}\right)  ^{\sigma}\left(
1-\left\vert z_{2}\right\vert ^{2}\right)  ^{\left(  \eta-1\right)  \alpha
}\left(  1-\left\vert z_{3}\right\vert ^{2}\right)  ^{\sigma}%
\end{align*}
are both dominated by $C\left(  1-\left\vert z_{1}\right\vert ^{2}\right)
^{\sigma}\left(  1-\left\vert z_{2}\right\vert ^{2}\right)  ^{\sigma}\left(
1-\left\vert z_{3}\right\vert ^{2}\right)  ^{\sigma}$ if $0<\sigma\leq\left(
\eta-1\right)  \alpha$. This completes the proof of (\ref{234simple}), and the
general case is similar.

The consequence we need from (\ref{impest}) is that if $m_{1}<m_{2}%
<...<m_{N}\leq k$ is an enumeration of the $m\in\Gamma$ such that $m\leq k$,
then
\begin{align}
\left\vert a_{k}\right\vert  &  \leq\left\vert \xi_{k}-\omega_{k}\right\vert
\label{consequence}\\
&  \leq\left\vert \xi_{k}\right\vert +\left\vert \omega_{k}\right\vert
\nonumber\\
&  \leq\left\vert \xi_{k}\right\vert +C\sum_{i=1}^{N}\left\{  \prod_{i\leq
\ell\leq N}\left(  1-\left\vert z_{m_{\ell}-1}\right\vert ^{2}\right)
^{\sigma}\right\}  \sum_{m_{i-1}\leq j\leq m_{i}}\left\vert \xi_{j-1}%
\right\vert +C\sum_{m_{N}\leq j\leq k}\left\vert \xi_{j-1}\right\vert
\nonumber
\end{align}
for $0\leq k\leq\ell$.

We now return our attention to the tree $\mathcal{Y}$. For each $\alpha
\in\mathcal{Y}$, with corresponding index $j\in\left\{  j\right\}  _{j=1}^{J}%
$, there are values $a\left(  \alpha\right)  =a_{j}$, $\xi\left(
\alpha\right)  =\xi_{j}$ and $m\left(  \alpha\right)  =z_{m\left(  j\right)
}\in\mathcal{Y}$. Define functions $f\left(  \alpha\right)  =\left|  a\left(
\alpha\right)  \right|  $ and $g\left(  \alpha\right)  =\left|  \xi\left(
\alpha\right)  \right|  $ on the tree $\mathcal{Y}$. Note that we are simply
relabelling the indices $\left\{  j\right\}  _{j=1}^{J}$ as $\alpha
\in\mathcal{Y}$ to emphasize the tree structure of $\mathcal{Y}$ when
convenient. If we define operators
\begin{align*}
J_{k}g\left(  \alpha\right)   &  =\sum_{m_{k-1}\left(  \alpha\right)
\leq\beta\leq m_{k}\left(  \alpha\right)  }g\left(  A\beta\right)  ,\\
J_{\infty}g\left(  \alpha\right)   &  =g\left(  \alpha\right)  +\sum
_{m_{N\left(  \alpha\right)  }\left(  \alpha\right)  \leq\beta\leq\alpha
}g\left(  A\beta\right)
\end{align*}
on the tree $\mathcal{Y}$, then inequality (\ref{consequence}) implies in
particular that
\begin{equation}
f\left(  \alpha\right)  \leq C\left(  J_{\infty}g\left(  \alpha\right)
+\sum_{k=1}^{N\left(  \alpha\right)  }\left\{  \prod_{k\leq\ell\leq N\left(
\alpha\right)  }\left(  1-\left|  z_{m_{\ell}-1}\right|  ^{2}\right)
^{\sigma}\right\}  J_{k}g\left(  \alpha\right)  \right)  ,\;\;\;\;\;\alpha
\in\mathcal{Y}. \label{fg}%
\end{equation}

Recall that we are assuming that the measure $d\mu=\sum_{\alpha\in\mathcal{Y}%
}\left(  \log\frac{1}{1-\left\vert z_{\alpha}\right\vert ^{2}}\right)  ^{-1}$,
where $z_{\alpha}=z_{j}\in\mathbb{D}$ if $\alpha$ corresponds to $j$,
satisfies the weak simple condition,
\begin{equation}
\beta\left(  0,t\right)  \sum_{j:z_{j}\in S\left(  t\right)  \text{ is
minimal}}\mu\left(  j\right)  \leq C,\;\;\;\;\;t\in\mathcal{T}. \label{simp}%
\end{equation}
Note that this last inequality refers to the tree $\mathcal{T}$ rather than to
$\mathcal{Y}$. Using the fact that $\beta\left(  0,\alpha\right)  \approx
\log\frac{1}{1-\left\vert z_{\alpha}\right\vert ^{2}}$, we obtain from this
weak simple condition that if $S\left(  t\right)  \approx V_{z_{k}}$, i.e.
$t\approx\left[  1-\left(  1-\left\vert z_{k}\right\vert ^{\beta}\right)
\right]  z_{k}$, then
\begin{align*}
\sum_{j:z_{j}\in S\left(  t\right)  \text{ is minimal}}\mu\left(  j\right)
&  \leq C\beta\left(  0,t\right)  ^{-1}\approx C\left(  \log\frac{1}{\left(
1-\left\vert z_{k}\right\vert ^{2}\right)  ^{\beta}}\right)  ^{-1}\\
&  \approx C\left(  \log\frac{1}{1-\left\vert z_{k}\right\vert ^{2}}\right)
^{-1}=C\mu\left(  z_{k}\right)  ,
\end{align*}
by the definition of the region $V_{z_{k}}$. To utilize this inequality on the
tree $\mathcal{Y}$ we need the following crucial property of the sequence $Z$:
if $\left[  \alpha,\beta\right]  $ is a geodesic in $\mathcal{Y}$ such that
$\gamma\notin\Gamma$ for all $\alpha<\gamma\leq\beta$, then the geodesic
$\left[  \alpha,\beta\right]  $, considered as a set of points in the tree
$\mathcal{T}$, is \emph{scattered in }$\mathcal{T}$ in the sense that no two
distinct points $\gamma,\gamma^{\prime}\in\left[  \alpha,\beta\right]  $ are
comparable in $\mathcal{T}$, i.e. neither $\gamma\leq\gamma^{\prime}$ nor
$\gamma^{\prime}\leq\gamma$ in $\mathcal{T}$. With this observation we obtain
that on the tree $\mathcal{Y}$, the adjoint $J_{k}^{\ast}$ of $J_{k}$
satisfies
\begin{equation}
J_{k}^{\ast}\mu\left(  \alpha\right)  \leq C\mu\left(  \alpha\right)
,\;\;\;\;\;\alpha\in\mathcal{Y}. \label{Istar'}%
\end{equation}
Now (\ref{control}) will follow from (\ref{fg}) together with the inequality
\begin{equation}
\sum_{\alpha\in\mathcal{Y}}J_{k}g\left(  \alpha\right)  ^{2}\mu\left(
\alpha\right)  \leq C\sum_{\alpha\in\mathcal{Y}}g\left(  \alpha\right)
^{2}\mu\left(  \alpha\right)  ,\;\;\;\;\;g\geq0, \label{doublemu'}%
\end{equation}
uniformly in $k$, and thus it suffices to show the equivalence of
(\ref{doublemu'}) and (\ref{Istar'}).

To see this we first claim that the inequality
\begin{equation}
\sum_{\alpha\in\mathcal{Y}}Ig\left(  \alpha\right)  ^{2}\mu\left(
\alpha\right)  \leq C\sum_{\alpha\in\mathcal{Y}}g\left(  \alpha\right)
^{2}\mu\left(  \alpha\right)  , \label{doublemu}%
\end{equation}
is equivalent to
\begin{equation}
I^{\ast}\mu\left(  \alpha\right)  \leq C\mu\left(  \alpha\right)
,\;\;\;\;\;\alpha\in\mathcal{Y}. \label{Istar}%
\end{equation}
Indeed, (\ref{Istar}) is obviously necessary for (\ref{doublemu}). To see the
converse, we use our more general tree theorem, Theorem 3 of \cite{ArRoSa},
for the tree $\mathcal{Y}$:
\[
\sum_{\alpha\in\mathcal{Y}}Ig\left(  \alpha\right)  ^{2}w\left(
\alpha\right)  \leq C\sum_{\alpha\in\mathcal{Y}}g\left(  \alpha\right)
^{2}v\left(  \alpha\right)  ,\;\;\;\;\;g\geq0,
\]
if and only if
\begin{equation}
\sum_{\beta\geq\alpha}I^{\ast}w\left(  \beta\right)  ^{2}v\left(
\beta\right)  ^{-1}\leq CI^{\ast}w\left(  \alpha\right)  <\infty
,\;\;\;\;\;\alpha\in\mathcal{Y}. \label{twoweightcond}%
\end{equation}
With $w=v=\mu$, (\ref{Istar}) yields condition (\ref{twoweightcond}) as
follows:
\[
\sum_{\beta\geq\alpha}I^{\ast}\mu\left(  \beta\right)  ^{2}\mu\left(
\beta\right)  ^{-1}\leq C\sum_{\beta\geq\alpha}\mu\left(  \beta\right)
^{2}\mu\left(  \beta\right)  ^{-1}=C\sum_{\beta\geq\alpha}\mu\left(
\beta\right)  =CI^{\ast}\mu\left(  \alpha\right)  ,
\]
and this completes the proof of the claim.

In general, condition (\ref{Istar'}), a consequence of the weak simple
condition, does not imply the simple condition (\ref{Istar}). However, we can
again exploit the crucial property of the sequence $Z$ mentioned above -
namely that if $\left[  \alpha,\beta\right]  $ is a geodesic in $\mathcal{Y}$
with $\left(  \alpha,\beta\right]  \cap\Gamma=\phi$, then $\left[
\alpha,\beta\right]  $, considered as a set of points in the tree
$\mathcal{T}$, is \emph{scattered in }$\mathcal{T}$. Now decompose the tree
$\mathcal{Y}$ into a family of pairwise disjoint forests $\mathcal{Y}_{\ell}$
as follows. Let $\mathcal{Y}_{1}$ consist of the root $o$ of $\mathcal{Y}$
together with all points $\beta>o$ having $\gamma\notin\Gamma$ for
$o<\gamma\leq\beta$. Then let $\mathcal{Y}_{2}$ consist of each minimal point
$\alpha$ in $\mathcal{Y}\setminus\mathcal{Y}_{1}$ together with all points
$\beta>\alpha$ having $\gamma\notin\Gamma$ for $\alpha<\gamma\leq\beta$, then
let $\mathcal{Y}_{3}$ consist of each minimal point $\alpha$ in $\mathcal{Y}%
\setminus\left(  \mathcal{Y}_{1}\cup\mathcal{Y}_{2}\right)  $ together with
all points $\beta>\alpha$ having $\gamma\notin\Gamma$ for $\alpha<\gamma
\leq\beta$, etc.

A key property of this decomposition is that on $\mathcal{Y}_{\ell}$ the
operator $J_{k}$ sees only the values of $g$ on $\mathcal{Y}_{\ell}$ itself. A
second key property is that since the geodesics in $\mathcal{Y}_{\ell}$ are
scattered, we see that the restriction $\mu_{\ell}$ of $\mu$ to the forest
$\mathcal{Y}_{\ell}$ satisfies the \emph{simple} condition, rather than just
the \emph{weak simple} condition. As a consequence, upon decomposing each
forest $\mathcal{Y}_{\ell}$ into trees and applying the above claim with
$\mu_{\ell}$ in place of $\mu$, i.e. (\ref{doublemu}) holds if and only if
(\ref{Istar}) holds, we conclude that
\[
\sum_{\alpha\in\mathcal{Y}_{\ell}}J_{k}g\left(  \alpha\right)  ^{2}\mu\left(
\alpha\right)  \leq C\sum_{\alpha\in\mathcal{Y}_{\ell}}g\left(  \alpha\right)
^{2}\mu\left(  \alpha\right)  ,\;\;\;\;\;g\geq0,
\]
uniformly in $k$ for each $\ell\geq1$. Summing in $\ell$ and using the finite
overlap, we obtain the sufficiency of (\ref{Istar'}) for (\ref{doublemu'}).

Finally, to see that (\ref{control}) now follows from (\ref{fg}), we use that
\[
\sum_{i=1}^{N}\left\{  \prod_{i\leq\ell\leq N}\left(  1-\left\vert z_{m_{\ell
}-1}\right\vert ^{2}\right)  ^{\sigma}\right\}  \leq C
\]
in (\ref{fg}) to obtain (\ref{control}). This completes the proof of Lemma
\ref{weightcon}.

\bigskip

Now we prove that the function $\varphi=\sum_{i=1}^{J}a_{i}\varphi_{i}$
constructed above comes close to interpolating the data $\left\{  \xi
_{j}\right\}  _{j=1}^{J}$ provided we choose $\varepsilon>0$ sufficiently
small in (\ref{logeps''}).

\begin{lemma}
\label{inductivecon}Suppose $s>-1$, that $\left\{  \xi_{j}\right\}  _{j=1}%
^{J}$ is a sequence of complex numbers, and let $0<\delta<1$. Let $\varphi
_{j}$, $g_{j}$ and $\gamma_{j}$ correspond to $z_{j}$ as in Lemma
\ref{analyticcon} and with the same $s$. Then for $\varepsilon>0$ sufficiently
small in (\ref{logeps''}), there is $\left\{  a_{i}\right\}  _{i=1}^{J}$ such
that $\varphi=\sum_{i=1}^{J}a_{i}\varphi_{i}$ satisfies
\begin{equation}
\left\Vert \left\{  \xi_{j}-\varphi\left(  z_{j}\right)  \right\}  _{j=1}%
^{J}\right\Vert _{\ell^{2}\left(  \mu\right)  }<\delta\left\Vert \left\{
\xi_{j}\right\}  _{j=1}^{J}\right\Vert _{\ell^{2}\left(  \mu\right)  }
\label{app}%
\end{equation}
and
\begin{equation}
\left\Vert \varphi\right\Vert _{B_{2}}\leq C\left\Vert \left\{  a_{j}\right\}
_{j=1}^{J}\right\Vert _{\ell^{2}\left(  \mu\right)  }. \label{bounds}%
\end{equation}

\end{lemma}

\begin{remark}
\label{linear}The construction in the proof below shows that both the sequence
$\left\{  a_{i}\right\}  _{i=1}^{J}$ and the function $\varphi$ depend
linearly on the data $\left\{  \xi_{j}\right\}  _{j=1}^{J}$.
\end{remark}

\textbf{Proof}: We now show that both (\ref{app}) and (\ref{bounds}) hold for
the function $\varphi=\sum_{i=1}^{J}a_{i}\varphi_{i}$ constructed above. Fix
an index $\ell\in\mathbb{N}$, and with notation as above, let $\mathcal{F}%
_{\ell}=\mathbb{N}\setminus\mathcal{G}_{\ell}$ and write using (\ref{amdef}),
\begin{align}
\varphi\left(  z_{\ell}\right)  -\xi_{\ell}  &  =\sum_{i=1}^{\infty}%
a_{i}\varphi_{i}\left(  z_{\ell}\right)  -\xi_{\ell}\label{match}\\
&  =\left(  \sum_{i\in\mathcal{G}_{P\left(  \ell\right)  }}a_{i}\varphi
_{i}\left(  z_{\ell}\right)  +a_{\ell}\varphi_{\ell}\left(  z_{\ell}\right)
+\sum_{i\in\mathcal{F}_{\ell}}a_{i}\varphi_{i}\left(  z_{\ell}\right)
\right)  -\left(  a_{\ell}+\sum_{i\in\mathcal{G}_{\ell}\setminus\left\{
0\right\}  }\beta_{\ell,i}a_{P\left(  i\right)  }\right) \nonumber\\
&  =\sum_{i\in\mathcal{G}_{\ell}\setminus\left\{  0\right\}  }a_{P\left(
i\right)  }\left(  \varphi_{P\left(  i\right)  }\left(  z_{\ell}\right)
-\beta_{\ell,i}\right)  +\sum_{i\in\mathcal{F}_{\ell}}a_{i}\varphi_{i}\left(
z_{\ell}\right) \nonumber\\
&  =\sum_{i\in\mathcal{F}_{\ell}}a_{i}\varphi_{i}\left(  z_{\ell}\right)
\equiv B_{\ell},\nonumber
\end{align}
since $\varphi_{\ell}\left(  z_{\ell}\right)  =1$ and $\varphi_{P\left(
i\right)  }\left(  z_{\ell}\right)  =\beta_{\ell,i}$.

We now claim that
\begin{equation}
\left\vert B_{\ell}\right\vert \leq C\sum_{i=1}^{J}\left\vert a_{i}\right\vert
\mu\left(  z_{i}\right)  . \label{logest'}%
\end{equation}
We first note that if $z_{\ell}\notin V_{z_{i}}$, then
\begin{equation}
\left\vert \varphi_{i}\left(  z_{\ell}\right)  \right\vert \leq C\left(
1-\left\vert z_{i}\right\vert ^{2}\right)  ^{\sigma},\;\;\;\;\;\sigma>0,
\label{better'}%
\end{equation}
by the third line in (\ref{satisfies}). On the other hand, if $z_{\ell}\in
V_{z_{i}}$, then $\left\vert z_{i}\right\vert <\left\vert z_{\ell}\right\vert
$, and if $\mathcal{G}_{\ell}=\left[  k_{0},k_{1},...,k_{m-1},k_{m}\right]  $,
then either $\left\vert z_{i}\right\vert <\left\vert z_{k_{0}}\right\vert $ or
there is $j$ such that $\left\vert z_{k_{j-1}}\right\vert <\left\vert
z_{i}\right\vert \leq\left\vert z_{k_{j}}\right\vert $. Note however that
equality cannot hold here by Lemma \ref{Mars}, and so we actually have
$\left\vert z_{k_{j-1}}\right\vert <\left\vert z_{i}\right\vert <\left\vert
z_{k_{j}}\right\vert $. From (\ref{daughter}) we obtain that no index
$m\in\left(  k_{j-1},k_{j}\right)  $ satisfies $V_{z_{k_{j}}}\subset V_{z_{m}%
}$. Since $i\notin\mathcal{G}_{\ell}$, we have $i\in\left(  k_{j-1}%
,k_{j}\right)  $ and thus we have both
\[
V_{z_{k_{j}}}\varsubsetneq V_{z_{i}}\text{ and }\left\vert z_{k_{j}%
}\right\vert >\left\vert z_{i}\right\vert .
\]
Now using Lemma \ref{Mars} and $\beta\eta>1$, we obtain
\[
\left(  1-\left\vert z_{k_{j}}\right\vert ^{2}\right)  ^{\beta}\leq\left(
1-\left\vert z_{i}\right\vert ^{2}\right)  ^{\beta\eta}\ll\left(  1-\left\vert
z_{i}\right\vert ^{2}\right)  .
\]
If we choose $w\in V_{z_{k_{j}}}\setminus V_{z_{i}}$, then $w,z_{\ell}\in
V_{z_{k_{j}}}$ implies $\left\vert z_{\ell}-w\right\vert \leq C\left(
1-\left\vert z_{k_{j}}\right\vert ^{2}\right)  ^{\beta}$ by definition, and
$w\notin V_{z_{i}}$ implies $\left\vert 1-\overline{w}\cdot Pz_{i}\right\vert
\geq C\left(  1-\left\vert z_{i}\right\vert ^{2}\right)  ^{\beta}$. Together
with the reverse triangle inequality we thus have
\begin{align*}
\left\vert 1-\overline{z_{\ell}}\cdot Pz_{i}\right\vert  &  \geq\left\vert
1-\overline{w}\cdot Pz_{i}\right\vert -\left\vert \overline{z_{\ell}}\cdot
Pz_{i}-\overline{w}\cdot Pz_{i}\right\vert \\
&  \geq C\left(  1-\left\vert z_{i}\right\vert ^{2}\right)  ^{\beta}-C\left(
1-\left\vert z_{i}\right\vert ^{2}\right)  ^{\beta\eta}\\
&  \geq\left(  1-\left\vert z_{i}\right\vert \right)  ^{\beta_{1}},
\end{align*}
for some $\beta_{1}\in\left(  \beta,\rho\right)  $ (again provided the
$\left\vert z_{i}\right\vert $ are large enough). Thus in the case $z_{\ell
}\in V_{z_{i}}$, estimate (\ref{better'}) again follows from the third line in
(\ref{satisfies}). Finally, the estimate $\left(  1-\left\vert z_{i}%
\right\vert ^{2}\right)  ^{\sigma}\leq C\mu\left(  z_{i}\right)  $ is trivial
and this yields (\ref{logest'}).

Combining (\ref{logest'}) and (\ref{logeps''}) we then have for the sequence
$\left\{  \xi_{j}-\varphi\left(  z_{j}\right)  \right\}  _{j=1}^{J}$,
\begin{align*}
\left\Vert \left\{  \xi_{j}-\varphi\left(  z_{j}\right)  \right\}  _{j=1}%
^{J}\right\Vert _{\ell^{2}\left(  d\mu\right)  }  &  \leq C\left\Vert \left\{
B_{j}\right\}  _{j=1}^{J}\right\Vert _{\ell^{2}\left(  d\mu\right)  }\\
&  \leq C\sum_{i=1}^{J}\left\vert a_{i}\right\vert \mu\left(  z_{i}\right)
\left\{  \sum_{j=1}^{J}\mu\left(  z_{j}\right)  \right\}  ^{\frac{1}{2}}\\
&  \leq C\left\{  \sum_{i=1}^{J}\left\vert a_{i}\right\vert ^{2}\mu\left(
z_{i}\right)  \right\}  ^{\frac{1}{2}}\left\{  \sum_{i=1}^{J}\mu\left(
z_{i}\right)  \right\}  ^{\frac{1}{2}}\left\{  \sum_{j=1}^{J}\mu\left(
z_{j}\right)  \right\}  ^{\frac{1}{2}}\\
&  =C\left\Vert \mu\right\Vert \left\Vert \left\{  a_{j}\right\}  _{j=1}%
^{J}\right\Vert _{\ell^{2}\left(  \mu\right)  }\\
&  <C\varepsilon\left\Vert \left\{  a_{j}\right\}  _{j=1}^{J}\right\Vert
_{\ell^{2}\left(  \mu\right)  }.
\end{align*}
This completes the proof of (\ref{app}).

We now prove the estimate $\left\Vert \varphi\right\Vert _{B_{2}}\leq C$ in
(\ref{bounds}) whenever $\left\Vert \left\{  a_{j}\right\}  _{j=1}%
^{J}\right\Vert _{\ell^{2}\left(  d\mu\right)  }=1$, independent of $J\geq1$.
Thus we must show that
\[
\int_{\mathbb{D}}\left\vert \nabla\varphi\left(  z\right)  \right\vert
^{2}dz\leq C,
\]
independent of $J\geq1$. Now
\[
\varphi=\sum_{i=1}^{J}a_{i}\varphi_{i}=\sum_{i=1}^{J}a_{i}\Gamma_{s}%
g_{i}=\Gamma_{s}g
\]
where $g=\sum_{i=1}^{J}a_{i}g_{i}$ with $\left\Vert \left\{  a_{i}\right\}
_{i=1}^{J}\right\Vert _{\ell^{2}\left(  \mu\right)  }=1$. Moreover,
\[
\left\vert \nabla\Gamma_{s}g\left(  z\right)  \right\vert \leq C_{s}^{\prime
}\widehat{T_{s}}\left\vert g\right\vert \left(  z\right)
\]
where the operator $\widehat{T_{s}}$ is given by
\[
\widehat{T_{s}}f\left(  z\right)  =C_{s}\int_{\mathbb{D}}\frac{f\left(
w\right)  \left(  1-\left\vert w\right\vert ^{2}\right)  ^{s}}{\left\vert
1-\overline{w}\cdot z\right\vert ^{1+s}}dw.
\]
Thus we must estimate $\int_{\mathbb{D}}\left\vert \widehat{T_{s}}\left\vert
g\right\vert \left(  z\right)  \right\vert ^{2}dz$. \ Now by Theorem 2.10 in
\cite{Zhu}, $\widehat{T_{s}}$ is bounded on $L^{2}$ if and only if
$s>-\frac{1}{2}$ so, for any such $s$
\[
\int_{\mathbb{D}}\left\vert \nabla\varphi\left(  z\right)  \right\vert
^{2}dz\leq C\int_{\mathbb{D}}\left\vert g\left(  z\right)  \right\vert ^{2}dz
\]
Since the supports of the $g_{i}$ are pairwise disjoint by the separation
condition, we obtain from (\ref{maincon'}) that $g=\sum_{i=1}^{J}a_{i}g_{i}$
satisfies
\begin{align*}
\int_{\mathbb{D}}\left\vert g\left(  z\right)  \right\vert ^{2}dz  &
=\sum_{j=1}^{J}\left\vert a_{i}\right\vert ^{2}\int_{\mathbb{D}}\left\vert
g_{i}\left(  z\right)  \right\vert ^{2}dz\\
&  \leq C\sum_{j=1}^{J}\left\vert a_{i}\right\vert ^{2}\left(  \log\frac
{1}{1-\left\vert z_{j}\right\vert ^{2}}\right)  ^{-1}\\
&  =C\left\Vert \left\{  a_{j}\right\}  _{j=1}^{J}\right\Vert _{\ell
^{2}\left(  d\mu\right)  }=C.
\end{align*}
This completes the proof of Lemma \ref{inductivecon}.

\bigskip

Now we finish proving the sufficiency portion of Theorem \ref{Boeonto}. Fix
$s>-1$, $0<\delta<1$ and $\left\{  \xi_{j}\right\}  _{j=1}^{J}$ with
$\left\Vert \left\{  \xi_{j}\right\}  _{j=1}^{J}\right\Vert _{\ell^{2}\left(
\mu\right)  }=1$. Then by Lemma \ref{inductivecon} there is $f_{1}=\sum
_{i=1}^{J}a_{i}^{1}\varphi_{i}\in B_{2}$ such that $\left\Vert \left\{
\xi_{j}-f_{1}\left(  z_{j}\right)  \right\}  _{j=1}^{J}\right\Vert _{\ell
^{2}\left(  \mu\right)  }<\delta$ and using Lemma \ref{weightcon} as well,
$\left\Vert \left\{  a_{i}^{1}\right\}  _{i=1}^{J}\right\Vert _{\ell
^{2}\left(  \mu\right)  },\left\Vert f_{1}\right\Vert _{B_{2}}\leq C$ where $C
$ is the product of the constants in (\ref{bounds}) and (\ref{control}). Now
apply Lemma \ref{inductivecon} to the sequence $\left\{  \xi_{j}-f_{1}\left(
z_{j}\right)  \right\}  _{j=1}^{J}$ to obtain the existence of $f_{2}%
=\sum_{i=1}^{\infty}a_{i}^{2}\varphi_{i}\in B_{2}$ such that $\left\Vert
\left\{  \xi_{j}-f_{1}\left(  z_{j}\right)  -f_{2}\left(  z_{j}\right)
\right\}  _{j=1}^{J}\right\Vert _{\ell^{2}\left(  \mu\right)  }<\delta^{2}$
and again using Lemma \ref{weightcon} as well, $\left\Vert \left\{  a_{i}%
^{2}\right\}  _{i=1}^{J}\right\Vert _{\ell^{2}\left(  \mu\right)  },\left\Vert
f_{2}\right\Vert _{B_{2}}\leq C\delta$ where $C$ is again the product of the
constants in (\ref{bounds}) and (\ref{control}). Continuing inductively, we
obtain $f_{m}=\sum_{i=1}^{J}a_{i}^{m}\varphi_{i}\in B_{2}$ such that
\begin{align*}
\left\Vert \left\{  \xi_{j}-\sum_{i=1}^{m}f_{i}\left(  z_{j}\right)  \right\}
_{j=1}^{J}\right\Vert _{\ell^{\infty}}  &  <\delta^{m},\\
\left\Vert \left\{  a_{i}^{m}\right\}  _{i=1}^{J}\right\Vert _{\ell^{2}\left(
\mu\right)  },\left\Vert f_{m}\right\Vert _{B_{2}}  &  \leq C\delta^{m-1}.
\end{align*}
If we now take%
\[
\varphi=\sum_{m=1}^{\infty}f_{m}=\sum_{m=1}^{\infty}\left\{  \sum_{i=1}%
^{J}a_{i}^{m}\varphi_{i}\right\}  =\sum_{i=1}^{J}a_{i}\varphi_{i},
\]
we have
\begin{align}
\xi_{j}  &  =\varphi\left(  z_{j}\right)  ,\;\;\;\;\;1\leq j\leq
J,\label{constructphi}\\
\left\Vert \left\{  a_{i}\right\}  _{i=1}^{J}\right\Vert _{\ell^{2}\left(
\mu\right)  }  &  \leq C,\nonumber\\
\left\Vert \varphi\right\Vert _{B_{2}}  &  \leq C,\nonumber
\end{align}
if $\varepsilon>0$ is chosen small enough in (\ref{logeps''}). A limiting
argument using $J\rightarrow\infty$ now completes the sufficiency proof of
Theorem \ref{Boeonto}.

\section{Theorem B}

We follow the pattern of the previous section. We will state a result more
general than Theorem B and prove the half of that result that contains Theorem
B. We return to the other half of the proof in Section \ref{neces}

In this section we drop the requirement that $\left\Vert \mu_{Z}\right\Vert
<\infty.$ \ The key is to let $Z$ be a subtree of $\mathcal{T}$, so that if
$z\in Z$ is a B\"{o}e child of $w\in Z$, then $z$ actually lies in the Bergman
successor set $S\left(  w\right)  $ of $w$, and hence the value of
$c_{\rho,\alpha}\left(  \gamma_{w}\left(  z\right)  \right)  $ in Lemma
\ref{analyticcon} is $1$, which is exploited in (\ref{dyadicexploit}) below.
The advantage when assuming (\ref{special}) is that we may dispense with the
complicated inductive definition of the coefficients $a_{k}$ in (\ref{amdef})
for the holomorphic function $\mathcal{S}\xi$ in (\ref{defofM}) approximating
$\xi$ on $Z$, and instead use the elementary construction in (\ref{elemconst})
below of a holomorphic function $M\xi$ approximating the \emph{integrated}
sequence $\left(  I\xi\right)  _{j}=\sum_{z_{i}\leq z_{j}}\xi_{i}$ on $Z$.
This permits us to interpolate the \emph{difference} sequence $\bigtriangleup
\xi$ using the operator $\bigtriangleup M$, whose kernel is better localized.
Of course in the absence of (\ref{special}), the values $c_{\rho,\alpha
}\left(  \gamma_{w}\left(  z\right)  \right)  $ may lie in $\left[
0,1\right)  $ and then $M\xi$ will \emph{not} be a good approximation to
$I\xi$ on $Z$.

\begin{theorem}
\label{onto}Suppose $Z=\left\{  z_{j}\right\}  _{j=1}^{\infty}\subset
\mathbb{D}$ is a subtree of $\mathcal{T}$ that satisfies the separation
condition (\ref{sep}), and if $C$ is the constant in (\ref{sep}), that there
is $\beta\in\left(  1-C/2,1\right)  $ satisfying (\ref{special}). Then $Z$ is
onto interpolating for the B\"{o}e space $B_{2,Z}$ if and only if the weak
simple condition (\ref{weaksimple}) holds.
\end{theorem}

\textbf{Proof }(of the sufficiency of the conditions in the Theorem): \ Fix
$\left\{  \xi_{j}\right\}  _{j=1}^{\infty}$ with
\[
\left\Vert \left\{  \frac{\xi_{j}}{\left\Vert k_{z_{j}}\right\Vert _{B_{2}}%
}\right\}  _{j=1}^{\infty}\right\Vert _{\ell^{2}}=1.
\]
Recall that $\left\Vert k_{z_{j}}\right\Vert _{B_{2}}\approx\left(  \log
\frac{1}{1-\left\vert z_{j}\right\vert ^{2}}\right)  ^{\frac{1}{2}}$ and that
we may suppose $Z\subset\mathcal{T}$. We note that (\ref{sigmasum}) holds here
- in fact the proof is simpler using the separation condition (\ref{sep}) and
the assumption that $Z$ is a subtree of $\mathcal{T}$ (and hence has branching
number at most $2$). We can adjoin the origin to $Z$ in which case
(\ref{sigmasum}) yields that $\sum_{j=1}^{\infty}\left(  1-\left\vert
z_{j}\right\vert \right)  ^{\sigma}<\infty$. Thus, as at the start of the
previous proof, given any $\sigma>0$, we can discard all points from $Z$ that
lie in some ball $B\left(  0,R\right)  $, $R<1$, and reorder the remaining
points so that
\begin{equation}
\left(  \log\frac{1}{1-R^{2}}\right)  ^{-1},\;1-R^{2},\;\sum_{j=1}^{\infty
}\left(  1-\left\vert z_{j}\right\vert \right)  ^{\sigma}<\varepsilon.
\label{logeps'}%
\end{equation}
We next, in addition, suppose that the sequence $Z=\left\{  z_{j}\right\}
_{j=1}^{J}$ is finite, and obtain an appropriate estimate independent of
$J\geq1$. Given a sequence of complex numbers $\xi=\left\{  \xi_{j}\right\}
_{j=1}^{J}$ we define a holomorphic function $M\xi$ on the ball by
\begin{equation}
M\xi\left(  z\right)  =\sum_{j=1}^{J}\xi_{j}\varphi_{z_{j}}\left(  z\right)
,\;\;\;\;\;z\in\mathbb{D}, \label{elemconst}%
\end{equation}
where $\varphi_{w}\left(  z\right)  $ is as in Lemma \ref{analyticcon}. View
$\mu$ as the measure assigning mass $\left(  \log\frac{1}{1-\left\vert
z_{j}\right\vert ^{2}}\right)  ^{-1}$ to the point $j\in\left\{
0,1,2,...,J\right\}  $. We have
\[
\left\Vert \left\{  \xi_{j}\right\}  _{j=1}^{J}\right\Vert _{\ell^{2}\left(
d\mu\right)  }\approx\left\Vert \left\{  \frac{\xi_{j}}{\left\Vert k_{z_{j}%
}\right\Vert _{B_{2}}}\right\}  _{j=1}^{J}\right\Vert _{\ell^{2}},
\]
for any complex sequence $\left\{  \xi_{j}\right\}  _{j=1}^{J}$. We will use
another useful consequence of Lemma \ref{Mars}:%
\begin{equation}
1-\left\vert Az_{j}\right\vert ^{2}\leq\left(  1-\left\vert z_{\ell
}\right\vert ^{2}\right)  ^{\eta},\ \ \ \ \ \text{for }z_{j}\in V_{z_{\ell}%
}\setminus\mathcal{C}\left(  z_{\ell}\right)  . \label{exclusion}%
\end{equation}
Indeed, if $z_{j}\in V_{z_{\ell}}\setminus\mathcal{C}\left(  z_{\ell}\right)
$, then $Az_{j}\neq z_{\ell}$ and $\left\vert Az_{j}\right\vert \geq\left\vert
z_{\ell}\right\vert $ by the construction in (\ref{daughter}). Then
$V_{z_{\ell}}\cap V_{Az_{j}}$ contains $z_{j}$ and is thus nonempty, and Lemma
\ref{Mars} now shows that $1-\left\vert Az_{j}\right\vert ^{2}\leq\left(
1-\left\vert z_{\ell}\right\vert ^{2}\right)  ^{\eta}$.

Now define a linear map $T$ from $\ell^{2}\left(  d\mu\right)  $ to $\ell
^{2}\left(  d\mu\right)  $ by
\[
T\xi=\bigtriangleup\left(  M\xi\right)  \mid_{Z}=\left\{  M\xi\left(
z_{k}\right)  -M\xi\left(  Az_{k}\right)  \right\}  _{j=1}^{J}=\left\{
\sum_{j=1}^{J}\xi_{j}\left[  \varphi_{z_{j}}\left(  z_{k}\right)
-\varphi_{z_{j}}\left(  Az_{k}\right)  \right]  \right\}  _{j=1}^{J},
\]
where $Az_{j}$ denotes the predecessor of $z_{j}$ in the forest structure on
$Z$ defined in (\ref{daughter}) above (we identify $z_{k}$ with $k$ here). Let
$\mathcal{R}$ denote the set of all roots of maximal trees in the forest. In
the event that $z_{k}\in\mathcal{R}$, then $Az_{k}$ isn't defined and our
convention is to define $\varphi_{z_{j}}\left(  Az_{k}\right)  =0$. We claim
that $T$ is a bounded invertible map on $\ell^{2}\left(  d\mu\right)  $ with
norms independent of $J\geq1$. To see this it is enough to prove that
$\mathbb{I}-T$ has small norm on $\ell^{2}\left(  d\mu\right)  $ where
$\mathbb{I}$ denotes the identity operator. We have
\begin{align*}
\left(  \mathbb{I}-T\right)  \xi &  =\left\{  \xi_{k}-\sum_{j=1}^{J}\xi
_{j}\left[  \varphi_{z_{j}}\left(  z_{k}\right)  -\varphi_{z_{j}}\left(
Az_{k}\right)  \right]  \right\}  _{k=1}^{J}\\
&  =\left\{  \xi_{k}\varphi_{z_{k}}\left(  Az_{k}\right)  \right\}  _{k=1}%
^{J}-\left\{  \sum_{j:j\neq k}\xi_{j}\left[  \varphi_{z_{j}}\left(
z_{k}\right)  -\varphi_{z_{j}}\left(  Az_{k}\right)  \right]  \right\}
_{k=1}^{J}%
\end{align*}
since $\varphi_{z_{k}}\left(  z_{k}\right)  =1$.

Now we estimate the kernel $K\left(  k,j\right)  $ of the operator
$\mathbb{I}-T$. We have on the diagonal,
\[
\left\vert K\left(  k,k\right)  \right\vert =\left\{
\begin{array}
[c]{llll}%
\left\vert \varphi_{z_{k}}\left(  Az_{k}\right)  \right\vert  & \leq & \left(
1-\left\vert z_{k}\right\vert ^{2}\right)  ^{\left(  \rho-\beta_{1}\right)
\left(  1+s\right)  } & \text{ if }z_{k}\notin\mathcal{R}\\
0 &  &  & \text{ if }z_{k}\in\mathcal{R}%
\end{array}
\right.  ,
\]
by the third estimate in (\ref{satisfies}).

Suppose now that $z_{k}\notin\mathcal{R}$ and $j\neq k$. Lemma \ref{fax} shows
that $\left\vert \varphi_{z_{j}}^{\prime}\left(  \zeta_{k}\right)  \right\vert
\leq\left(  1-\left\vert z_{j}\right\vert ^{2}\right)  ^{-\alpha} $ and the
definition of $V_{z_{j}}$ shows that $\left\vert z_{k}-Az_{k}\right\vert
\leq\left(  1-\left\vert Az_{k}\right\vert ^{2}\right)  ^{\beta}$. Thus if
$1-\left\vert Az_{k}\right\vert ^{2}\leq\left(  1-\left\vert z_{j}\right\vert
^{2}\right)  ^{\eta}$, then
\begin{align*}
\left\vert K\left(  k,j\right)  \right\vert  &  =\left\vert \varphi_{z_{j}%
}\left(  z_{k}\right)  -\varphi_{z_{j}}\left(  Az_{k}\right)  \right\vert
\leq\left\vert \varphi_{z_{j}}^{\prime}\left(  \zeta_{k}\right)  \right\vert
\left\vert z_{k}-Az_{k}\right\vert \\
&  \leq C\left(  1-\left\vert z_{j}\right\vert ^{2}\right)  ^{-\alpha}\left(
1-\left\vert Az_{k}\right\vert ^{2}\right)  ^{\beta}\\
&  \leq C\left(  1-\left\vert z_{j}\right\vert ^{2}\right)  ^{\eta\left(
\beta-\delta\right)  -\alpha}\left(  1-\left\vert Az_{k}\right\vert
^{2}\right)  ^{\delta},
\end{align*}
where the exponent $\eta\left(  \beta-\delta\right)  -\alpha$ is positive if
we choose $\delta$ small enough, since $\alpha<1<\beta\eta$ by Lemma
\ref{Mars}.

Suppose instead that $1-\left\vert Az_{k}\right\vert ^{2}>\left(  1-\left\vert
z_{j}\right\vert ^{2}\right)  ^{\eta}$. Then $Az_{k}\notin V_{z_{j}}$ by Lemma
\ref{Mars}. If $z_{k}\notin\mathcal{C}\left(  z_{j}\right)  $, then
$z_{k}\notin V_{z_{j}}$ by (\ref{exclusion}), and this time we use the third
estimate in (\ref{satisfies}) to obtain
\begin{align*}
\left\vert K\left(  k,j\right)  \right\vert  &  =\left\vert \varphi_{z_{j}%
}\left(  z_{k}\right)  -\varphi_{z_{j}}\left(  Az_{k}\right)  \right\vert
\leq\left\vert \varphi_{z_{j}}\left(  z_{k}\right)  \right\vert +\left\vert
\varphi_{z_{j}}\left(  Az_{k}\right)  \right\vert \\
&  \leq C\left(  1-\left\vert z_{j}\right\vert ^{2}\right)  ^{\left(
\rho-\beta_{1}\right)  \left(  1+s\right)  }\\
&  \leq C\left(  1-\left\vert z_{j}\right\vert ^{2}\right)  ^{\left(
\rho-\beta_{1}\right)  \left(  1+s\right)  -\delta}\left(  1-\left\vert
Az_{k}\right\vert ^{2}\right)  ^{\frac{\delta}{\rho}}.
\end{align*}
On the other hand, if $z_{k}\in\mathcal{C}\left(  z_{j}\right)  $, then
$\left\vert z_{k}\right\vert \geq\left\vert z_{j}\right\vert $ and our
hypothesis (\ref{special}) implies that $z_{k}\in S\left(  z_{j}\right)  $.
Then we have
\begin{equation}
\left\vert K\left(  k,j\right)  \right\vert =\left\vert \varphi_{z_{j}}\left(
z_{k}\right)  -\varphi_{z_{j}}\left(  z_{j}\right)  \right\vert \leq\left(
\log\frac{1}{1-\left\vert z_{j}\right\vert ^{2}}\right)  ^{-1}
\label{dyadicexploit}%
\end{equation}
by the first two estimates in (\ref{satisfies}) since $c_{\rho,\alpha}\left(
\gamma_{z_{j}}\left(  z_{k}\right)  \right)  =1$ in Lemma \ref{analyticcon} if
$z_{k}\in S\left(  z_{j}\right)  $.

Finally, we consider the case when $z_{k}\in\mathcal{R}$ and $j\neq k$. The
third estimate in (\ref{satisfies}) shows that
\[
\left\vert K\left(  k,j\right)  \right\vert =\left\vert \varphi_{z_{j}}\left(
z_{k}\right)  \right\vert \leq C\left(  1-\left\vert z_{j}\right\vert
^{2}\right)  ^{\left(  \rho-\beta_{1}\right)  \left(  1+s\right)  },
\]
where the exponent $\left(  \rho-\beta_{1}\right)  \left(  1+s\right)  $ can
be made as large as we wish by taking $s$ sufficiently large. Combining all
cases we have in particular the following estimate for some $\sigma_{1}%
,\sigma_{2}>0$:
\[
\left\vert K\left(  k,j\right)  \right\vert \leq C\left\{
\begin{array}
[c]{ll}%
\left(  1-\left\vert z_{j}\right\vert ^{2}\right)  ^{\sigma_{1}}\left(
1-\left\vert Az_{k}\right\vert ^{2}\right)  ^{\sigma_{2}}, & \text{ if }%
z_{k}\notin\mathcal{R}\text{ and }z_{k}\notin\mathcal{C}\left(  z_{j}\right)
\\
\left(  \log\frac{1}{1-\left\vert z_{j}\right\vert ^{2}}\right)  ^{-1} &
\text{ if }z_{k}\notin\mathcal{R}\text{ and }z_{k}\in\mathcal{C}\left(
z_{j}\right) \\
\left(  1-\left\vert z_{j}\right\vert ^{2}\right)  ^{3}, & \text{ if }z_{k}%
\in\mathcal{R}%
\end{array}
\right.  .
\]

Now we obtain the boundedness of $\mathbb{I}-T$ on $\ell^{2}\left(
d\mu\right)  $ with small norm by Schur's test. It is here that we use the
assumption that $\mu$ satisfies the weak simple condition (\ref{weaksimple}).
With $\xi\in\ell^{2}\left(  d\mu\right)  $ and $\eta\in\ell^{2}\left(
d\mu\right)  $, we have
\begin{align*}
\left\vert \left\langle \left(  \mathbb{I}-T\right)  \xi,\eta\right\rangle
_{\mu}\right\vert  &  =\left\vert \sum_{k}\left(  \sum_{j}K\left(  k,j\right)
\xi_{j}\right)  \overline{\eta_{k}}\mu\left(  k\right)  \right\vert \\
&  \leq C\sum_{j}\sum_{k\notin\mathcal{R},z_{k}\notin\mathcal{C}\left(
z_{j}\right)  }\left(  1-\left\vert z_{j}\right\vert ^{2}\right)  ^{\sigma
_{1}}\left(  1-\left\vert Az_{k}\right\vert ^{2}\right)  ^{\sigma_{2}%
}\left\vert \xi_{j}\right\vert \left\vert \eta_{k}\right\vert \mu\left(
k\right) \\
&  +C\sum_{j}\sum_{k\notin\mathcal{R},z_{k}\in\mathcal{C}\left(  z_{j}\right)
}\left(  \log\frac{1}{1-\left\vert z_{j}\right\vert ^{2}}\right)
^{-1}\left\vert \xi_{j}\right\vert \left\vert \eta_{k}\right\vert \mu\left(
k\right) \\
&  +C\sum_{j}\sum_{k\in\mathcal{R}}\left(  1-\left\vert z_{j}\right\vert
^{2}\right)  ^{3}\left\vert \xi_{j}\right\vert \left\vert \eta_{k}\right\vert
\mu\left(  k\right)  ,
\end{align*}
and since $\mu\left(  j\right)  \leq\left(  1-\left\vert z_{j}\right\vert
^{2}\right)  ^{\varepsilon}$, we have with $\sigma_{1}^{\prime}=\sigma
_{1}-\varepsilon$,
\begin{align*}
\left\vert \left\langle \left(  \mathbb{I}-T\right)  \xi,\eta\right\rangle
_{\mu}\right\vert  &  \leq C\sum_{j}\sum_{k\notin\mathcal{R},z_{k}%
\notin\mathcal{C}\left(  z_{j}\right)  }\left(  1-\left\vert z_{j}\right\vert
^{2}\right)  ^{\sigma_{1}^{\prime}}\left(  1-\left\vert Az_{k}\right\vert
^{2}\right)  ^{\sigma_{2}}\left\vert \xi_{j}\right\vert \mu\left(  j\right)
\left\vert \eta_{k}\right\vert \mu\left(  k\right) \\
&  +C\sum_{j}\sum_{k\notin\mathcal{R},z_{k}\in\mathcal{C}\left(  z_{j}\right)
}\left\vert \xi_{j}\right\vert \mu\left(  j\right)  \left\vert \eta
_{k}\right\vert \mu\left(  k\right) \\
&  +C\sum_{j}\sum_{k\in\mathcal{R}}\left(  1-\left\vert z_{j}\right\vert
^{2}\right)  ^{2}\left\vert \xi_{j}\right\vert \mu\left(  j\right)  \left\vert
\eta_{k}\right\vert \mu\left(  k\right)  .
\end{align*}
By Schur's test it suffices to show
\begin{align}
\mu\left(  Ak\right)  +\sum_{j=1}^{J}\left(  1-\left\vert z_{j}\right\vert
^{2}\right)  ^{\sigma_{1}^{\prime}}\mu\left(  j\right)   &  <C\varepsilon
<1,\label{Schurtest}\\
\sum_{k:z_{k}\in\mathcal{C}\left(  z_{j}\right)  }\mu\left(  k\right)
+\sum_{k\notin\mathcal{R}}\left(  1-\left\vert Az_{k}\right\vert ^{2}\right)
^{\sigma_{2}}\mu\left(  k\right)   &  <C\varepsilon<1,\nonumber\\
\sum_{k\in\mathcal{R}}\left(  1-\left\vert z_{j}\right\vert ^{2}\right)
^{3}\mu\left(  k\right)   &  <C\varepsilon<1.\nonumber
\end{align}
Now (\ref{logeps'}) yields
\[
\sum_{j=1}^{J}\left(  1-\left\vert z_{j}\right\vert ^{2}\right)  ^{\sigma
_{1}^{\prime}}\mu\left(  j\right)  \leq C\sum_{j}\left(  1-\left\vert
z_{j}\right\vert ^{2}\right)  ^{\sigma_{1}^{\prime\prime}}<C\varepsilon,
\]
and combined with the weak simple condition (\ref{weaksimple}), we have
\begin{align*}
\sum_{k=1}^{J}\left(  1-\left\vert Az_{k}\right\vert ^{2}\right)  ^{\sigma
_{2}}\mu\left(  k\right)   &  =\sum_{\ell}\left(  1-\left\vert z_{\ell
}\right\vert ^{2}\right)  ^{\sigma_{2}}\left(  \sum_{z_{k}\in\mathcal{C}%
\left(  z_{\ell}\right)  }\mu\left(  k\right)  \right) \\
&  \leq C\sum_{\ell}\left(  1-\left\vert z_{\ell}\right\vert ^{2}\right)
^{\sigma_{2}}\mu\left(  \ell\right) \\
&  \leq C\sum_{\ell}\left(  1-\left\vert z_{\ell}\right\vert ^{2}\right)
^{\sigma_{2}^{\prime}}<C\varepsilon.
\end{align*}
Finally we write the annulus $B\left(  0,1\right)  \setminus B\left(
0,R\right)  $ as a pairwise disjoint union $\cup_{i=1}^{N}B_{i}$ of Carleson
boxes of \textquotedblleft size\textquotedblright\ $R$ where $N\approx\left(
1-R^{2}\right)  ^{-1}$. Then
\[
\sum_{z_{k}\in B_{i}:k\in\mathcal{R}}\mu\left(  k\right)  \leq C\left(
1+\log\frac{1}{1-R^{2}}\right)  ^{-1}\leq C
\]
by the weak simple condition (\ref{weaksimple}), and thus the left side of the
final estimate in (\ref{Schurtest}) satisfies
\begin{align*}
\sum_{k\in\mathcal{R}}\left(  1-\left\vert z_{j}\right\vert ^{2}\right)
^{2}\mu\left(  k\right)   &  \leq\left(  1-R^{2}\right)  ^{2}\sum_{i=1}%
^{N}\sum_{z_{k}\in B_{i}:k\in\mathcal{R}}\mu\left(  k\right) \\
&  \leq C\left(  1-R^{2}\right)  ^{2}N\\
&  \leq C\left(  1-R^{2}\right)  <C\varepsilon,
\end{align*}
by (\ref{logeps'}) as required.

Thus $T^{-1}$ exists uniformly in $J$. Now we take $\xi\in\ell^{2}\left(
d\mu\right)  $ and set $\eta=\bigtriangleup\xi$. Here we use the convention
that $\xi\left(  A\alpha\right)  =0$ if $\alpha$ is a root of a tree in the
forest $Z$. By the weak simple condition we have the estimate
\begin{align}
\left\Vert \eta\right\Vert _{\ell^{2}\left(  d\mu\right)  }^{2}  &  =\sum
_{j}\left\vert \eta_{j}\right\vert ^{2}\mu\left(  j\right)  =\sum
_{j}\left\vert \xi_{j}-\xi_{Aj}\right\vert ^{2}\mu\left(  j\right)
\label{etazeta}\\
&  \leq C\sum_{j}\left\vert \xi_{j}\right\vert ^{2}\mu\left(  j\right)
+C\sum_{\ell}\left\vert \xi_{\ell}\right\vert ^{2}\left(  \sum_{z_{j}%
\in\mathcal{C}\left(  z_{\ell}\right)  }\mu\left(  j\right)  \right)
\nonumber\\
&  \leq C\sum_{j}\left\vert \xi_{j}\right\vert ^{2}\mu\left(  j\right)
+C\sum_{\ell}\left\vert \xi_{\ell}\right\vert ^{2}\mu\left(  \ell\right)
\nonumber\\
&  \leq C\left\Vert \xi\right\Vert _{\ell^{2}\left(  d\mu\right)  }%
^{2}.\nonumber
\end{align}
Then let $h=M\left(  T^{-1}\eta\right)  $ so that
\[
\bigtriangleup h\mid_{Z}=\bigtriangleup\left(  MT^{-1}\eta\right)  \mid
_{Z}=TT^{-1}\eta=\eta=\bigtriangleup\xi.
\]
Thus the holomorphic function $h$ satisfies
\[
h\mid_{Z}=\xi.
\]
Finally, from (\ref{maincon'}) and then (\ref{etazeta}) we have the Besov
space estimate (\cite{ArRoSa2}),
\begin{align*}
\left\Vert h\right\Vert _{B_{2}}^{2}  &  \leq C\sum_{j=1}^{J}\left\vert
\left(  T^{-1}\eta\right)  _{j}\right\vert ^{2}\int_{\mathbb{D}}\left\vert
\left(  1-\left\vert \zeta\right\vert ^{2}\right)  g_{w}\left(  \zeta\right)
\right\vert ^{2}d\lambda_{1}\left(  \zeta\right)  d\zeta\\
&  \leq C\sum_{j=1}^{J}\left\vert \left(  T^{-1}\eta\right)  _{j}\right\vert
^{2}\left(  \log\frac{1}{1-\left\vert w\right\vert ^{2}}\right)  ^{-1}\\
&  \leq C\left\Vert T^{-1}\eta\right\Vert _{\ell^{2}\left(  d\mu\right)  }%
^{2}\leq C\left\Vert \eta\right\Vert _{\ell^{2}\left(  d\mu\right)  }^{2}\leq
C\left\Vert \xi\right\Vert _{\ell^{2}\left(  d\mu\right)  }^{2}.
\end{align*}
Since all of this is uniform in $J$ we may let $J\rightarrow\infty$ and use a
standard normal families argument to complete the proof of the sufficiency of
the hypotheses in Theorem \ref{onto}.

\section{Relations Between Conditions\label{examples disk}}

\subsection{An Example Covered by Theorem B but not Theorem A}

Let $a,b>1$ satisfy $\left[  a^{k+1}b\right]  \geq\left[  a^{k}b\right]  +1$
for all $k\geq0$ (in particular this will hold if $\left(  a-1\right)  b\geq2
$), and define a Cantor-like sequence $_{Z}$%
\[
Z=Z_{a,b}=\cup_{k=0}^{\infty}\left\{  z_{j}^{k}\right\}  _{j=1}^{2^{k}}%
\subset\mathcal{T}%
\]
as follows. Pick a point $z_{1}^{0}$ of $\mathcal{T}$ satisfying $d\left(
z_{1}^{0}\right)  =\left[  b\right]  $. Then choose $2^{1}$ points $\left\{
z_{1}^{1},z_{2}^{1}\right\}  \subset\mathcal{T}$ that are successors to
distinct children of $z_{1}^{0}$ and having $d\left(  z_{j}^{1}\right)
=\left[  ab\right]  $, $1\leq j\leq2^{1}$, and, recalling that $\beta$ is the
hyperbolic distance on the disk, $\beta\left(  z_{i}^{1},z_{j}^{1}\right)
\gtrapprox\left[  ab\right]  $ for $i\neq j$. Then choose $2^{2}$ points
$\left\{  z_{1}^{2},z_{2}^{2},z_{3}^{2},z_{4}^{2}\right\}  \subset\mathcal{T}$
that are successors to distinct children of the points in $\left\{  z_{1}%
^{1},z_{2}^{1}\right\}  $ and having $d\left(  z_{j}^{2}\right)  =\left[
a^{2}b\right]  $, $1\leq j\leq2^{2}$, and $\beta\left(  z_{i}^{2},z_{j}%
^{2}\right)  \gtrapprox\left[  a^{2}b\right]  $ for $i\neq j$. Having
constructed $2^{k}$ points $\left\{  z_{j}^{k}\right\}  _{j=1}^{2^{k}}%
\subset\mathcal{T}$ in this way, we then choose $2^{k+1}$ points $\left\{
z_{j}^{k+1}\right\}  _{j=1}^{2^{k+1}}\subset\mathcal{T}$ that are successors
to distinct children of the points in $\left\{  z_{j}^{k}\right\}
_{j=1}^{2^{k}}$ and having $d\left(  z_{j}^{k+1}\right)  =\left[
a^{k+1}b\right]  $, $1\leq j\leq2^{k}$, and $\beta\left(  z_{i}^{k+1}%
,z_{j}^{k+1}\right)  \gtrapprox\left[  a^{k+1}b\right]  $ for $i\neq j$. Note
that the condition $\left[  a^{k+1}b\right]  \geq\left[  a^{k}b\right]  +1$
allows for the existence of such points. Then $Z=\cup_{k=0}^{\infty}\left\{
z_{j}^{k}\right\}  _{j=1}^{2^{k}}$ satisfies the separation condition
(\ref{sep}) with constant roughly $a-1$ and condition (\ref{special}) with
$\beta$ close to $1$, and the associated measure $\mu_{Z}$ satisfies the weak
simple condition (\ref{weaksimple}) with constant $2 $. Thus Theorem
\ref{onto} applies to show that $Z$ is onto interpolating for $B_{2}$. Yet the
total mass of the measure $\mu_{Z}$ satisfies
\[
\left\Vert \mu_{Z}\right\Vert =\sum_{k=0}^{\infty}2^{k}\left[  a^{k}b\right]
^{-1}\approx\frac{1}{b}\sum_{k=0}^{\infty}\left(  \frac{2}{a}\right)
^{k}=\infty
\]
if $a\leq2$.

\subsection{An Example Covered by Theorem A but not Theorem 3}

We now use a similar construction to give a separated sequence $W$ in the disk
with \emph{finite} measure $\mu=\mu_{W}$ satisfying the weak simple condition
but not the simple condition. This yields an example of a sequence which fails
the simple condition, but to which Theorem A applies.

We continue the notation of the previous example. We choose $a=2$ for
convenience, let $b,N$ be large integers, and replace the sequence $Z_{2,b}$
above with the truncated sequence $Z_{2,b,N}=\cup_{k=0}^{N}\left\{  z_{j}%
^{k}\right\}  _{j=1}^{2^{k}}$. Then $Z_{2,b,N}$ satisfies the separation
condition (\ref{sep}) with constant roughly $1$, the associated measure
$\mu_{Z_{2,b,N}}$ satisfies the weak simple condition (\ref{weaksimple}) with
constant $2$, and the total mass of $\mu_{Z_{2,b,N}}$ is about $\frac{N}{b}$:
\[
\left\Vert \mu_{Z_{2,b,N}}\right\Vert =\sum_{k=0}^{N}2^{k}\left[
2^{k}b\right]  ^{-1}\approx\frac{1}{b}\sum_{k=0}^{N}\left(  \frac{2}%
{2}\right)  ^{k}\approx\frac{N}{b}.
\]
On the other hand the constant $C\left(  \mu_{Z_{2,b,N}}\right)  $ in the
simple condition (\ref{simp}) for $\mu_{Z_{2,b,N}}$ satisfies%
\begin{equation}
C\left(  \mu_{Z_{2,b,N}}\right)  \gtrsim N, \label{ss}%
\end{equation}
since%
\[
\frac{N}{b}\approx\left\Vert \mu_{Z_{2,b}^{\ast}}\right\Vert =\sum_{\alpha\geq
z_{0}}\mu_{Z_{2,b}^{\ast}}\left(  \alpha\right)  \leq C\frac{1}{d\left(
z_{0}\right)  }=\frac{C}{b}.
\]

It is now an easy exercise to choose sequences of parameters $\left\{
b\left(  n\right)  \right\}  _{n=1}^{\infty}$ and $\left\{  N\left(  n\right)
\right\}  _{n=1}^{\infty}$, and initial points $\left\{  z_{1}^{0}\left(
n\right)  \right\}  _{n=1}^{\infty}$ so that the corresponding sequences
$Z_{2,b\left(  n\right)  ,N\left(  n\right)  }=\cup_{k=0}^{N\left(  n\right)
}\left\{  z_{j}^{k}\left(  n\right)  \right\}  _{j=1}^{2^{k}}$ satisfy%
\begin{equation}
\left\Vert \mu_{Z_{2,b\left(  n\right)  ,N\left(  n\right)  }}\right\Vert
\approx\frac{N\left(  n\right)  }{b\left(  n\right)  }\leq2^{-n} \label{fi}%
\end{equation}
and%
\begin{equation}
\lim_{n\rightarrow\infty}N\left(  n\right)  =\infty, \label{Ngoes}%
\end{equation}
along with the nested property%
\begin{equation}
z_{1}^{0}\left(  n+1\right)  \geq z_{1}^{b\left(  n\right)  }\left(  n\right)
,\ \ \ \ \ n\geq1. \label{nes}%
\end{equation}
Then the union $W=\cup_{n=1}^{\infty}Z_{2,b\left(  n\right)  ,N\left(
n\right)  }$ satisfies the separation condition and the associated measure
$\mu_{W}$ is finite by (\ref{fi}), satisfies the weak simple condition by
(\ref{nes}), yet fails the simple condition by (\ref{ss}) and (\ref{Ngoes}).

\section{Tree Interpolation and Theorem C\label{tree}}

\subsection{Reduction to a Basic Construction}

In this section we prove Theorem C.

Previously our tree $\mathcal{T}$ was constructed to contain our given
sequence $Z.$ In this section we regard $\mathcal{T}$ as a given, fixed,
Bergman tree and we will be interested in subsets of $\mathcal{T}.$ Our tree
has a bounded branching number but to keep the notation simple we suppose it
is a dyadic tree. We want to characterize the subsequences $Z=\left\{
\alpha_{j}\right\}  _{j=1}^{\infty}$ of $\mathcal{T}$ which are onto
interpolating sequence for $B_{2}\left(  \mathcal{T}\right)  $. We had defined
that class of sequences using the weighted restricting operator $R_{Z};$
however, taking note of the estimate $d\left(  \alpha_{j}\right)
\sim\left\Vert k_{z_{j}}\right\Vert _{B_{2}}^{2}$, we can alternatively
characterize the sequences by
\begin{align}
&  \text{for every sequence }\left\{  \xi_{j}\right\}  _{j=1}^{\infty}\text{
with }\left\Vert \left\{  d\left(  \alpha_{j}\right)  ^{-1/2}\xi_{j}\right\}
_{j=1}^{\infty}\right\Vert _{\ell^{2}}=1\text{,}\label{ontotree}\\
&  \text{there is }f\in B_{2}\left(  \mathcal{T}\right)  \text{ with
}\left\Vert f\right\Vert _{B_{2}\left(  \mathcal{T}\right)  }\leq C\text{ and
}f\left(  \alpha_{j}\right)  =\xi_{j}\text{, }j\geq1.\nonumber
\end{align}

Suppose $Z\subset\mathcal{T}$ \ is given and (\ref{treecap}) holds. One
immediate consequence is the tree separation condition,
\begin{equation}
\forall z,w\in Z,\text{ \ }d\left(  z,w\right)  \geq cd\left(  o,z\right)
\label{treesep}%
\end{equation}
This holds because the left hand side of (\ref{treecap}) a majorant for
\[
\inf\left\{  \sum\nolimits_{\zeta\in\mathcal{T}}\left\vert \bigtriangleup
f\left(  \zeta\right)  \right\vert ^{2}:f\left(  z\right)  =1,f\left(
w\right)  =0\right\}  =\frac{1}{d\left(  z,w\right)  -1}.
\]
This separation is certainly similar to (\ref{sep'}), however, because of
"edge effects" this condition is weaker.

We will prove lemma about the existence of certain almost extremal functions
in $B_{2}\left(  \mathcal{T}\right)  .$ They only take the values $0,1$ on
$Z,$ their discrete derivatives have disjoint support, and they have
controlled norms. Using the lemma the proof of Theorem C is immediate.

Let $E,F$ be disjoint subsets of $\mathcal{T}$. The capacity $Cap(E,F)$ of the
condenser $(E,F)$ is defined as
\[
Cap(E,F)=\inf\{\Vert h\Vert_{\ell^{2}}^{2}:\ Ih|_{E}=1,\ Ih|_{F}=0\}.
\]
For $Z$ a sequence in $\mathcal{T}$ set $\gamma(z,Z)=Cap(z,Z\setminus\{z\}).$
With this notation our tree capacity condition (\ref{treecap}) can be written
as
\begin{equation}
\gamma(z,Z)\leq Cd(z)^{-1}. \label{eqcap}%
\end{equation}
We say that $S\subset\mathcal{T}$ is a stopping region if every pair of
distinct points in $S$ are incomparable in $\mathcal{T}$.

\begin{lemma}
\label{shrunk}Given a subset $Z$ of $\mathcal{T}$ that satisfies the tree
separation condition (\ref{treesep}), there are functions $H_{w}=Ih_{w}$ on
$\mathcal{T}$ and a constant $C$ depending only on $C$ in (\ref{treesep}) satisfying

\begin{enumerate}
\item $H_{w}\left(  z\right)  =\delta_{w,z}^{{}}$ for $w,z\in Z$

\item $supp\left(  h_{w}\right)  \cap supp\left(  h_{z}\right)  =\phi$ for
$z,w\in Z$, $z\neq w$,

\item $\left\Vert H_{w}\right\Vert _{B_{2}\left(  \mathcal{T}\right)  }%
\leq2C\gamma\left(  w,Z\right)  $ for $w\in Z$,

\item If $S$ is a stopping region in $\mathcal{T}$, then $\sum_{\alpha\in S}$
$\left\vert h_{w}\left(  \alpha\right)  \right\vert \leq2\gamma\left(
w,Z\right)  $ for $w\in Z$.
\end{enumerate}
\end{lemma}

\textbf{Proof of Theorem C given the Lemma}: We've already mentioned that
(\ref{ontotree}) implies (\ref{treecap}). Conversely, if (\ref{treecap})
holds, then so does (\ref{treesep}) and hence also properties 1, 2 and 3 of
Lemma \ref{shrunk}. Let $Z=\left\{  z_{j}\right\}  _{j=1}^{\infty}$. If
$\xi=\left\{  \xi_{j}\right\}  _{j=1}^{\infty}\in\ell^{2}\left(  \mu\right)
$, then $f=\sum_{j=1}^{\infty}\xi_{j}H_{z_{j}}\in B_{2}\left(  \mathcal{T}%
\right)  $ with $B_{2}\left(  \mathcal{T}\right)  $ norm at most $C\left\Vert
\xi\right\Vert _{\ell^{2}\left(  \mu\right)  }$ by properties 2 and 3; and
$f\left(  z_{j}\right)  =\xi_{j}$ for all $j\geq1$ by property 1.

\subsection{Extremal Functions}

We now prove a string of results on capacity that will culminate in the proof
of Lemma \ref{shrunk}. More precisely, properties 1, 2 and 3 of Lemma
\ref{shrunk} will follow from Proposition \ref{propontint} below, and property
4 will follow from Proposition \ref{stoppingtime}. We use the following
notation. If $x$ is an element of the tree $T$, $x^{-1}$ denotes its immediate
predecessor in $T$. If $z$ is an element of the sequence $Z\subset T$, \ $Pz$
denotes its predecessor in $Z$: $Pz\in Z$ is the maximum element of
$Z\cap\left[  o,z\right)  $ Let $\Omega\subseteq T$. A point $x\in T$ is in
the interior of $\Omega$ if $x,x^{-1},x_{+},x_{-}\in\Omega$. A function $H$ is
\textit{harmonic} in $\Omega$ if
\begin{equation}
H(x)=\frac{1}{3}[H(x^{-1})+H(x_{+})+H(x_{-})] \label{mvp}%
\end{equation}
for every point $x$ which is interior in $\Omega$. If $H=Ih$ is harmonic in
$\Omega$, then we have that
\begin{equation}
h(x)=h(x_{+})+h(x_{-}) \label{addition}%
\end{equation}
whenever $x$ is in the interior of $\Omega$.

\subsubsection{Basic Properties}

\begin{proposition}
Let $T$ be a dyadic tree.

\begin{enumerate}
\item If $E$ and $F$ are finite, there is an extremal function $H=Ih$ such
that $Cap(E,F)=\Vert h\Vert_{\ell^{2}}^{2}$.

\item The function $H$ is harmonic on $T\setminus(E\cup F)$.

\item Let $E=\{z\}$, $F=Z-\{z\}$. Then the support of $h$ consists of (at
most) three connected components.. The upper support consists of the segment
$(Pz,z]$ and of all segments $[\zeta(w),w)$, where $w\in Z$ and $(\zeta(w),w]
$ has some intersection with the component of $T\setminus Z$ lying above $z$
(i.e., the component containing $z^{-1}$). The lower support consists of all
segments $[\zeta(w),w)$, where $w\in Z$ and $(\zeta(w),w]$ has some
intersection with one of the (at most) two components of $T\setminus Z$ lying
below $z$ (i.e., the components containing, respectively, $z_{+}$ and $z_{-}$).

\item The function $h$ is positive on $(Pz,z]$, negative on the segments
$[\zeta(w),w)$ and vanishes everywhere else.
\end{enumerate}
\end{proposition}

In classical potential theory capacities can be recovered from the derivative
of the Green potential. An analogous result holds here for the capacity of
$(z,Z\setminus\{z\})$ with $h$ playing the role of that derivative.

\begin{proposition}
Let $\gamma_{\pm}$, $\gamma_{P}$ be the $\ell^{2}$-sum of $h$ over the lower
and upper components of its support, respectively. Then,
\[
\gamma_{\pm}=-h(z_{\pm}),\ \gamma_{P}=h(z^{-1}).
\]

\end{proposition}

\textbf{Proof of both propositions}: We consider first the case of $\gamma
_{+}$. Let $T_{z}^{+}$ be the component of $T\setminus Z$ containing $x_{+}$
and let $\overline{T}_{z}^{+}$ be its forward closure, which is obtained by
adding to $T_{z}^{+}$ all points $w\in Z$, $w\neq z$, such that $d(w,T_{z}%
^{+})=1$. For each $x\in T_{z}^{+}$, $x_{\pm}\in\overline{T}_{z}^{+}$. We
proceed by induction on the cardinality of $Z\cap\overline{T}_{z}^{+}$.

If $\sharp(Z\cap\overline{T}_{z}^{+})=1,$ $Z\cap\overline{T}_{z}^{+}=\{w\} $,
then $supp(h)=(z,w]$ and $h=-d(z,w)^{-1}$ over its support, so that
$\gamma_{+}=d(z,w)\cdot d(z,w)^{-2}=d(z,w)^{-1}=-h(z_{+})$.

Suppose we know that $N\geq1$ and that the property holds when $\sharp
(Z\cap\overline{T}_{z}^{+})\leq N$ and suppose now that $\sharp(Z\cap
\overline{T}_{z}^{+})=N+1$. Consider $W_{z}=\{w\wedge w^{\prime}%
:\ w,w^{\prime}\in Z\cap\overline{T}_{z}^{+}\}$, the subtree of $\overline
{T}_{z}^{+}$ generated by $Z\cap\overline{T}_{z}^{+}$. Then, $W_{z}$ has a
minimal element $\zeta\in T_{z}^{+}$ ($\zeta\notin Z$ because $N+1\geq2$). Let
$U_{+}=S(\zeta_{+})\cap Z\cap\overline{T}_{z}^{+}$ and $U_{-}=S(\zeta_{-})\cap
Z\cap\overline{T}_{z}^{+}$. Then, the function $H$ goes from $1$ to
$H(\zeta)>0$ on $[z,\zeta] $ and from $H(\zeta)$ to $0$ as $x$ moves from
$\zeta$ to $U_{\pm}=S(\zeta_{\pm})\cap Z\cap\overline{T}_{z}^{+}$ in
$\{\zeta\}\cup X_{\pm}=\{\zeta\}\cup S(\zeta_{\pm})\cap\overline{T}_{z}^{+}$.
Note that
\[
H(\zeta)=1-d(z,\zeta)h(z_{+})=1-d(z,\zeta)h(\zeta).
\]
The function $h|_{X_{+}}$ has minimal $\ell^{2}$ norm with the property that
$\sum_{x=\zeta_{+}}^{w}=-H(\zeta)$ whenever $w\in U_{+}$, otherwise we could
modify it to obtain a global $H$ with better $\ell^{2}$ norm, contradicting
the hypothesis that $H$ was optimal. Also, observe that $\sharp(U_{\pm})\leq
N$.

Let $H_{+}=Ih_{+}$ be the function such that $H_{+}(\zeta)=1$, $H_{+}|_{U_{+}%
}=0$ and which has minimal $\ell^{2}$ norm with these properties. Similarly
define $H_{-}=Ih_{-}$ and $H^{\prime}=Ih^{\prime}$, the latter with the
conditions $H^{\prime}(z)=1$ and $H^{\prime}(\zeta)=0$. By the induction
hypothesis,
\[
\Vert h_{+}\Vert_{\ell^{2}}^{2}=h_{+}(\zeta_{+}),\ \Vert h_{-}\Vert_{\ell^{2}%
}^{2}=h_{-}(\zeta_{-}),\ \Vert h^{\prime}\Vert_{\ell^{2}}^{2}=h^{\prime}%
(z_{+})=h^{\prime}(\zeta).
\]
By uniqueness of the extremal function $H$ and by homogeneity of the
minimization problem, we have that, with constants $a_{+}=H(\zeta
)=a_{-},\ a^{\prime}=1-H(\zeta)$,
\[
h=a_{\pm}\cdot h_{\pm},\ h=a^{\prime}\cdot h^{\prime}.
\]
The norm of $h$ is then
\begin{align*}
\Vert h\Vert_{\ell^{2}}^{2}  &  =a_{+}^{2}\Vert h_{+}\Vert_{\ell^{2}}%
^{2}+a_{-}^{2}\Vert h_{-}\Vert_{\ell^{2}}^{2}+{a^{\prime}}^{2}\Vert h^{\prime
}\Vert_{\ell^{2}}^{2}\\
&  =H(\zeta)^{2}\cdot\lbrack h_{+}(\zeta_{+})+(h_{-}(\zeta_{-})]+(1-H(\zeta
))^{2}\cdot h^{\prime}(\zeta)\\
&  =H(\zeta)\cdot\lbrack h(\zeta_{+})+h(\zeta_{-})]+(1-H(\zeta))\cdot
h(\zeta)\\
&  =h(\zeta)=h(z_{+}),
\end{align*}
as wished. In the fourth equality, we used the harmonicity of $H=Ih$.

Exactly the same argument works for $\gamma_{-}$ and a variation thereof gives
the desired formula for $\gamma_{P}$.

\begin{remark}
The proof gives a useful formula for computing capacities. Given $z$,$\zeta
\in\mathcal{T},$ $z<\zeta$ and given $U_{\pm}\subset S(\zeta_{\pm})$ we have
\begin{equation}
Cap(z,U_{+}\cup U_{-})=\frac{Cap(\zeta,U_{+})+Cap(\zeta,U_{-})}{1+d(z,\zeta
)[Cap(\zeta,U_{+})+Cap(\zeta,U_{-})]}. \label{formula}%
\end{equation}
To see this note that,
\begin{align*}
Cap(z,U_{+}\cup U_{-})  &  =h(z_{+})=h(\zeta)=h(\zeta_{+})+h(\zeta_{-})\\
&  =H(\zeta)[h_{+}(\zeta_{+})+h_{-}(\zeta_{-})]\\
&  =(1-d(z,\zeta)h(\zeta_{+}))[Cap(\zeta,U_{+})+Cap(\zeta,U_{-})]\\
&  =(1-d(z,\zeta)\cdot Cap(z,U_{+}\cup U_{-}))[Cap(\zeta,U_{+})+Cap(\zeta
,U_{-})].
\end{align*}

\end{remark}

Later in this section we will use this formula to develop a simple, computable
and geometric algorithm for calculating capacities.

\subsubsection{Disjoint supports}

From now on, we consider a \emph{finite} sequence $Z$ satisfying
(\ref{treesep}). We want to show that the functions $h=h_{z}$ can be replaced
by near extremal functions $k_{z}$, with the extra property that
$supp(k_{z})\cap supp(k_{w})=\phi$ if $z\neq w$.

We will assume that $Z=\{z_{j}:\ j\geq0\}$ is ordered in such a way that
$d(z_{n})\leq d(z_{n+1})$. We will also assume that $z_{0}=o=0$ belongs to
$Z$. We define
\[
Z_{n}=\{z_{0},\dots,z_{n}\}
\]
and
\[
\mathcal{T}_{n}=\cup_{j=1}^{n}[o,z_{j}]
\]
We also need ${\mathcal{T}}_{\infty}=\cup_{j=1}^{\infty}[o,z_{j}]$, the
minimum subset of ${\mathcal{T}}$ containing $Z$ which is geodesically
connected. The \textit{landing point} of $z=z_{n+1}$ is, by definition,
\[
\xi(z)=\max\left(  [o,z]\cap\mathcal{T}_{n}\right)
\]
By construction, $Z\cap(\xi(z),z)=\phi$ and, if $z\neq w$, $[\xi(z),z]$ and
$[\zeta(w),w]$ are either disjoint, or they intersect in $\xi(z)$, or in
$\xi(w)$. If $i_{j}\geq i_{k}$, then $\xi(z_{i_{k}})\in\lbrack\xi(z_{i_{j}%
}),z_{i_{j}}]$ if and only if $\xi(z_{i_{j}})=\xi(z_{i_{k}})$.

\begin{lemma}
\label{lemmasep}(\cite{ArRoSa}; Lemma 27) Let $Z$ be a sequence satisfying
(\ref{treesep}). Then, for some positive constant $\eta$,
\begin{equation}
d(\xi(z),z)\geq\eta d(z) \label{eqlanding}%
\end{equation}
for all $z$ in $Z$.
\end{lemma}

As a consequence, by removing a finite number of points from $Z$ we can assume
that $d(\xi(z),z)\geq3$ for all $z\neq o$.

We record some further properties of the functions $H=H_{z}$.

\begin{proposition}
\label{lemmacap} The function $H$ is increasing on $[Pz,z]$ and decreasing on
all segments of the form $[\zeta(w),w]$. Moreover, $H$ is convex on all
intervals of the form $[\zeta(w),w]$ and on $[Pz,z]$.
\end{proposition}

For $z$ in $Z$, let
\[
{\mathcal{N}}(z)=\bigcup\limits_{\substack{w\geq P(z) \\(w,P(z))\cap
Z\subseteq\{z\}}}[w,P(z))
\]
The set ${\mathcal{N}}(z)$ is the \textquotedblleft downward
closure\textquotedblright\ of the connected component ${\mathcal{N}}_{0}(z)$
of ${\mathcal{T}}_{\infty}\setminus\left(  Z\setminus\left\{  z\right\}
\right)  $ containing $z$. We have
\[
{\mathcal{N}}_{0}(z)={\mathcal{N}}(z)\setminus\left(  Z\setminus\left\{
z\right\}  \right)  .
\]
All of the interesting action takes place inside ${\mathcal{N}}(z)$.

The main tool we need is the following. We write $DK=k$ if $K=Ik$.

\begin{proposition}
\label{propontint} If (\ref{eqcap}) holds, then to each $z$ in $Z$ we can
associate a function $K_{z}$ in $B_{2}$ such that $K_{z}\left(  w\right)
=\delta_{z}\left(  w\right)  $ for $w\in Z$, $K_{z}$ is harmonic on
$supp\left(  DK_{z}\right)  \setminus Z$, $\Vert K_{z}\Vert_{B_{2}}\leq
C\gamma\left(  z,Z\right)  $ and such that $\mbox{supp}(DK_{z})\cap
\mbox{supp}(DK_{w})=\phi$ if $z\neq w$. In fact the $K^{\prime}s$ are pairwise
disjoint as well: $\mbox{supp}(K_{z})\cap\mbox{supp}(K_{w})=\phi$ if $z\neq w$.
\end{proposition}

\begin{remark}
Because the $K^{\prime}s$ have disjoint support we have an immediate solution,
$\sum_{j=1}^{\infty}\xi_{j}K_{z_{j}},$ to the problem of interpolating
$\xi=\left\{  \xi_{j}\right\}  _{j=1}^{\infty}$ on the tree $T$.
\end{remark}

\textbf{Proof of the proposition:} Fix $z\in Z$. By Proposition \ref{lemmacap}
there is a function $H$ such that $\Vert H\Vert_{B_{2}}^{2}=\gamma(z,Z)$ and
$H(z)=1$, $H|_{Z\setminus\{z\}}\equiv0$, $H\geq0$, $\mbox{supp}(DH)\subseteq
{\mathcal{N}}(z)$ and $H$ is convex on intervals of the form $[w,\xi(w)]$. Let
$Q_{z}$ be the first point $x$ on $[P(z),z]$ such that $d(x,z)\leq\frac{1}%
{3}d(P(z),z)$. For each $w\in{\mathcal{N}}(z)\cap(Z\setminus\{z\})\cap
S(\xi(z))$, let $Q_{w}$ be the first point $y$ in $[\xi(w),w)$ such that
$d(y,\xi(w))\geq\frac{1}{3}d(\xi(w),w)$. As in the proof of Proposition
\ref{lemmacap}, ${\mathcal{N}}(z)\cap(Z\setminus\{z\})\cap S(\xi
(z))=\{z_{i_{j}}:\ j\geq0\}$, where $z=z_{i_{0}}$, $d(z_{i_{j+1}})\geq
d(z_{i_{j}})$ and $\xi(z_{i_{j+1}})\in{\mathcal{T}}_{i_{j}}$. The function
$K_{z}$ is constructed inductively, in such a way that
\[
\mbox{supp}(k_{z})\subseteq{\mathcal{N}}_{r}\left(  z\right)  \equiv\cup
_{j}[Q_{z_{i_{0}}},Q_{z_{i_{j}}}],\ \ \ \ \ k_{z}=DK_{z}.
\]
It is clear that if $z\neq w$ are points in $Z$, then ${\mathcal{N}}%
_{r}(z)\cap{\mathcal{N}}_{r}(w)=\phi$.

\textbf{Step $j=0$}. Construct a new function $K_{0}$ as follows. On
$[P(z),Q_{z})$, set $K_{0}=0$. For $x\in\lbrack Q_{z},z]$, set
\[
K_{0}(x)=\frac{H(x)-H(Q_{z})}{H(z)-H(Q_{z})}%
\]
If $y$ is such that $y\wedge z\in\lbrack P(z),z]$, set
\[
K_{0}(y)=\frac{K_{0}(y\wedge z)}{H(y\wedge z)}H(y)
\]
Set $K_{0}(y)=H(y)$ otherwise. The function $K_{0}$ has the following
properties (with $\eta$ from (\ref{eqlanding}) and with the parameter $\ell=0$),

\begin{enumerate}
\item $K_{\ell}$ is admissible for $z$, $K_{\ell}\geq0$, $K_{\ell}$ is
harmonic on ${\mathcal{N}}_{r}\left(  z\right)  $, $K_{\ell}$ is convex on
$[\xi(z_{i_{j}}),z_{i_{j}}]$ for $j\geq1$;

\item $K_{\ell}\leq H$, pointwise;

\item if $x\notin\lbrack Q_{z},z]$, then $|DK_{\ell}(x)|\leq|DH(x)|$;

\item if $x\in\lbrack Q_{z},z]$, then $|DK_{\ell}(x)|\leq C|DH(x)|$, where
$C=3/\eta$;

\item $\mbox{supp}(k_{\ell})\subset\lbrack Q_{z},z]\cup(\cup_{j\geq1}%
[z_{i_{j}},\xi(z_{i_{j}})])$.
\end{enumerate}

Properties 1, 2, 3 and 5 immediately follow from the construction. We show 4.
By convexity,
\[
\frac{1-H(Q_{z})}{d(z)-d(Q_{z})}=\frac{H(z)-H(Q_{z})}{d(z)-d(Q_{z})}\geq
\frac{H(z)-H(o)}{d(z)-d(o)}=\frac{1}{d(z)-1}%
\]
hence,
\[
H(z)-H(Q_{z})\geq\frac{d(z)-d(Q_{z})}{d(z)-1}\geq\frac{\eta}{3}%
\]
by Lemma \ref{lemmasep} and by the definition of $Q_{z}$. Thus, if
$x\in\lbrack Q_{z},z]$, then $|DK_{0}(x)|=|DH(x)|/(H(z)-H(Q_{z}))\leq C|DH(x)|
$.

\textbf{Induction.} Pick now $z_{i_{1}}$. If $\xi(z_{i_{1}})\in\lbrack
P(z),Q_{z})$, let $K_{1}=K_{0}$. If $\xi(z_{i_{1}})\in(z,Q_{z}]$, construct
$K_{1}$ as we did above, changing first the values of $K_{0}$ on
$[\xi(z_{i_{1}}),z_{i_{1}}]$ in such a way $K_{1}$ vanishes on $(Q_{z_{i_{1}}%
},z_{i_{1}}]$, then adjusting the values elsewhere to make $K_{1}$ admissible
for $z$. Observe that $\{x:K_{1}\neq K_{0}\}\cap\lbrack P(z),z]=\phi$. The
function $K_{1}$ has ($\ell=1$) properties 1, 2 (that is, $K_{\ell}\leq
K_{\ell-1}$), 3 (that is, if $x\notin\lbrack Q_{z_{i_{\ell}}},z_{i_{\ell}}]$,
then $|DK_{\ell}(x)|\leq|DK_{\ell-1}(x)|$), 4 (that is, if $x\in\lbrack
Q_{z_{i_{j}}},z_{i_{j}}]$, then $|DK_{\ell}(x)|\leq C|DK_{\ell-1}(x)|\leq
C|DH(x)|$, and this time $C=3$ suffices). Property 5 becomes
\[
\mbox{supp}(K_{\ell})\subset\lbrack Q_{z},z]\cup(\cup_{1\leq j\leq\ell
}[Q_{z_{i_{j}}},z_{i_{j}}])\cup(\cup_{j\geq\ell+1}[z_{i_{j}},\xi(z_{i_{j}})])
\]
Moreover, we have that%
\[
DK_{\ell}(x)=DK_{\ell-1}(x)\text{on the set}B_{\ell-1}=\cup_{1\leq j\leq
\ell-1}[\xi(z_{i_{j}}),z_{i_{j}}].
\]

Passing to the limit, we find a function $K_{z}$ with all the desired
properties, since the estimates for $|DK_{\ell}|$ add up nicely.

\subsubsection{Stopping estimates}

Let $S$ be a \textit{stopping region} in the tree $T$. By this, we mean that
there are no $x,y\in S$ such that $x>y;$ equivalently $S(x)\cap S(y)=\phi$ if
$x\neq y$.

\begin{proposition}
\label{stoppingtime}Let $\tilde{K}_{z}=I\tilde{k}_{z}$ be the function in
Proposition \ref{propontint} and let $S$ be a stopping region. Then,
\[
\sum_{x\in S}|\tilde{k}_{z}(x)|\leq2\gamma\left(  z,Z\right)  .
\]

\end{proposition}

\textbf{Proof}: Let $k=\tilde{k}_{z}$. We know that $k\geq0$ on $[Q_{z},z]$
and that $k\leq0$ elsewhere. Let $S$ be a stopping time in $T$. Then
$S\cap\lbrack Q_{z},z]$ consists of at most one element $x$ and
\[
0\leq k(x)\leq k(z)=\gamma_{P},
\]
since $k$ is convex on $[Q_{z},z]$.

Let $x_{0}=Q_{z},x_{1},\dots,x_{N}=z$ be an ordered enumeration of the points
in $[Q_{z},z^{-1}]$. Let $S_{j}=S\cap(S(x_{j})-S(x_{j+1}))$, $j=0,\dots,N-1$.
Without loss of generality, we can assume that $x_{j+1}={x_{j}}_{+}$, so that
$S(x_{j})-S(x_{j+1})=S({x_{j}}_{-})$. By harmonicity of $k$, (\ref{addition})
and easy induction, we have that
\[
\sum_{x\in S_{j}}|k(x)|=-\sum_{x\in S_{j}}k(x)=-k\left(  {x_{j}}_{-}\right)
=-\left[  k\left(  x_{j}\right)  -k\left(  {x_{j}}_{+}\right)  \right]  .
\]
Summing over $j$, we have
\[
\sum_{j}\sum_{x\in S_{j}}|k(x)|=k(z)-k(Q_{z})\leq k(z)=\gamma_{P}.
\]
Let $S_{\pm}=S\cap S(z_{\pm})$. Induction and harmonicity show that
\[
\sum_{j}\sum_{x\in S_{\pm}}|k(x)|=|k(z_{\pm})|=\gamma_{\pm}.
\]
All points in the support of $k$ fall in $S(Q_{z})\cup\dots\cup S(z)$, hence
\[
\sum_{x\in S}|k(x)|\leq2\gamma_{P}+\gamma_{+}+\gamma_{-}\leq2\gamma(z,Z).
\]

\section{More Relations Between Conditions\label{examples tree}}

\subsection{The Weak Simple Condition is not Necessary\label{failweak}}

Fix a point $z_{0}\in\mathcal{T}$ with $d\left(  z_{0}\right)  $ large. Let $N
$ and $b$ be positive integers. Fix a geodesic segment $\left[  z_{0}%
=w_{0},w_{1},w_{2},...,w_{N}\right]  $ in $\mathcal{T}$ and choose for each
$1\leq n\leq N$ a point $z_{n}$ satisfying%
\begin{align*}
d\left(  w_{n},z_{n}\right)   &  =b,\\
w_{n+1}  &  \notin\left[  w_{n},z_{n}\right]  ,\ \ \ \ \ 1\leq n<N.
\end{align*}
Thus $z_{n}$ and $w_{n+1}$ lie on different branches below $w_{n}$ and
$z_{n}\wedge w_{n+1}=w_{n}$. Set $Z_{N}=\left\{  z_{n}\right\}  _{n=1}^{N}$
and $\beta=\frac{1}{b}$, so that $Cap\left(  w_{n},\left\{  z_{n}\right\}
\right)  =\beta$ for each $n$. Then the formula in (\ref{formula}) above shows
that%
\begin{align*}
Cap\left(  z_{0},Z_{N}\right)   &  =\frac{1}{1+\frac{1}{Cap_{w_{1}}\left(
\left\{  z_{1}\right\}  \right)  +\text{ \ \ }\cdot\cdot\cdot\cdot\text{
\ \ +}\frac{1}{1+\frac{1}{Cap_{w_{N}}\left(  \left\{  z_{N}\right\}  \right)
}}}}\\
&  =\frac{1}{1+\frac{1}{\beta+\frac{1}{1+\frac{1}{\beta+\text{ \ \ \ }%
\cdot\cdot\cdot\text{ \ \ }\frac{1}{1+\frac{1}{\beta}}}}}}\equiv\gamma_{N}.
\end{align*}
The function $\varphi_{\beta}\left(  x\right)  =\frac{1}{1+\frac{1}{\beta+x}}$
is strictly increasing so if we take $\gamma_{0}=0$ and note that
$\gamma_{N+1}=\varphi_{\beta}\left(  \gamma_{N}\right)  $ we see that
$\gamma_{N}<\gamma_{N+1}$ for all $N\geq0.$ If $\gamma_{\infty}=\lim
_{N\rightarrow\infty}\gamma_{N}$ denotes the corresponding infinite continued
fraction,%
\[
\gamma_{\infty}=\frac{1}{1+\frac{1}{\beta+\frac{1}{1+\frac{1}{\beta+\frac
{1}{\ddots}}}}},
\]
then $Cap_{z_{0}}\left(  Z_{N}\right)  =\gamma_{N}<\gamma_{\infty}$. Since
$\gamma_{\infty}=\varphi_{\beta}\left(  \gamma_{\infty}\right)  =\frac
{1}{1+\frac{1}{\beta+\gamma_{\infty}}}$ we compute that
\[
\gamma_{\infty}=\frac{\sqrt{\beta^{2}+4\beta}-\beta}{2}<\sqrt{\beta}%
\]
since $0<\beta<1$. Altogether then,%
\[
Cap\left(  z_{0},Z_{N}\right)  <\sqrt{\beta},\ \ \ \ \ N\geq1.
\]

Now we fix $N=b=d\left(  z_{0}\right)  ^{2}$ so that
\[
Cap\left(  z_{0},Z_{N}\right)  <\sqrt{\frac{1}{d\left(  z_{0}\right)  ^{2}}%
}=\frac{1}{d\left(  z_{0}\right)  }.
\]
From this we obtain that%
\[
Z\equiv\left\{  z_{0}\right\}  \cup Z_{N}=\left\{  z_{n}\right\}  _{n=0}^{N}%
\]
satisfies the tree capacity condition (\ref{treecap}) with constant $C$ (for
$n\geq1$ the capacity $Cap\left(  z_{n},\left\{  z_{m}\right\}  _{m\neq
n}\right)  $ is easily seen to be bounded by $Cd\left(  z_{0}\right)  ^{-2}$
since the distance from $z_{n}$ to $w_{n}$ is $d\left(  z_{0}\right)  ^{2}$,
and the geodesics from $z_{n}$ to another point of $Z$ must pass through
$w_{n}$). The separation condition (\ref{treesep}) holds with constant close
to $\frac{1}{2}$ since $d\left(  z_{i},z_{j}\right)  \geq d\left(
z_{0}\right)  ^{2}$ and $d\left(  z_{n}\right)  \leq d\left(  z_{0}\right)
+2d\left(  z_{0}\right)  ^{2}$.

\begin{remark}
It is possible to choose the points $z_{j}$ above so that they are separated
with the same constant in the Bergman metric in the disk $\mathbb{D}$.
\end{remark}

On the other hand, the weak simple condition constant for $Z$ is quite large
since $\left\{  z_{n}\right\}  _{n=1}^{N}$ are the children of $z_{0}$ and%
\begin{align*}
\sum_{n=1}^{N}\mu\left(  z_{n}\right)   &  =\sum_{n=1}^{N}\frac{1}{d\left(
z_{n}\right)  }=\sum_{n=1}^{d\left(  z_{0}\right)  ^{2}}\frac{1}{d\left(
z_{0}\right)  +n+d\left(  z_{0}\right)  ^{2}}\\
&  \approx\log\frac{d\left(  z_{0}\right)  +2d\left(  z_{0}\right)  ^{2}%
}{d\left(  z_{0}\right)  +1+d\left(  z_{0}\right)  ^{2}}\approx\log2,
\end{align*}
which is much larger than $\mu\left(  z_{0}\right)  =\frac{1}{d\left(
z_{0}\right)  }$ if $d\left(  z_{0}\right)  $ is large.

Finally, even the constant in the \emph{enveloping weak simple condition}
(i.e. there is $Z^{\prime}\supset Z$ such that (\ref{weaksimple}) holds for
$Z^{\prime}$) must remain quite large if the separation condition constant is
to remain under control, i.e. \emph{not} go to zero. Indeed, suppose we can
add points $\left\{  v_{k}\right\}  _{k=1}^{K}$ to $Z$ so that the weak simple
condition for $Z^{\prime}=\left\{  v_{k}\right\}  _{k=1}^{K}\cup Z$ holds with
constant $C$. Without loss of generality we may assume that the points
$\left\{  v_{k}\right\}  _{k=1}^{K}$ lie along the geodesic segment $\left[
w_{1},w_{2},...,w_{d\left(  z_{0}\right)  ^{2}}\right]  $. If we consider the
weak simple condition at $w_{m}$ where $m=d\left(  z_{0}\right)  ^{2}-RC$ (and
where $R$ is a sufficiently large positive integer: $R=10$ surely works), then
we must have a point $v_{K}$ lying below $w_{m}$ since otherwise%
\begin{align*}
\sum_{n=m+1}^{N}\mu\left(  z_{n}\right)   &  =\sum_{n=d\left(  z_{0}\right)
^{2}-RC+1}^{d\left(  z_{0}\right)  ^{2}}\frac{1}{d\left(  z_{0}\right)
+n+d\left(  z_{0}\right)  ^{2}}\\
&  \approx\log\frac{d\left(  z_{0}\right)  +2d\left(  z_{0}\right)  ^{2}%
}{d\left(  z_{0}\right)  +2d\left(  z_{0}\right)  ^{2}-RC+1}\approx\frac
{RC}{d\left(  z_{0}\right)  ^{2}},
\end{align*}
which is \emph{not} bounded by%
\[
C\mu\left(  w_{m}\right)  =\frac{C}{d\left(  w_{m}\right)  }=\frac{C}{d\left(
z_{0}\right)  +d\left(  z_{0}\right)  ^{2}-RC}\approx\frac{C}{d\left(
z_{0}\right)  ^{2}}%
\]
provided $RC$ is much smaller than $d\left(  z_{0}\right)  ^{2}$ and $R$ is
sufficiently large. The same argument shows that there must be a point
$v_{K-1}$ lying in the segment $\left[  w_{p},w_{m}\right)  $ where
$p=m-RC=d\left(  z_{0}\right)  ^{2}-2RC$. But then $d\left(  v_{K-1}%
,v_{K}\right)  \leq RC$ while
\[
\min\left\{  d\left(  v_{K-1}\right)  ,d\left(  v_{K}\right)  \right\}
=d\left(  z_{0}\right)  +d\left(  z_{0}\right)  ^{2}-2RC,
\]
which shows that the separation constant is at most%
\[
\frac{d\left(  v_{K-1},v_{K}\right)  }{\min\left\{  d\left(  v_{K-1}\right)
,d\left(  v_{K}\right)  \right\}  }=\frac{RC}{d\left(  z_{0}\right)  +d\left(
z_{0}\right)  ^{2}-2RC}\approx\frac{RC}{d\left(  z_{0}\right)  ^{2}},
\]
a spectacularly small number if $d\left(  z_{0}\right)  $ is large.

\subsection{Separation Plus Weak Simple Implies Tree Capacity}

\begin{proposition}
\label{weaksimtree}If $Z$ satisfies both (\ref{treesep}) and (\ref{weaksimple}%
), then $Z$ satisfies the tree capacity condition (\ref{treecap}).
\end{proposition}

\textbf{Proof}: Fix $z_{0}\in Z$ and let $Z^{\prime}=\left\{  z_{j}\right\}
_{j=1}^{\infty}$ be those points in $Z$ whose geodesic to $z_{0}$ contains no
other points of $Z$. Consider for the moment the case where $Z^{\prime}\subset
S\left(  z_{0}\right)  $. Arrange the sequence $Z^{\prime}$ so that
$\left\vert z_{j+1}\right\vert \geq\left\vert z_{j}\right\vert $ for all
$j\geq1$. Define $\tau_{k}$ to be the smallest connected subset of
$\mathcal{T}$ containing $\left\{  z_{j}\right\}  _{j=0}^{k}$. We then define
the landing point $\xi_{k}$ of $z_{k}$ on $\tau_{k-1}$ as the maximal point on
the geodesic segment $\left[  o,z_{k}\right]  \cap\tau_{k-1}$. We now claim
that%
\begin{equation}
d\left(  z_{k},\xi_{k}\right)  \geq\frac{C}{2}d\left(  z_{k}\right)
,\ \ \ \ \ k\geq1, \label{landinglengths}%
\end{equation}
where $C$ is the constant in (\ref{treesep}). Indeed, there is $z_{j}\in
Z\cap\tau_{k-1}$ with $j<k$ such that $z_{j}\geq\xi_{k}$, and it follows that%
\begin{align*}
Cd\left(  z_{k}\right)   &  \leq d\left(  z_{j},z_{k}\right)  \leq d\left(
z_{j},\xi_{k}\right)  +d\left(  \xi_{k},z_{k}\right) \\
&  =d\left(  z_{j}\right)  -d\left(  \xi_{k}\right)  +d\left(  z_{k}\right)
-d\left(  \xi_{k}\right) \\
&  \leq2d\left(  z_{k}\right)  -2d\left(  \xi_{k}\right)  =2d\left(  z_{k}%
,\xi_{k}\right)  .
\end{align*}

Now define $h_{k}=\frac{1}{d\left(  z_{k},\xi_{k}\right)  }\chi_{\left(
\xi_{k},z_{k}\right]  }$ for $k\geq1$ and $h_{0}=\frac{1}{d\left(
z_{k}\right)  }\chi_{\left(  o,z_{k}\right]  }$. If we set $h=h_{0}-\sum
_{k=1}^{\infty}h_{k} $ and $f=Ih$, then we have $f\left(  o\right)  =0$,
$f\left(  z_{0}\right)  =1$, and $f\left(  z_{j}\right)  =0$ for $j\geq1$.
Since $Z\subset S\left(  z_{0}\right)  $ by our momentary assumption, we
actually have $f\left(  z\right)  =0$ for all $z\in Z\setminus\left\{
z_{0}\right\}  $. We also have the norm estimate%
\begin{align*}
\left\Vert f\right\Vert _{B_{2}\left(  \mathcal{T}\right)  }^{2}  &
=\sum_{k=0}^{\infty}\left\Vert h_{k}\right\Vert _{\ell^{2}}^{2}=\sum
_{\alpha\in\left(  o,z_{0}\right]  }\left(  \frac{1}{d\left(  z_{0}\right)
}\right)  ^{2}+\sum_{k=1}^{\infty}\sum_{\alpha\in\left(  \xi_{k},z_{k}\right]
}\left(  \frac{1}{d\left(  z_{k},\xi_{k}\right)  }\right)  ^{2}\\
&  =\frac{1}{d\left(  z_{0}\right)  }+\sum_{k=1}^{\infty}\frac{1}{d\left(
z_{k},\xi_{k}\right)  }\\
&  \leq\frac{1}{d\left(  z_{0}\right)  }+\frac{2}{C}\sum_{k=1}^{\infty}%
\frac{1}{d\left(  z_{k}\right)  }\leq C\frac{1}{d\left(  z_{0}\right)  },
\end{align*}
by (\ref{landinglengths}) and the weak simple condition (\ref{weaksimple}).
Thus we obtain $\gamma\left(  z_{0},Z^{\prime}\right)  \leq C\mu\left(
z_{0}\right)  $, and so (\ref{treecap}) holds for $\alpha=z_{0}$.

The general case, where not all points in $Z^{\prime}$ lie in $S\left(
z_{0}\right)  $, is handled as follows. Let $w$ be the point on the geodesic
$\left(  o,z_{0}\right]  $ satisfying $d\left(  w,z_{0}\right)  =\left[
Cd\left(  z_{0}\right)  \right]  $. Note that by (\ref{treesep}) there are no
points of $Z$ in the segment $\left(  w,z_{0}\right)  $. Now redefine
$h_{0}=\frac{1}{d\left(  w,z_{0}\right)  }\chi_{\left(  w,z_{0}\right]  }$.
For each point $\xi$ in the segment $\left(  w,z_{0}\right)  $, let $\xi+$ be
the child of $\xi$ on the geodesic $\left(  w,z_{0}\right]  $, and consider
the subset
\[
Z_{\xi}^{\prime}=Z^{\prime}\cap\left(  S\left(  \xi\right)  \setminus S\left(
\xi+\right)  \right)  .
\]
Just as we did for $Z^{\prime}\cap S\left(  z_{0}\right)  $ above, we
construct a function $h_{\xi}=\sum_{z_{k}\in Z_{\xi}^{\prime}}h_{k}$ such that
$Ih_{\xi}\left(  \xi\right)  =0$ and $Ih_{\xi}\left(  z_{k}\right)  =1$ for
$z_{k}\in Z_{\xi}^{\prime}$. Now define
\[
h=h_{0}-\sum_{z_{k}\in S\left(  z_{0}\right)  }h_{k}-\sum_{\xi\in\left(
w,z_{0}\right)  }Ih_{0}\left(  \xi\right)  \sum_{z_{k}\in Z_{\xi}^{\prime}%
}h_{k}.
\]
Then $f=Ih$ satisfies $f\left(  o\right)  =0$, $f\left(  z_{0}\right)  =1$,
and $f\left(  z\right)  =0$ for all $z\in Z\setminus\left\{  z_{0}\right\}  $.
Moreover we have the norm estimate%
\begin{align*}
\left\Vert f\right\Vert _{B_{2}\left(  \mathcal{T}\right)  }^{2}  &
=\left\Vert h_{0}\right\Vert _{\ell^{2}}^{2}+\sum_{z_{k}\in S\left(
z_{0}\right)  }\left\Vert h_{k}\right\Vert _{\ell^{2}}^{2}+\sum_{\xi\in\left(
w,z_{0}\right)  }\left\vert Ih_{0}\left(  \xi\right)  \right\vert ^{2}%
\sum_{z_{k}\in Z_{\xi}^{\prime}}\left\Vert h_{k}\right\Vert _{\ell^{2}}^{2}\\
&  \leq\frac{1}{d\left(  z_{0}\right)  }+\frac{2}{C}\sum_{z_{k}\in S\left(
z_{0}\right)  }\frac{1}{d\left(  z_{k}\right)  }+\frac{2}{C}\sum_{\xi
\in\left(  w,z_{0}\right)  }\left\vert Ih_{0}\left(  \xi\right)  \right\vert
^{2}\sum_{z_{k}\in Z_{\xi}^{\prime}}\frac{1}{d\left(  z_{k}\right)  }\\
&  \leq\frac{1}{d\left(  z_{0}\right)  }+\frac{2}{C}\sum_{z_{k}\in S\left(
z_{0}\right)  }\frac{1}{d\left(  z_{k}\right)  }+\frac{2}{C}\sum_{z_{k}\in
S\left(  w\right)  }\frac{1}{d\left(  z_{k}\right)  }\\
&  \leq\frac{1}{d\left(  z_{0}\right)  }+C\frac{1}{d\left(  z_{0}\right)
}+C\frac{1}{d\left(  w\right)  }\\
&  \leq C\frac{1}{d\left(  z_{0}\right)  }.
\end{align*}
Thus again $\gamma\left(  z_{0},Z^{\prime}\right)  \leq C\mu\left(
z_{0}\right)  $ and (\ref{treecap}) holds for $\alpha=z_{0}$.

\subsection{Separation Plus Finite Measure Doesn't Imply Interpolation}

We use the sequence $Z_{N}$ in constructed in Subsection \ref{failweak} to
obtain a separated sequence $Z$ with $\left\Vert \mu_{Z}\right\Vert <\infty$
that fails the tree capacity condition (\ref{treecap}), and hence by Theorem C
and (\ref{every}) fails to be onto interpolating sequence for the tree or the disk.

Let $N=b$ in the construction. example $Z_{N}$ above. Recall from that
construction that
\[
Cap\left(  z_{0},Z_{N}\right)  =\gamma_{N}<\gamma_{\infty}<\sqrt{\beta}%
=\frac{1}{\sqrt{N}}%
\]
We claim that for $0<\beta<\frac{1}{36}$, the affine function%
\[
\psi_{\beta}\left(  x\right)  =\frac{\beta}{2}+\left(  1-3\sqrt{\beta}\right)
x
\]
satisfies%
\[
\psi_{\beta}\left(  x\right)  <\varphi_{\beta}\left(  x\right)  <\gamma
_{\infty},\ \ \ \ \ 0\leq x<\gamma_{\infty}.
\]
Indeed, $\varphi_{\beta}-\psi_{\beta}$ is positive at $0$ and has positive
derivative in the interval $\left(  0,\sqrt{\beta}\right)  $ provided
$0<\beta<\frac{1}{36}$. If we now let $\delta_{n}=\psi_{\beta}\left(
\delta_{n-1}\right)  $, $n\geq1$, and $\delta_{0}=0$, then we have by
induction%
\[
\delta_{n}<\gamma_{n}<\gamma_{\infty},\ \ \ \ \ n\geq1.
\]
Indeed, $\delta_{n}=\psi_{\beta}\left(  \delta_{n-1}\right)  <\varphi_{\beta
}\left(  \delta_{n-1}\right)  <\varphi_{\beta}\left(  \gamma_{n-1}\right)
=\gamma_{n}<\gamma_{\infty}$.

Now we compute%
\begin{align*}
\delta_{n}-\delta_{n-1}  &  =\psi_{\beta}\left(  \delta_{n-1}\right)
-\psi_{\beta}\left(  \delta_{n-2}\right)  =\left(  1-3\sqrt{\beta}\right)
\left(  \delta_{n-1}-\delta_{n-2}\right) \\
&  =\left(  1-3\sqrt{\beta}\right)  ^{n-1}\left(  \delta_{1}-\delta
_{0}\right)  =\frac{\beta}{2}\left(  1-3\sqrt{\beta}\right)  ^{n-1},
\end{align*}
and so%
\begin{align*}
\delta_{N}  &  =\sum_{n=1}^{N}\left(  \delta_{n}-\delta_{n-1}\right)
=\frac{\beta}{2}\sum_{n=1}^{N}\left(  1-3\sqrt{\beta}\right)  ^{n-1}\\
&  =\frac{\beta}{2}\frac{1-\left(  1-3\sqrt{\beta}\right)  ^{N}}{1-\left(
1-3\sqrt{\beta}\right)  }=\frac{\sqrt{\beta}}{6}\left\{  1-\left(
1-3\sqrt{\beta}\right)  ^{N}\right\} \\
&  =\frac{1}{6\sqrt{N}}\left\{  1-\left(  1-\frac{3}{\sqrt{N}}\right)
^{N}\right\}  >\frac{1}{12\sqrt{N}},
\end{align*}
for $N$ large since $\left(  1-\frac{3}{\sqrt{N}}\right)  ^{N}\rightarrow0$ as
$N\rightarrow\infty$ by l'Hospital's rule. Altogether we have%
\[
Cap_{z_{0}}\left(  Z_{N}\right)  =\gamma_{N}>\delta_{N}>\frac{1}{12\sqrt{N}}%
\]
for large $N$. If we now take $N=d\left(  z_{0}\right)  ^{\theta}$, $d\left(
z_{0}\right)  $ large and $1\leq\theta<2$, we obtain that the separation
constant $C$ of $Z_{N}$ in (\ref{treesep}) is at least $1$, and that%
\[
\left\Vert \mu_{Z_{N}}\right\Vert =\frac{1}{d\left(  z_{0}\right)  }%
+\sum_{n=1}^{N}\frac{1}{d\left(  w_{n}\right)  }=\frac{1}{N^{\frac{1}{\theta}%
}}+\sum_{n=1}^{N}\frac{1}{d\left(  z_{0}\right)  +n+N}\leq2,
\]
yet
\[
Cap_{z_{0}}\left(  Z_{N}\right)  >\frac{1}{12}\left(  \frac{1}{d\left(
z_{0}\right)  }\right)  ^{\frac{\theta}{2}}>>\frac{1}{d\left(  z_{0}\right)
}.
\]

\subsection{The Simple Condition and Interpolation in the B\"{o}e
Space\label{simcond}}

Suppose that $Z\subset\mathbb{D}$ satisfies the separation condition
(\ref{sep}) and that the associated measure $\mu$ is finite. Here we show that
if the simple condition (\ref{simp}) holds then $R\left(  B_{2,Z}\right)
\subset\ell^{2}\left(  \mu\right)  $, and in the other direction, if $R\left(
B_{2,Z}\right)  \subset\ell^{2}\left(  \mu\right)  $ then a weaker version
(\ref{restsimple}) of condition (\ref{simp}) holds. To see this we fix
$f=\sum_{i=1}^{\infty}a_{i}\varphi_{z_{i}}\in B_{2,Z}$, $z_{j}\in Z$,
and,\qquad\ as Subsection \ref{sp}, we let $\mathcal{Y}$ be the B\"{o}e tree
containing $j$ and%
\[
\mathcal{G}_{j}=\left[  j_{0},j\right]  =\left\{  j_{0},j_{1},...,j_{m-1}%
,j_{m}=j\right\}
\]
be the geodesic $\mathcal{G}_{j}$ in $\mathcal{Y}$ joining $j_{0}$ to $j$.
Then we have%
\[
f\left(  z_{j}\right)  =\sum_{k=0}^{m}a_{j_{k}}\varphi_{z_{j_{k}}}\left(
z_{j}\right)  +\sum_{i\notin\left\{  j_{0},j_{1},...,j_{m}\right\}  }%
a_{i}\varphi_{z_{i}}\left(  z_{j}\right)  .
\]
From H\"{o}lder's inequality, (\ref{atomicBoe}) (which follows from
Proposition \ref{Riesz}) and the third estimate in (\ref{satisfies}) we obtain%
\begin{align*}
\left\vert \sum_{i\notin\left\{  j_{0},j_{1},...,j_{m}\right\}  }a_{i}%
\varphi_{z_{i}}\left(  z_{j}\right)  \right\vert ^{2}  &  \leq C\left\{
\sum_{i}\left\vert a_{i}\right\vert ^{2}\mu\left(  z_{i}\right)  \right\}
\left\{  \sum_{i\notin\left\{  j_{0},j_{1},...,j_{\ell}\right\}  }\left\vert
\varphi_{z_{i}}\left(  z_{j}\right)  \right\vert ^{2}\mu\left(  z_{i}\right)
^{-1}\right\} \\
&  \leq C\left\Vert f\right\Vert _{B_{2,Z}}^{2}\left\{  \sum\nolimits_{i\neq
j_{0}}\left\vert d\left(  z_{i}\right)  ^{-1}\left(  1-\left\vert
z_{i}\right\vert ^{2}\right)  ^{\sigma}\right\vert ^{2}d\left(  z_{i}\right)
\right\} \\
&  \leq C\left\Vert f\right\Vert _{B_{2,Z}}^{2}.
\end{align*}
We also have%
\[
\left\vert \sum\nolimits_{k=0}^{m}a_{j_{k}}\varphi_{z_{j_{k}}}\left(
z_{j}\right)  \right\vert \leq C\sum\nolimits_{k=0}^{m}\left\vert a_{j_{k}%
}\right\vert =CI\left\vert a\right\vert \left(  z_{j}\right)
\]
where $I$ denotes the summation operator on the B\"{o}e tree $\mathcal{Y}$.
Thus we have%
\begin{equation}
\left\Vert Rf\right\Vert _{\ell^{2}\left(  \mu\right)  }\leq C\left\Vert
I\left\vert a\right\vert \right\Vert _{\ell^{2}\left(  \mu\right)
}+C\left\Vert f\right\Vert _{B_{2,Z}}\left\Vert \mu\right\Vert ^{\frac{1}{2}}.
\label{Rmap}%
\end{equation}
By Theorem 3 in \cite{ArRoSa} $I$ is bounded on $\ell^{2}\left(  \mu\right)  $
if and only if%
\begin{equation}
\sum_{\beta\geq\alpha}\frac{I^{\ast}\mu\left(  \beta\right)  ^{2}}{\mu\left(
\beta\right)  }\leq CI^{\ast}\mu\left(  \alpha\right)  ,\ \ \ \ \ \alpha
\in\mathcal{Y}. \label{2weighttree}%
\end{equation}
Now if $\mu$ satisfies the simple condition (\ref{simp}) then $I^{\ast}%
\mu\left(  \beta\right)  \leq C\mu\left(  \beta\right)  $ for $\beta
\in\mathcal{Y}$, and we see that (\ref{2weighttree}) holds. Thus $\left\Vert
I\left\vert a\right\vert \right\Vert _{\ell^{2}\left(  \mu\right)  }\leq
C\left\Vert a\right\Vert _{\ell^{2}\left(  \mu\right)  }\approx\left\Vert
f\right\Vert _{B_{2,Z}}$ and this combined with (\ref{Rmap}) completes the
proof that $R$ maps $B_{2,Z}$ boundedly into $\ell^{2}\left(  \mu\right)  $.

Conversely, if $R$ is bounded from $B_{2,Z}$ to $\ell^{2}\left(  \mu\right)
$, then we have%
\begin{equation}
\sum_{z_{k}\in V_{z_{j}}^{\alpha}}\mu\left(  z_{k}\right)  \leq\left\Vert
R\varphi_{z_{j}}\right\Vert _{\ell^{2}\left(  \mu\right)  }^{2}\leq
C\left\Vert \varphi_{z_{j}}\right\Vert _{B_{2,Z}}^{2}=C\mu\left(
z_{j}\right)  \label{restsimple}%
\end{equation}
for all $z_{j}\in Z$, a weaker version of the simple condition (\ref{simp}).

\section{Converse Results for B\"{o}e Space Interpolating
Sequences\label{neces}}

\subsubsection{Riesz bases of B\"{o}e functions}

The proof that the weak simple condition (\ref{weaksimple}) is necessary for
onto interpolation for the B\"{o}e space $B_{2,Z}$ requires additional tools,
including the fact that the B\"{o}e functions $\left\{  \varphi_{z_{j}%
}\right\}  _{j=1}^{\infty}$ corresponding to a separated sequence $Z=\left\{
z_{j}\right\}  _{j=1}^{\infty}$ in the disk $\mathbb{D}$ form a Riesz basis
for the B\"{o}e space $B_{2,Z}$, at least in the presence of a mild
summability condition on $Z$. It is interesting to note that for a separated
sequence $Z$ in $\mathbb{D}$, the set of Dirichlet reproducing kernels
$\left\{  k_{z_{j}}\right\}  _{j=1}^{\infty}$ form a Riesz basis if and only
if $\mu_{Z}$ is $B_{2}$-Carleson (\cite{Bo}), a condition much stronger than
the mild summability used for the B\"{o}e functions. This points to an
essential advantage of the set of B\"{o}e functions $\left\{  \varphi_{z_{j}%
}\right\}  _{j=1}^{\infty}$ over the set of corresponding normalized
reproducing kernels $\left\{  k_{z_{j}}\left(  z_{j}\right)  ^{-1}k_{z_{j}%
}\right\}  _{j=1}^{\infty}$. The feature of B\"{o}e functions responsible for
this advantage is\ the fact that the supports of the functions $g_{z_{i}}$ are
pairwise disjoint.

\begin{proposition}
\label{Riesz}Let $Z=\left\{  z_{j}\right\}  _{j=1}^{\infty}\subset\mathbb{D}$
satisfy the separation condition (\ref{sep}) and the mild summability
condition $\sum_{j=1}^{\infty}(1-\left\vert z_{j}\right\vert ^{2})^{\sigma
}<\infty$ for all $\sigma>0$. Then there is a finite subset $S$ of $Z\,$such
that $\left\{  \varphi_{z_{j}}\right\}  _{z_{j}\in Z\setminus S}$ is a Riesz
basis for the closed linear span $B_{2,Z}$ of $\left\{  \varphi_{z_{j}%
}\right\}  _{j=1}^{\infty}$ in the Dirichlet space $B_{2}.$
\end{proposition}

\textbf{Proof}: A sequence of B\"{o}e functions $\left\{  \varphi_{z_{j}%
}\right\}  _{j=1}^{\infty}$ is a Riesz basis if
\begin{equation}
C^{-1}\left\Vert \left\{  a_{j}\right\}  _{j=1}^{\infty}\right\Vert _{\ell
^{2}\left(  \mu\right)  }^{2}\leq\left\Vert \sum\nolimits_{j=1}^{\infty}%
a_{j}\varphi_{z_{j}}\right\Vert _{B_{2}}^{2}\leq C\left\Vert \left\{
a_{j}\right\}  _{j=1}^{\infty}\right\Vert _{\ell^{2}\left(  \mu\right)  }^{2}
\label{RB}%
\end{equation}
holds for all sequences $\left\{  a_{j}\right\}  _{j=1}^{\infty}$ with a
positive constant $C$ independent of $\left\{  a_{j}\right\}  _{j=1}^{\infty}%
$. Here $\mu=\sum_{j=1}^{\infty}\left\Vert \varphi_{z_{j}}\right\Vert _{B_{2}%
}^{-2}\delta_{z_{j}}$ and $\mu\left(  z_{j}\right)  =\left\Vert \varphi
_{z_{j}}\right\Vert _{B_{2}}^{-2}\approx d\left(  z_{j}\right)  ^{-1}$. The
inequality on the right follows from (\ref{maincon'}) and the disjoint
supports of the $g_{z_{i}}$ - see the argument use to prove (\ref{bounds})
above - so we concentrate on proving the leftmost inequality in (\ref{RB}) for
an appropriate set of B\"{o}e functions. We begin with%
\begin{align*}
\left\Vert \sum\nolimits_{j=1}^{\infty}a_{j}\varphi_{z_{j}}\right\Vert
_{B_{2}}^{2}  &  =\int_{\mathbb{D}}\left\vert \sum\nolimits_{j=1}^{\infty
}a_{j}\varphi_{z_{j}}^{\prime}\left(  z\right)  \right\vert ^{2}dz\\
&  =\sum_{j,k=1}^{\infty}a_{j}\overline{a_{k}}\int_{\mathbb{D}}\varphi_{z_{j}%
}^{\prime}\left(  z\right)  \overline{\varphi_{z_{k}}^{\prime}\left(
z\right)  }dz\\
&  =\sum_{j=1}^{\infty}\left\vert a_{j}\right\vert ^{2}\mu\left(
z_{j}\right)  +\sum_{j\neq k}a_{j}\overline{a_{k}}\int_{\mathbb{D}}%
\varphi_{z_{j}}^{\prime}\left(  z\right)  \overline{\varphi_{z_{k}}^{\prime
}\left(  z\right)  }dz.
\end{align*}
We now claim that by discarding finitely many points of $Z$, we have
\begin{equation}
\left\vert \sum_{j\neq k}a_{j}\overline{a_{k}}\int_{\mathbb{D}}\varphi_{z_{j}%
}^{\prime}\left(  z\right)  \overline{\varphi_{z_{k}}^{\prime}\left(
z\right)  }dz\right\vert <\frac{1}{2}\sum_{j=1}^{\infty}\left\vert
a_{j}\right\vert ^{2}\mu\left(  z_{j}\right)  . \label{smallmass}%
\end{equation}

Indeed, we will estimate (\ref{smallmass}) using the following derivative
estimates for B\"{o}e functions in the unit disk.

\begin{lemma}
\label{fax}Let $\varphi_{w}\left(  z\right)  $ be as in Lemma
\ref{analyticcon}. Then we have
\[
\left\{
\begin{array}
[c]{lllll}%
\left\vert \varphi_{w}^{\prime}\left(  z\right)  \right\vert  & \leq &
C\left(  1-\left\vert w\right\vert ^{2}\right)  ^{-\alpha}, &  & z\in
V_{w}^{\alpha}\\
\left\vert \varphi_{w}^{\prime}\left(  z\right)  \right\vert  & \leq &
C\left\vert z-w\left\vert w\right\vert ^{-1}\right\vert ^{-1} & \leq\left(
1-\left\vert w\right\vert ^{2}\right)  ^{-\alpha}, & z\in V_{w}^{\rho
}\setminus V_{w}^{\alpha}\\
\left\vert \varphi_{w}^{\prime}\left(  z\right)  \right\vert  & \leq &
C\frac{\left(  1-\left\vert w\right\vert ^{2}\right)  ^{\rho\left(
1+s\right)  }}{\left\vert z-w\left\vert w\right\vert ^{-1}\right\vert ^{2+s}}
& \leq\left(  1-\left\vert w\right\vert ^{2}\right)  ^{-\rho}, & z\notin
V_{w}^{\rho}%
\end{array}
\right.  ,
\]
where $V_{w}^{\beta}=\left\{  z\in\mathbb{D}:\gamma_{w}\left(  z\right)
\geq\beta\right\}  $ and $\gamma_{w}\left(  z\right)  $ is given by
\[
\left\vert z-w\left\vert w\right\vert ^{-1}\right\vert =\left(  1-\left\vert
w\right\vert ^{2}\right)  ^{\gamma_{w}\left(  z\right)  }.
\]

\end{lemma}

%

\proof
This follows readily from the formula
\[
\varphi_{w}\left(  z\right)  =\Gamma_{s}g_{w}\left(  z\right)  =\int
_{\mathbb{D}}\frac{g_{w}\left(  \zeta\right)  \left(  1-\left\vert
\zeta\right\vert ^{2}\right)  ^{s}}{\left(  1-\overline{\zeta}z\right)
^{1+s}}d\zeta,
\]
together with the estimate in \cite[(5.45)]{ArRoSa2},
\[
\left\vert g_{w}\left(  \zeta\right)  \right\vert \leq C\left(  \log\frac
{1}{1-\left\vert w\right\vert ^{2}}\right)  ^{-1}\left\vert \zeta-w\left\vert
w\right\vert ^{-1}\right\vert ^{-1},\;\;\;\;\;\zeta\in\mathbb{D},
\]
and the fact that the support of $g_{w}$ lives in the annular sector
$\mathcal{S}$ centred at $w\left\vert w\right\vert ^{-1}$ given as the
intersection of the annulus
\[
\mathcal{A}=\mathcal{A}_{w}=\left\{  \zeta\in\mathbb{D}:\left(  1-\left\vert
w\right\vert ^{2}\right)  ^{\alpha}\leq\left\vert \zeta-w\left\vert
w\right\vert ^{-1}\right\vert \leq\left(  1-\left\vert w\right\vert
^{2}\right)  ^{\rho}\right\}
\]
and the $45^{\circ}$ angle cone $\mathcal{C}_{w}$ with vertex at $w\left\vert
w\right\vert ^{-1}$. Note that the cone $\mathcal{C}_{w}$ corresponds to the
geodesic in the Bergman tree $\mathcal{T}$ joining the root to the
\textquotedblleft boundary point\textquotedblright\ $w\left\vert w\right\vert
^{-1}$.

\bigskip

The estimate we will prove is, for $j\neq k$,
\begin{align}
\left\vert \left\langle \varphi_{z_{j}},\varphi_{z_{k}}\right\rangle
\right\vert  &  =\left\vert \int_{\mathbb{D}}\varphi_{z_{j}}^{\prime}\left(
z\right)  \overline{\varphi_{z_{k}}^{\prime}\left(  z\right)  }dz\right\vert
\label{want}\\
&  \leq C\left(  1-\left\vert z_{j}\right\vert ^{2}\right)  ^{\sigma}\left(
1-\left\vert z_{k}\right\vert ^{2}\right)  ^{\sigma}\mu\left(  z_{j}\right)
\mu\left(  z_{k}\right)  ,\nonumber
\end{align}
for some $\sigma>0$. We may assume that $1-\left\vert z_{j}\right\vert
^{2}\leq1-\left\vert z_{k}\right\vert ^{2}$ and write%
\begin{align*}
\left\vert \int_{\mathbb{D}}\varphi_{z_{j}}^{\prime}\left(  z\right)
\overline{\varphi_{z_{k}}^{\prime}\left(  z\right)  }dz\right\vert  &
\leq\int_{\mathbb{D}}\left\vert \varphi_{z_{j}}^{\prime}\left(  z\right)
\right\vert \left\vert \varphi_{z_{k}}^{\prime}\left(  z\right)  \right\vert
dz\\
&  =\left\{  \int_{V_{z_{j}}}+\int_{\mathbb{D}\smallsetminus V_{z_{j}}%
}+\right\}  \left\vert \varphi_{z_{j}}^{\prime}\left(  z\right)  \right\vert
\left\vert \varphi_{z_{k}}^{\prime}\left(  z\right)  \right\vert dz\\
&  =I+II.
\end{align*}
To estimate $II$ we use Lemma \ref{fax} to obtain%
\begin{align*}
II  &  \leq C\int_{\mathbb{D}\smallsetminus V_{z_{j}}}\frac{\left(
1-\left\vert z_{j}\right\vert ^{2}\right)  ^{\rho\left(  1+s\right)  }%
}{\left\vert 1-\overline{z}z_{j}\right\vert ^{2+s}}\left(  1-\left\vert
z_{k}\right\vert ^{2}\right)  ^{-\alpha}dz\\
&  =C\left(  1-\left\vert z_{k}\right\vert ^{2}\right)  ^{-\alpha}\left(
1-\left\vert z_{j}\right\vert ^{2}\right)  ^{\rho\left(  1+s\right)  }%
\int_{\mathbb{D}\smallsetminus V_{z_{j}}}\frac{dz}{\left\vert 1-\overline
{z}z_{j}\right\vert ^{2+s}}\\
&  \leq C\left(  1-\left\vert z_{k}\right\vert ^{2}\right)  ^{-\alpha}\left(
1-\left\vert z_{j}\right\vert ^{2}\right)  ^{\rho\left(  1+s\right)  }\left(
1-\left\vert z_{j}\right\vert ^{2}\right)  ^{-\beta s}\\
&  \leq C\left(  1-\left\vert z_{j}\right\vert ^{2}\right)  ^{\rho\left(
1+s\right)  -\beta s-\alpha}.
\end{align*}
Using (\ref{parameters}) we see that the exponent $\rho\left(  1+s\right)
-\beta s-\alpha$ is positive, and using $1-\left\vert z_{j}\right\vert
^{2}\leq1-\left\vert z_{k}\right\vert ^{2}$ we easily obtain (\ref{want}).

To estimate $I$ we consider two cases. In the case that $V_{z_{j}}\cap
V_{z_{k}}\neq\phi$, we have from Lemma \ref{fax} and the estimate $\left\vert
V_{z_{j}}\right\vert \leq C\left(  1-\left\vert z_{j}\right\vert ^{2}\right)
^{2\beta}$ that%
\begin{align*}
I  &  \leq C\sup_{\mathbb{D}}\left\vert \varphi_{z_{j}}^{\prime}\right\vert
\sup_{\mathbb{D}}\left\vert \varphi_{z_{k}}^{\prime}\right\vert \left\vert
V_{z_{j}}\right\vert \\
&  \leq C\left(  1-\left\vert z_{j}\right\vert ^{2}\right)  ^{-\alpha}\left(
1-\left\vert z_{k}\right\vert ^{2}\right)  ^{-\alpha}\left(  1-\left\vert
z_{j}\right\vert ^{2}\right)  ^{2\beta}\\
&  \leq C\left(  1-\left\vert z_{j}\right\vert ^{2}\right)  ^{2\beta-\alpha
}\left(  1-\left\vert z_{k}\right\vert ^{2}\right)  ^{-\alpha},
\end{align*}
and now Lemma \ref{Mars} yields%
\[
I\leq C\left(  1-\left\vert z_{j}\right\vert ^{2}\right)  ^{2\beta-\alpha
}\left(  1-\left\vert z_{j}\right\vert ^{2}\right)  ^{-\frac{\alpha}{\eta}%
}=C\left(  1-\left\vert z_{j}\right\vert ^{2}\right)  ^{2\beta-\alpha
-\frac{\alpha}{\eta}}.
\]
Now using (\ref{parameters}) we see that the exponent $2\beta-\alpha
-\frac{\alpha}{\eta}$ is positive and we again obtain (\ref{want}). On the
other hand, if $V_{z_{j}}\cap V_{z_{k}}=\phi$, then we use%
\begin{align*}
I  &  \leq C\sup_{\mathbb{D}}\left\vert \varphi_{z_{j}}^{\prime}\right\vert
\sup_{V_{z_{j}}}\left\vert \varphi_{z_{k}}^{\prime}\right\vert \left\vert
V_{z_{j}}\right\vert \\
&  \leq C\left(  1-\left\vert z_{j}\right\vert ^{2}\right)  ^{-\alpha}%
\frac{\left(  1-\left\vert z_{k}\right\vert ^{2}\right)  ^{\rho\left(
1+s\right)  }}{\left\vert 1-\overline{z_{j}}\cdot z_{k}\right\vert ^{2+s}%
}\left(  1-\left\vert z_{j}\right\vert ^{2}\right)  ^{2\beta}\\
&  \leq C\left(  1-\left\vert z_{j}\right\vert ^{2}\right)  ^{2\beta-\alpha
}\left(  1-\left\vert z_{k}\right\vert ^{2}\right)  ^{\rho\left(  1+s\right)
-\beta\left(  2+s\right)  }\\
&  \leq C\left(  1-\left\vert z_{j}\right\vert ^{2}\right)  ^{\varepsilon
}\left(  1-\left\vert z_{k}\right\vert ^{2}\right)  ^{\rho\left(  1+s\right)
-\beta\left(  2+s\right)  +2\beta-\alpha-\varepsilon},
\end{align*}
upon using $\left(  1-\left\vert z_{j}\right\vert ^{2}\right)  ^{2\beta
-\alpha-\varepsilon}\leq\left(  1-\left\vert z_{k}\right\vert ^{2}\right)
^{2\beta-\alpha-\varepsilon}$ in the last line. Now choosing%
\[
s>\frac{\alpha-\rho+\varepsilon}{\rho-\beta},
\]
the exponent $\rho\left(  1+s\right)  -\beta\left(  2+s\right)  +2\beta
-\alpha-\varepsilon$ is positive, and once more we obtain (\ref{want}).

Now we can estimate the left side of (\ref{smallmass}) by (\ref{want}) and
(\ref{sigmasum}) to obtain%
\begin{align*}
\left\vert \sum\nolimits_{j\neq k}a_{j}\overline{a_{k}}\int_{\mathbb{D}%
}\varphi_{z_{j}}^{\prime}\left(  z\right)  \overline{\varphi_{z_{k}}^{\prime
}\left(  z\right)  }dz\right\vert  &  \leq C\sum_{j\neq k}\left\vert
a_{j}a_{k}\right\vert \left(  1-\left\vert z_{j}\right\vert ^{2}\right)
^{\sigma}\left(  1-\left\vert z_{k}\right\vert ^{2}\right)  ^{\sigma}%
\mu\left(  z_{j}\right)  \mu\left(  z_{k}\right) \\
&  \leq C\left\{  \sum\nolimits_{k}\left(  1-\left\vert z_{k}\right\vert
^{2}\right)  ^{\sigma}\mu\left(  z_{k}\right)  \right\}  \sum_{j=1}^{\infty
}\left\vert a_{j}\right\vert ^{2}\mu\left(  z_{j}\right) \\
&  <\frac{1}{2}\sum_{j=1}^{\infty}\left\vert a_{j}\right\vert ^{2}\mu\left(
z_{j}\right)
\end{align*}
if $\sum_{k}\left(  1-\left\vert z_{k}\right\vert ^{2}\right)  ^{\sigma}%
\mu\left(  z_{k}\right)  $ is sufficiently small, which can be achieved by
discarding a sufficiently large finite subset $F$ from $Z$. This shows that
$\left\{  \varphi_{z_{j}}\right\}  _{z_{j}\in Z\setminus F}$ is a Riesz basis.
However, if $w\in F$ is \emph{not} in the closed linear span of the Riesz
basis $\left\{  \varphi_{z_{j}}\right\}  _{z_{j}\in Z\smallsetminus F}$, then
it is immediate that $\left\{  \varphi_{z_{j}}\right\}  _{z_{j}\in
Z\smallsetminus F}\cup\left\{  \varphi_{w}\right\}  $ is also a Riesz basis.
We can continue adding B\"{o}e functions $\varphi_{w}$ with $w\in G\subset F $
so that $\left\{  \varphi_{z_{j}}\right\}  _{z_{j}\in Z\setminus F}%
\cup\left\{  \varphi_{w}\right\}  _{w\in G}$ is a Riesz basis, and such that
\emph{all} of the remaining Boe functions $\varphi_{w}$ with $w\in
F\smallsetminus G$ lie in the closed linear span of the Riesz basis $\left\{
\varphi_{z_{j}}\right\}  _{z_{j}\in Z\setminus F}\cup\left\{  \varphi
_{w}\right\}  _{w\in G}$. This completes the proof of Proposition \ref{Riesz}
with $S=F\smallsetminus G$.

\subsubsection{Completion of the Proofs of Theorems \ref{Boeonto} and
\ref{onto}}

Now we consider the necessity of the two conditions (\ref{sep}), or
equivalently (\ref{sep'}), and (\ref{weaksimple}) in Theorem \ref{Boeonto}. We
noted when we introduced (\ref{sep'}) that it is necessary. To see that
(\ref{weaksimple}) is necessary, we note that by Proposition \ref{Riesz} above
(the summability hypothesis there is a consequence of $\left\Vert
\mu\right\Vert <\infty$), we can remove a finite subset $S$ from $Z$ so that
$B_{2,Z\smallsetminus S}\left(  \mathbb{D}\right)  =B_{2,Z}$ and $\left\{
\varphi_{z_{j}}\right\}  _{z_{j}\in Z\smallsetminus S}$ is a Riesz basis. We
can obviously add finitely many points to a sequence satisfying the weak
simple condition and obtain a new sequence satisfying the weak simple
condition. Thus we may assume that (\ref{atomicBoe}) holds for $Z$.

Now let $e_{j}$ be the function on $Z$ that is $1$ at $z_{j}$ and vanishes on
the rest of $Z$. Denote the collection of all children of $z_{j}$ in the
forest structure $\mathcal{F}$ by $\mathcal{C}\left(  z_{j}\right)  $, and let
$\mu=\mu_{Z}$. We now claim that for $j$ sufficiently large,
\begin{equation}
\mathcal{S}e_{j}=\varphi_{z_{j}}-\sum_{z_{i}\in\mathcal{C}\left(
z_{j}\right)  }\varphi_{z_{j}}\left(  z_{i}\right)  \varphi_{z_{i}}+f_{j},
\label{claimS}%
\end{equation}
where $f_{j}\in B_{2,Z}$ has the form%
\[
f_{j}=\sum_{i=1}^{\infty}a_{i}\varphi_{z_{i}}%
\]
with $\left\{  a_{i}\right\}  _{i=1}^{\infty}\in\ell^{2}\left(  \mu\right)  $
and%
\[
\left\vert a_{j}\right\vert <\frac{1}{2};\text{ }\left\vert a_{i}\right\vert
<\frac{1}{2}\ \text{if }z_{i}\in\mathcal{C}\left(  z_{j}\right)  .
\]
Indeed, by (\ref{atomicBoe}) we have%
\begin{equation}
\mathcal{S}e_{j}=\sum_{i=1}^{\infty}b_{i}\varphi_{z_{i}} \label{have}%
\end{equation}
with $\left\{  b_{i}\right\}  _{i=1}^{\infty}\in\ell^{2}\left(  \mu\right)  $
and $\left\Vert \left\{  b_{i}\right\}  _{i=1}^{\infty}\right\Vert _{\ell
^{2}\left(  \mu\right)  }^{2}\approx\mu\left(  z_{j}\right)  $.

Now let $\mathcal{Y}$ be the B\"{o}e tree containing $j$ and%
\[
\mathcal{G}_{j}=\left[  j_{0},j\right]  =\left\{  j_{0},j_{1},...,j_{m-1}%
,j_{m}=j\right\}
\]
be the geodesic $\mathcal{G}_{j}$ in $\mathcal{Y}$ joining $j_{0}$ to $j$. If
we evaluate both sides of (\ref{have}) at $z_{j_{\ell}}$ where $0\leq\ell<m$,
we have%
\begin{equation}
0=\mathcal{S}e_{j}\left(  z_{j_{\ell}}\right)  =\sum_{k=0}^{\ell}b_{j_{k}%
}\varphi_{z_{j_{k}}}\left(  z_{j_{\ell}}\right)  +\sum_{i\notin\left\{
j_{0},j_{1},...,j_{\ell}\right\}  }b_{i}\varphi_{z_{i}}\left(  z_{j_{\ell}%
}\right)  . \label{evaluate}%
\end{equation}
Subtracting the cases $\ell$ and $\ell+1$ in (\ref{evaluate}) we obtain%
\begin{align*}
0  &  =\mathcal{S}e_{j}\left(  z_{j_{\ell+1}}\right)  -\mathcal{S}e_{j}\left(
z_{j_{\ell}}\right) \\
&  =\sum_{k=0}^{\ell-1}b_{j_{k}}\left[  \varphi_{z_{j_{k}}}\left(
z_{j_{\ell+1}}\right)  -\varphi_{z_{j_{k}}}\left(  z_{j_{\ell}}\right)
\right]  +b_{j_{\ell}}\left(  \varphi_{z_{j_{\ell}}}\left(  z_{j_{\ell+1}%
}\right)  -1\right) \\
&  +b_{j_{\ell+1}}+\sum_{i\notin\left\{  j_{0},j_{1},...,j_{\ell+1}\right\}
}b_{i}\varphi_{z_{i}}\left(  z_{j_{\ell+1}}\right)  -\sum_{i\notin\left\{
j_{0},j_{1},...,j_{\ell}\right\}  }b_{i}\varphi_{z_{i}}\left(  z_{j_{\ell}%
}\right)  .
\end{align*}
From H\"{o}lder's inequality and the third estimate in (\ref{satisfies}) we
obtain%
\begin{align}
\left\vert \sum_{i\notin\left\{  j_{0},j_{1},...,j_{\ell}\right\}  }%
b_{i}\varphi_{z_{i}}\left(  z_{j_{\ell}}\right)  \right\vert  &  \leq\left\{
\sum_{i}\left\vert b_{i}\right\vert ^{2}\mu\left(  z_{i}\right)  \right\}
^{\frac{1}{2}}\left\{  \sum_{i\notin\left\{  j_{0},j_{1},...,j_{\ell}\right\}
}\left\vert \varphi_{z_{i}}\left(  z_{j_{0}}\right)  \right\vert ^{2}%
\mu\left(  z_{i}\right)  ^{-1}\right\}  ^{\frac{1}{2}}\label{boundsum}\\
&  \leq C\mu\left(  z_{j}\right)  ^{\frac{1}{2}}\left\{  \sum\nolimits_{i\neq
j_{0}}\left\vert d\left(  z_{i}\right)  ^{-1}\left(  1-\left\vert
z_{i}\right\vert ^{2}\right)  ^{\sigma}\right\vert ^{2}d\left(  z_{i}\right)
\right\}  ^{\frac{1}{2}}\nonumber\\
&  \leq C_{0}\mu\left(  z_{j}\right)  ^{\frac{1}{2}},\nonumber
\end{align}
where the final term in braces is bounded by hypothesis. We also have from
(\ref{satisfies})
\[
\left\vert \varphi_{z_{j_{\ell}}}\left(  z_{j_{\ell+1}}\right)  -1\right\vert
\leq\left(  1+C\mu\left(  z_{j_{\ell}}\right)  \right)
\]
and%
\begin{align*}
&  \left\vert \sum\nolimits_{k=0}^{\ell-1}b_{j_{k}}\left[  \varphi_{z_{j_{k}}%
}\left(  z_{j_{\ell+1}}\right)  -\varphi_{z_{j_{k}}}\left(  z_{j_{\ell}%
}\right)  \right]  \right\vert \\
&  \leq C\sum_{k=0}^{\ell-1}\left\vert b_{j_{k}}\right\vert \mu\left(
z_{j_{k}}\right)  \leq C\left\Vert \mu\right\Vert ^{\frac{1}{2}}\left\{
\sum\nolimits_{k=0}^{\ell-1}\left\vert b_{j_{k}}\right\vert ^{2}\mu\left(
z_{j_{k}}\right)  \right\}  ^{\frac{1}{2}}.
\end{align*}
Altogether then we have%
\begin{align*}
\left\vert b_{j_{\ell+1}}\right\vert  &  \leq\left\vert b_{j_{\ell}%
}\right\vert \left(  1+C\mu\left(  z_{j_{\ell}}\right)  \right)  +C\left\Vert
\mu\right\Vert ^{\frac{1}{2}}\left\{  \sum\nolimits_{k=0}^{\ell-1}\left\vert
b_{j_{k}}\right\vert ^{2}\mu\left(  z_{j_{k}}\right)  \right\}  ^{\frac{1}{2}%
}+2C_{0}\mu\left(  z_{j}\right)  ^{\frac{1}{2}}\\
&  \leq\left\vert b_{j_{\ell}}\right\vert \left(  1+C\mu\left(  z_{j_{\ell}%
}\right)  \right)  +C_{1}\mu\left(  z_{j}\right)  ^{\frac{1}{2}}.
\end{align*}
Now the case $\ell=0$ of (\ref{evaluate}) together with (\ref{boundsum})
yields%
\[
\left\vert b_{j_{0}}\right\vert =\left\vert \sum_{i\notin\left\{  j_{0}%
,j_{1},...,j_{\ell}\right\}  }b_{i}\varphi_{z_{i}}\left(  z_{j_{\ell}}\right)
\right\vert \leq C_{0}\mu\left(  z_{j}\right)  ^{\frac{1}{2}},
\]
and now by induction on $\ell$ we obtain that for $0\leq\ell\leq m-1$,
$\left\vert b_{j_{\ell}}\right\vert $ is dominated by%
\[
C_{0}\mu\left(  z_{j}\right)  ^{\frac{1}{2}}\left\{
{\displaystyle\prod\limits_{k=0}^{\ell-1}}
\left(  1+Cd\left(  z_{j_{k}}\right)  ^{-1}\right)  +%
{\displaystyle\prod\limits_{k=1}^{\ell-1}}
\left(  1+Cd\left(  z_{j_{k}}\right)  ^{-1}\right)  +...+\left(  1+Cd\left(
z_{j_{\ell-1}}\right)  ^{-1}\right)  \right\}  .
\]
In particular,
\begin{equation}
\left\vert b_{j_{\ell}}\right\vert \leq C_{0}\mu\left(  z_{j}\right)
^{\frac{1}{2}}\ell\exp\left(  C\sum\nolimits_{k=0}^{\ell-1}d\left(  z_{j_{k}%
}\right)  ^{-1}\right)  \label{bjl}%
\end{equation}
for $0\leq\ell\leq m-1$.

Now evaluate both sides of (\ref{have}) at $z_{j}=z_{j_{m}}$ to obtain%
\[
1=b_{j}+\sum_{k=0}^{m-1}b_{j_{k}}\varphi_{z_{j_{k}}}\left(  z_{j}\right)
+\sum_{i\notin\left\{  j_{0},j_{1},...,j_{m}\right\}  }b_{i}\varphi_{z_{i}%
}\left(  z_{j}\right)  ,
\]
which by the argument above yields%
\[
\left\vert b_{j}-1\right\vert \leq C_{0}\mu\left(  z_{j}\right)  ^{\frac{1}%
{2}}m\exp\left(  C\sum\nolimits_{k=0}^{m}d\left(  z_{j_{k}}\right)
^{-1}\right)  .
\]
Similarly, for $z_{i}\in\mathcal{C}\left(  z_{j}\right)  $ we obtain%
\[
\left\vert b_{i}-b_{j}\varphi_{z_{j}}\left(  z_{i}\right)  \right\vert \leq
C_{0}\mu\left(  z_{j}\right)  ^{\frac{1}{2}}\left(  m+1\right)  \exp\left(
Cd\left(  z_{i}\right)  ^{-1}+C\sum\nolimits_{k=0}^{m}d\left(  z_{j_{k}%
}\right)  ^{-1}\right)  .
\]
Now the separation condition (\ref{sep'}) yields $d\left(  z_{j_{k}}\right)
\geq\left(  1+C\right)  d\left(  z_{j_{k-1}}\right)  $ for $1\leq k\leq m$ and
it follows that
\begin{equation}
\sum_{k=0}^{m}d\left(  z_{j_{k}}\right)  ^{-1}\leq C \label{sumC}%
\end{equation}
independent of $\dot{j}$. Thus we see that%
\[
\left\vert b_{i}-b_{j}\varphi_{z_{j}}\left(  z_{i}\right)  \right\vert \leq
C\left(  m+1\right)  \mu\left(  z_{j}\right)  ^{\frac{1}{2}},\ \ \ \ \ z_{i}%
\in\mathcal{C}\left(  z_{j}\right)  ,
\]
with a constant $C$ independent of $\dot{j}$. If we take $j_{0}$ large enough,
then since $d\left(  z_{j}\right)  =d\left(  z_{j_{m}}\right)  \geq\left(
1+C\right)  ^{m}d\left(  z_{j_{0}}\right)  $, we have
\[
\left\vert b_{j}-1\right\vert \leq Cm\mu\left(  z_{j}\right)  ^{\frac{1}{2}%
}=Cmd\left(  z_{j}\right)  ^{-\frac{1}{2}}\leq C\frac{m}{\left(  1+C\right)
^{\frac{m}{2}}}d\left(  z_{j_{0}}\right)  ^{-\frac{1}{2}}<\frac{1}{2}.
\]
It follows that
\[
\left\vert b_{i}-\varphi_{z_{j}}\left(  z_{i}\right)  \right\vert
\leq\left\vert b_{i}-b_{j}\varphi_{z_{j}}\left(  z_{i}\right)  \right\vert
+\left\vert b_{j}-1\right\vert \left\vert \varphi_{z_{j}}\left(  z_{i}\right)
\right\vert <\frac{1}{2},\ \ \ \ \ z_{i}\in\mathcal{C}\left(  z_{j}\right)  ,
\]
which proves (\ref{claimS}).

By (\ref{atomicBoe}) we then have using (\ref{claimS}) and the fact that
$\varphi_{z_{j}}\left(  z_{i}\right)  =1$ for $z_{i}\in\mathcal{C}\left(
z_{j}\right)  \cap V_{z_{j}}^{\alpha}$:
\[
\left\Vert \mathcal{S}e_{j}\right\Vert _{B_{2,Z}}\approx\left\{  \sum
_{i}\left\vert b_{i}\right\vert ^{2}\mu\left(  z_{i}\right)  \right\}
^{\frac{1}{2}}\geq\frac{1}{2}\left\{  \sum_{z_{i}\in\mathcal{C}\left(
z_{j}\right)  \cap V_{z_{j}}^{\alpha}}\mu\left(  z_{i}\right)  \right\}
^{\frac{1}{2}}.
\]
It follows that%
\[
\mu\left(  z_{j}\right)  =\left\Vert e_{j}\right\Vert _{\ell^{2}\left(
\mu\right)  }^{2}\geq C^{2}\left\Vert \mathcal{S}e_{j}\right\Vert _{B_{2,Z}%
}^{2}\geq C^{\prime}\sum_{z_{i}\in\mathcal{C}\left(  z_{j}\right)  \cap
V_{z_{j}}^{\alpha}}\mu\left(  z_{i}\right)  ,
\]
which yields (\ref{weaksimple}) for $\alpha=z_{j}\in Z$ with $j$ large, and
hence for all $j$ with a worse constant.

\bigskip

Now we suppose that $\alpha\in\mathcal{T}\setminus Z$. We claim that with
either $z_{0}=\alpha$ or $z_{0}=A^{M}\alpha$, where $M=\left[  \frac{C}%
{10}d\left(  \alpha\right)  \right]  $ and $C$ is as in (\ref{sep'}), the set
$Z^{\prime}=Z\cup\left\{  z_{0}\right\}  $ is separated with separation
constant in (\ref{sep'}) at least $\ C/100$. Indeed, if $Z\cup\left\{
\alpha\right\}  $ fails to satisfy (\ref{sep'}) with separation constant
$C/100$, then there is some $w$ in $Z$ such that%
\[
\beta\left(  \alpha,w\right)  <\frac{C}{50}\left(  1+\beta\left(  o,w\right)
\right)  .
\]
From this we obtain that%
\begin{align*}
\beta\left(  A^{M}\alpha,w\right)   &  \geq\beta\left(  A^{M}\alpha
,\alpha\right)  -\beta\left(  \alpha,w\right) \\
&  >\frac{C}{10}\left(  1+\beta\left(  o,w\right)  \right)  -\frac{C}%
{50}\left(  1+\beta\left(  o,w\right)  \right) \\
&  >\frac{C}{20}\left(  1+\beta\left(  o,w\right)  \right)  ,
\end{align*}
and then for any $z\in Z\setminus\left\{  w\right\}  $,%
\begin{align*}
\beta\left(  A^{M}\alpha,z\right)   &  \geq\beta\left(  w,z\right)
-\beta\left(  A^{M}\alpha,w\right) \\
&  >\beta\left(  w,z\right)  -\left\{  \beta\left(  A^{M}\alpha,\alpha\right)
+\beta\left(  \alpha,w\right)  \right\} \\
&  >C\left(  1+\beta\left(  o,w\right)  \right)  -\left\{  \frac{C}{10}\left(
1+\beta\left(  o,w\right)  \right)  +\frac{C}{50}\left(  1+\beta\left(
o,w\right)  \right)  \right\} \\
&  >\frac{C}{2}\left(  1+\beta\left(  o,w\right)  \right)  ,
\end{align*}
which shows that $Z\cup\left\{  A^{M}\alpha\right\}  $ satisfies (\ref{sep'})
with separation constant C/2. Now we associate a B\"{o}e function
$\varphi_{z_{0}}$ with $z_{0}$, but take the parameters $\beta,\beta_{1}%
,\rho,\alpha$ so close to $1$ for this additional function $\varphi_{z_{0}}$
that the extended set of B\"{o}e functions $\left\{  \varphi_{z}\right\}
_{z\in Z^{\prime}}=\left\{  \varphi_{z}\right\}  _{z\in Z}\cup\left\{
\varphi_{z_{0}}\right\}  $ satisfy the property that the supports of the
associated functions $g_{z}$ are pairwise disjoint for $z\in Z^{\prime}$.

Now we define a bounded linear operator $S^{\prime}$ from $\ell^{2}\left(
\mu_{Z^{\prime}}\right)  $ into $B_{2,Z^{\prime}}\left(  \mathcal{T}\right)  $
by%
\[
S^{\prime}\left[  \xi^{\prime}\right]  =S\xi+\left(  \xi_{0}-S\xi\left(
z_{0}\right)  \right)  \left\{  \varphi_{z_{0}}-S\left[  \varphi_{z_{0}}%
\mid_{Z}\right]  \right\}  ,
\]
where $\xi^{\prime}=\left(  \xi_{0},\xi\right)  =\left(  \xi_{0},\xi
_{1},...\right)  $. For $j\geq1$ we have%
\begin{align*}
S^{\prime}\left[  \xi^{\prime}\right]  \left(  z_{j}\right)   &  =\xi
_{j}+\left(  \xi_{0}-S\xi\left(  z_{0}\right)  \right)  \left\{
\varphi_{z_{0}}\left(  z_{j}\right)  -S\left[  \varphi_{z_{0}}\mid_{Z}\right]
\left(  z_{j}\right)  \right\} \\
&  =\xi_{j}+\left(  \xi_{0}-S\xi\left(  z_{0}\right)  \right)  \left\{
0\right\}  =\xi_{j},
\end{align*}
and for $j=0$,%
\begin{align*}
S^{\prime}\left[  \xi^{\prime}\right]  \left(  z_{0}\right)   &  =S\xi\left(
z_{0}\right)  +\left(  \xi_{0}-S\xi\left(  z_{0}\right)  \right)  \left\{
1-S\left[  \varphi_{z_{0}}\mid_{Z}\right]  \left(  z_{0}\right)  \right\} \\
&  =\xi_{0}-S\left[  \varphi_{z_{0}}\mid_{Z}\right]  \left(  z_{0}\right)
\left(  \xi_{0}-S\xi\left(  z_{0}\right)  \right)  .
\end{align*}
Now $S\left[  \varphi_{z_{0}}\mid_{Z}\right]  \left(  z_{0}\right)  $ is small
by the argument used to prove (\ref{claimS}) above, and in fact (\ref{bjl})
and (\ref{sumC}) of that argument yield%
\[
\left\vert S\left[  \varphi_{z_{0}}\mid_{Z}\right]  \left(  z_{0}\right)
\right\vert \leq C\mu\left(  z_{0}\right)  ^{\frac{1}{2}}.
\]
At this point we may assume that $C\mu\left(  z_{0}\right)  ^{\frac{1}{2}%
}<\varepsilon$ since there are only finitely many (depending on $\varepsilon
>0$) points $\alpha$ in the tree $\mathcal{T}$ having such a point $z_{0}$
that fails this condition. Thus $S^{\prime}$ is an \emph{approximate} bounded
right inverse to the restriction map $\mathcal{U}$, and in fact,%
\[
\mathcal{U}S^{\prime}\xi^{\prime}-\xi^{\prime}=S\left[  \varphi_{z_{0}}%
\mid_{Z}\right]  \left(  z_{0}\right)  \left(  \xi_{0}-S\xi\left(
z_{0}\right)  \right)  e_{z_{0}},
\]
so that%
\[
\left\Vert \mathcal{U}S^{\prime}\xi^{\prime}-\xi^{\prime}\right\Vert
_{\ell^{2}\left(  \mu\right)  }\leq\varepsilon C\left\Vert \xi^{\prime
}\right\Vert _{\ell^{2}\left(  \mu\right)  }<\frac{1}{2}\left\Vert \xi
^{\prime}\right\Vert _{\ell^{2}\left(  \mu\right)  }%
\]
if $\varepsilon>0$ is small enough. Then $\mathcal{U}S^{\prime}$ is invertible
on $\ell^{2}\left(  \mu\right)  $, and so the operator $S^{\prime\prime
}=S^{\prime}\left(  \mathcal{U}S^{\prime}\right)  ^{-1}$ is an exact bounded
right inverse to the restriction map $\mathcal{U}$ since $\mathcal{U}%
S^{\prime\prime}=\mathcal{U}S^{\prime}\left(  \mathcal{U}S^{\prime}\right)
^{-1}=\mathbb{I}_{\ell^{2}\left(  \mu\right)  }$. Then the result proved in
the previous paragraph shows that the weak simple condition (\ref{weaksimple})
holds at $z_{0}$ with a controlled constant, and thus also at $\alpha$ with a
controlled constant. This completes the proof of Theorem \ref{Boeonto}.

It remains to show the necessity of (\ref{weaksimple}) in the context of
Theorem \ref{onto}. For that situation, when $Z$ is onto interpolating for the
B\"{o}e space $B_{2,Z}$, we note that a subtree of a dyadic tree has branching
number at most $2$, and it follows easily from the separation condition that%
\[
\sum_{j=1}^{\infty}\left(  1-\left\vert z_{j}\right\vert ^{2}\right)
^{\sigma}<\infty
\]
for all $\sigma>0$. Thus Proposition \ref{Riesz} can be applied together with
the argument used above to prove necessity of (\ref{weaksimple}) in the case
$\left\Vert \mu_{Z}\right\Vert <\infty$.

\end{document}